# Erasure Techniques in MRD codes

W. B. Vasantha Kandasamy
Florentin Smarandache
R. Sujatha
R. S. Raja Durai

# Erasure Techniques in MRD codes

**W. B. Vasantha Kandasamy**
**Florentin Smarandache**
**R. Sujatha**
**R. S. Raja Durai**

# CONTENTS

# PREFACE

In this book the authors study the erasure techniques in concatenated Maximum Rank Distance (MRD) codes. The authors for the first time in this book introduce the new notion of concatenation of MRD codes with binary codes, where we take the outer code as the RD code and the binary code as the inner code. The concatenated code consists of the codewords of the outer code expressed in terms of the alphabets of the inner code. These new class of codes are defined as CRM codes. This concatenation techniques helps one to construct any CRM code of desired minimum distance which is not enjoyed by any other class of codes.

Also concatenation of several binary codes are introduced using the newly defined notion of special blanks. These codes can be used in bulk transmission of a message into several channels and the completed work is again consolidated and received.

Finally the notion of integer rank distance code is introduced. This book is organized into six chapters. The first chapter introduces the basic algebraic structures essential to make this book a self contained one. Algebraic linear codes and their basic properties are discussed in chapter two. In chapter three the authors study the basic properties of erasure decoding in maximum rank distance codes.

Some decoding techniques about MRD codes are described and discussed in chapter four of this book. Rank distance codes with complementary duals and MRD codes with complementary duals are introduced and their applications are discussed. Chapter five introduces the notion of integer rank distance codes. The final chapter introduces some concatenation techniques.

We thank Dr. K.Kandasamy for proof reading and being extremely supportive.

W.B.VASANTHA KANDASAMY
FLORENTIN SMARANDACHE
R. SUJATHA
R. S. RAJA DURAI

## Chapter One

# BASIC CONCEPTS

In this chapter we give the basic concepts to make this book a self contained one. Basically we need the notion of vector spaces over finite fields and the notion of irreducible polynomials over finite fields. Also the notion of cosets of groups for error correction is needed. We will briefly recall only those facts, which is essential for a beginner to under stand algebraic coding theory.

**DEFINITION 1.1:** *Let G be a non empty set '*' a closed associative binary operation defined on G such that;*

*(i)* *there exists a unique element e in G with* $g * e = e * g = g$ *for all* $g \in G$.

*(ii)* *For every* $g \in G$ *there exists a unique* $g'$ *in G such that* $g * g' = g' * g = e$; $g'$ *is called the inverse of g in G under the binary operation *.*

*Then (G, *) is defined as a group.*

*If in G for every pair a, b ∈ G; a*b = b*a then we say G is an abelian group.*

If the number of distinct elements in G is finite then G is said to be a finite group otherwise G is said to be an infinite group.

***Example 1.1:*** $Z_2 = \{0, 1\}$ is an additive abelian group given by the following table:

| + | 0 | 1 |
|---|---|---|
| 0 | 0 | 1 |
| 1 | 1 | 0 |

$Z_2$ under multiplication is not a group.

***Example 1.2:*** $(Z_3, +)$ is an abelian group of order 3 given by the following table:

| + | 0 | 1 | 2 |
|---|---|---|---|
| 0 | 0 | 1 | 2 |
| 1 | 1 | 2 | 0 |
| 2 | 2 | 0 | 1 |

***Example 1.3:*** $(Z_n, +)$, $(1 < n < \infty)$ is an abelian group of order n.

***Example 1.4:*** $G = Z_2 \times Z_2 \times Z_2 = \{(a, b, c) / a, b, c \in Z_2\}$ is a group under addition modulo 2. G is of order 8.

***Example 1.5:*** Let $G = Z_2 \times Z_2 \times \ldots \times Z_2$, n times with

$$G = \{(x_1, x_2, \ldots, x_n) / x_i \in Z_2; 1 \leq i \leq n\}$$

is an abelian group under addition modulo 2 of order $2^n$.

***Example 1.6:*** $H = Z_p \times Z_p \times \ldots \times Z_p$ - t times is an abelian group under addition modulo p of order $p^t$.

Now we will just define the notion of a subgroup.

**DEFINITION 1.2:** *Let (G, *) be a group. H a proper subset of G. If (H, *) is a group then, we call H to be a subgroup of G.*

***Example 1.7:*** Let G = ($Z_{12}$, +) be a group. H = {0, 2, 4, 6, 8, 10} $\subseteq Z_{12}$ under + is a subgroup of G.

***Example 1.8:*** Let (Z, +) be a group under addition (Z the set of positive, negative integers with zero) (3Z, +) $\subseteq$ (Z, +) is a subgroup of G.

***Example 1.9:*** Let (Q, +) be a group under addition. Q the set of rationals. Q is an abelian group of infinite order.

(Z, +) $\subseteq$ (Q, +) is a subgroup, (Z, +) is also of infinite order. (Q, +) has no subgroup of finite order.

***Example 1.10:*** Let (R, +) be the group under addition. (R, +) is an infinite abelian group. (R, +) has infinitely many subgroups of infinite order.

***Example 1.11:*** Let (R \ {0}, ×) be a group under multiplication of infinite order. Clearly R \ {0} is an abelian group;
H = {1, –1} $\subseteq$ R \ {0} is a subgroup of R \ {0} of order two.

***Example 1.12:*** G = (Q \ {0}, ×) is an abelian group of infinite order. H = {1, –1} $\subseteq$ (Q \ {0}, ×) is a subgroup of finite order in G.

Interested reader can refer and read the classical results on groups [2].

Now we proceed onto define normal subgroup of a group.

**DEFINITION 1.3**: *Let (G, *) be a group. (H, *) be a subgroup of G we say, H is a normal subgroup of G if gH = Hg or $gHg^{-1}$ = H for all g $\in$ G.*

We will illustrate this with examples. First it is important to note that in case of commutative groups every subgroup is normal.

We have so far given only examples of groups which are commutative.

***Example 1.13:*** Let

$$G = \left\{ \begin{pmatrix} a & b \\ c & d \end{pmatrix} \middle| \; ad - bc \neq 0,\, a, b, c, d \in Q^+ \cup \{0\} \right\}.$$

G is a non commutative group under multiplication of matrices.

It is easily verified all subgroups of G are not normal in G,

However $H = \left\{ \begin{pmatrix} x & 0 \\ 0 & x \end{pmatrix} \middle| \; x \in Q^+ \right\}$ is a subgroup of G and H is not a normal in G.

Also $P = \left\{ \begin{pmatrix} x & 0 \\ 0 & y \end{pmatrix} \middle| \; x, y \in Q^+ \right\}$ is a subgroup of G and is P not a normal subgroup of G.

We now proceed onto recall the definition of cosets of a group.

**DEFINITION 1.4:** *Let G be a group H a subgroup of G. For $a \in G$, $aH = \{ah \mid h \in H\}$ is defined as the right coset of H in G.*

We have the following properties to be true.

(i) There is a one to one correspondence between any two right cosets of H in G.

(ii) In case H is a normal subgroup of G we have the right coset to be equal with the left coset.

If H is a normal subgroup of G we can define the quotient group or the factor group of H in G as

$$G/H = \{Hb / b \in G\}$$
$$= \{bH / b \in G\} \qquad (Hb = bH \text{ for all } b \in G\}.$$

G/H is a group and if G is finite o (G) / o (H) = o (G/H).

Now we proceed onto recall definition of rings and fields.

**DEFINITION 1.5:** *Let R be a non empty set with the two binary operations + and ×, such that*

*(i) (R, +) is an abelian group*

*(ii) (R, ×) is such that '×' is a closed operation on R and '×' is an associative binary operation on R*

*(iii) If for all a, b, c ∈ R we have*
$a \times (b+c) = a \times b + a \times c$ *and* $(b+c) \times a = b \times a + c \times a$;
*then (R, +, ×) is defined to be a ring.*

*If a × b = b × a for all a, b ∈ R, we define R to be a commutative ring.*

*If R contains 1 called the multiplicative identity such that a × 1 = 1 × a = a for all a ∈ R we define R to be a ring with unit.*

If is important to note that all rings need not be rings with unit.

We will give examples of them.

***Example 1.14:*** Let R = (Z, +, ×); R is a commutative ring with unit.

***Example 1.15:*** Consider $Z_{15}$; $Z_{15}$ under modulo addition and modulo multiplication is a commutative ring with unit.

***Example 1.16:*** Consider $Z_2 = \{0, 1\}$; $Z_2$ is a commutative ring with unit.

***Example 1.17:*** Let $Z_3 = \{0, 1, 2\}$ be the commutative ring with unit.

***Example 1.18:*** Let $3Z = \{3n / n \in Z\}$, 3Z be a commutative ring but 3Z does not contain the unit 1.

***Example 1.19:*** Let 5Z be the commutative ring, but 5Z is a ring and does not contain the unit.

***Example 1.20:*** (Q, +, ×) is a commutative ring with unit.

***Example 1.21:*** R = $\left\{\begin{pmatrix} a & b \\ c & d \end{pmatrix} \middle| a, b, c, d \in Z, +, \times\right\}$ is a non commutative ring with identity $\begin{pmatrix} 1 & 0 \\ 0 & 1 \end{pmatrix}$ under matrix multiplication ×.

***Example 1.22:*** Let P = {n × n matrices with entries from Q}, P is a non commutative ring with unit under matrix addition and multiplication.

We will now recall the definitions of subrings and ideals of a ring.

**DEFINITION 1.6:** *Let (R, +, ×) be a ring. Suppose $\phi \neq S \subseteq R$ be a proper subset of R. If (S, +, ×) is a ring, we define S to be a subring of R.*

We will illustrate this situation by some examples.

***Example 1.23:*** Let R = (Z, +, $\times$) be a ring; S = (3Z, +, $\times$) is a subring of R.

***Example 1.24:*** Let S = (Z, +, $\times$) be the ring of integers, P = {12Z, +, $\times$} is a subring of S.

It is important and interesting to note that a subring of S need not in general have the unit however the ring S has unit.

This is a marked difference between a subgroup of a group, for a subgroup must contain the identity but a subring need not contain the unit 1 with respect to multiplication. Examples 1.23 and 1.24 give subrings of (Z, +, $\times$) which do not contain 1, the unit of Z.

***Example 1.25:*** Let R = (Q, +, $\times$) be a ring. (Z, +, $\times$) is a subring and it contains the unit. However (8Z, +, $\times$) is also a subring of (Q, +, $\times$) but (8Z, +, $\times$) does not contain the unit 1 of R.

We have seen subrings of a ring. All rings given are of infinite order except those given in examples 1.15, 1.16 and 1.17.

***Example 1.26:*** Let $Z_{30}$ be the ring. S = {0, 10, 20} $\subseteq Z_{30}$ is a subring and $1 \notin$ S.

***Example 1.27:*** Let $Z_{12}$ be the ring. Consider P = {0, 2, 4, 6, 8, 10}, P is a subring of $Z_{12}$. P does not contain the unit 1.

***Example 1.28:*** Let $Z_5$ be the ring. $Z_5$ has no proper subrings. Only 0 and $Z_5$ are the subrings of $Z_5$.

***Example 1.29:*** Let $Z_8$ = {0, 1, 2, 3, …, 7} be a ring. Consider H = {0, 4} $\subseteq Z_8$, H is a subring of $Z_8$ and H does not contain the unit.

***Example 1.30:*** $Z_n$ = {0, 1, 2, …, n–1} (n < ∞) is the commutative ring with unit having n elements. If n is a prime, $Z_n$ has no proper subrings.

**DEFINITION 1.7:** *Let R be a ring. H be a subring of R. If for all r ∈ R and h ∈ H; rh ∈ H then we define H to be left ideal.*

*Similarly if for all r ∈ R and h ∈ H, hr ∈ H we define H to be a right ideal. If H is both a left ideal and a right ideal then we define H to be an ideal of R.*

The following observations are important.

(i) If R is a commutative ring the notion of right ideal and left ideal coincides.

(ii) If I is an ideal of a ring R then (0) ∈ I.

(iii) (0) is called the zero ideal of R.

(iv) If I is an ideal of the ring R then 1 the unit of R is not in I.
For if 1 ∈ I then I = R.

We will now give examples of ideals.

***Example 1.31***: Let (Z, +, ×) be the ring of integers. 2Z is an ideal of Z.

***Example 1.32:*** Let (Q, +, ×) be the ring the rationals, Q has no proper ideals.

***Example 1.33:*** Let (R, +, ×) be the ring of reals. R has no proper ideals.

***Example 1.34:*** Let (Z, +, ×) be the ring of integers, nZ is an ideal of Z for n = 2, 3, …, ∞.

***Example 1.35:*** Let $Z_{12}$ = {0, 1, 2, …, 11} be the ring of integers modulo 12.

I = {0, 6} is an ideal of $Z_{12}$. J = {0, 3, 6, 9} is an ideal of $Z_{12}$. K = {0, 2, 4, 6, …, 10} is an ideal of $Z_{12}$.

***Example 1.36:*** Let $Z_7$ = {0, 1, 2, …, 6} be the ring of integers modulo 7. $Z_7$ has no proper ideals.

***Example 1.37:*** Let $Z_{23}$ = {0, 1, 2, …, 23} be the ring of integers modulo 23. $Z_{23}$ has no ideals.

***Example 1.38:*** Let S = $\left\{\begin{pmatrix} a & b \\ c & d \end{pmatrix} \middle| a, b, c, d \in Z\right\}$ be the ring. Clearly S is a non commutative ring.

Consider I = $\left\{\begin{pmatrix} a & 0 \\ b & 0 \end{pmatrix} \middle| a, b \in Z\right\} \subseteq$ S. I is only a left ideal of S.

For consider $\begin{pmatrix} a & b \\ c & d \end{pmatrix}\begin{pmatrix} x & 0 \\ y & 0 \end{pmatrix} = \begin{pmatrix} ax + by & 0 \\ cx + dy & 0 \end{pmatrix} \in$ I for every

$\begin{pmatrix} a & b \\ c & d \end{pmatrix} \in$ S and $\begin{pmatrix} x & 0 \\ y & 0 \end{pmatrix} \in$ I. I is not a right ideal of S for

$\begin{pmatrix} x & 0 \\ y & 0 \end{pmatrix}\begin{pmatrix} a & b \\ c & d \end{pmatrix} = \begin{pmatrix} ax & bx \\ ya & yd \end{pmatrix} \notin$ I for any $\begin{pmatrix} a & b \\ c & d \end{pmatrix} \in$ S and

$\begin{pmatrix} x & 0 \\ y & 0 \end{pmatrix} \in$ I.

Thus we see in general in a non commutative ring R a left ideal is not a right ideal of R.

Now just like quotient groups we can define quotient rings using two sided ideals of the ring R. Let R be a ring and I a two sided ideal of R. R/I = {a + I / a ∈ R} is a ring defined as the quotient ring of R [2].

We will illustrate this situation by some examples.

***Example 1.39:*** Let Z be the ring. I = 2Z is an ideal of Z. Z/2Z=I = {I, 1 + I}. Clearly I acts as the additive identity for Z/I. Since 1 + I + I = 1 + (I+I) = 1 + I, as I is an ideal, I + I = I. Further (1+I) is the unit in Z/I for (1+I) (1+I) = 1+I + 1.I = 1+I.

Thus Z/I is a ring and Z/2Z = Z/I is isomorphic with $Z_2$. For $0 \mapsto$ I and $1 \mapsto$ 1+I so $Z_2 \cong$ Z/I = {I, 1+I}.

***Example 1.40:*** Let Z be the ring. I = 5Z be the ideal of Z. Z/I = Z/5Z = {I, 1+I, 2+I, 3+I, 4+I} is the quotient ring.

***Example 1.41:*** Let Z be the ring. nZ = I be the ideal of Z. Z/I = {I, 1+I, 2+I, …, n–1 + I} is the quotient ring.

***Example 1.42:*** Let $Z_{12}$ be the ring of modulo integers.

I = {0, 6} be an ideal of $Z_{12}$.

$Z_{12}$/I = {I, 1+I, 2+I, 3+I, 4+I, 5+I} is the quotient ring.

Now we have seen examples of quotient rings.

We proceed onto recall the definition of a field.

**DEFINITION 1.8:** *Let (F, +, ×) be such that F is a non empty set. If*

*(i) (F, +) is an abelian group under addition,*

*(ii) {F \ {0}, ×} is an abelian group under multiplication and*

*(iii) if $a \times (b+c) = ab + ac$ for all a, b, c $\in$ F, then we define F to be a field.*

*If in F we have $nx = 0$ for all $x \in F$ only if $n = 0$ where n is any positive integer then we say F is a field of characteristic zero.*

*If for some n; $nx = 0$ for all $x \in F$ where n is the smallest such number then we say F is a field of characteristic n. Clearly n will only be a prime.*

We will now proceed onto give examples of fields.

***Example 1.43:*** Let $Z_5 = \{0, 1, 2, 3, 4\}$ be a field and $Z_5$ is a field of characteristic 5.

***Example 1.44:*** Let $Z_{43} = \{0, 1, 2, \ldots, 42\}$ be a field and is of characteristic 43.

***Example 1.45:*** $(Q, +, \times)$ is a field and is a field of characteristic zero.

***Example 1.46:*** $(R, +, \times)$ is a field and is of characteristic zero and we see the set of integers, Z is not a field, only a ring.

***Example 1.47:*** $Z_{15}$ is not a field only a ring as $3.5 \equiv 0 \pmod{15}$ so $Z_{15} \setminus \{0\}$ is not a group.

It is easily verified that every field is a ring and a ring in general is not a field.

We recall the definition of a subfield.

**DEFINITION 1.9:** *Let $(F, +, \times)$ be a field. A proper subset $P \subseteq F$ such that if $(P, +, \times)$ is a field, then we define P to be a subfield, of F. If F has no subfield then we call F to be a prime field.*

We will give examples of prime fields.

***Example 1.48:*** Let Q be the field. Q has no proper subfield. Hence Q is the prime field of characteristic zero.

***Example 1.49:*** Let (R, +, $\times$) be the field of reals, R is not a prime field of characteristic zero as (Q, +, $\times$) $\subseteq$ (R, +, $\times$) is a proper subfield of R.

***Example 1.50:*** $Z_2 = \{0, 1\}$ is the prime field of characteristic two as $Z_2$ has no proper subfields.

***Example 1.51:*** $Z_3 = \{0, 1, 2\}$ is the prime field of characteristic three. $Z_3$ has no proper subfields.

***Example 1.52:*** $Z_{17} = \{0, 1, 2, \ldots, 16\}$ is the prime field of characteristic 17.

***Example 1.53:*** $Z_p = \{0, 1, 2, \ldots, p-1\}$ is the prime field of characteristic p, p a prime.

Now it is an interesting fact that fields cannot have ideals in them. This work is left as an exercise to the reader [2].

Now we see we do not have non prime fields of finite characteristic p, p a prime; we proceed on to give methods for constructing them.

We will recall the definition of polynomial rings.

**DEFINITION 1.10:** *Let R be a commutative ring with unit. x an indeterminate;*

$$R[x] = \left\{ \sum_{i=0}^{\infty} a_i x^i \,\middle|\, a_i \in R \right\},$$

*R[x] under usual polynomial addition and multiplication is a ring called the polynomial ring. R[x] is a again a commutative ring with unit and R $\subseteq$ R[x].*

We will give examples of them.

***Example 1.54:*** Let R be the field of reals. R[x] is the polynomial ring in the indeterminate x.

Clearly R[x] is not a field; R[x] is only a commutative ring with unit.

***Example 1.55:*** Let Q be the field of rationals. Q[x] is the polynomial ring.

***Example 1.56:*** Let $Z_2$ be the field of characteristic two.

$Z_2[x]$ is the polynomial ring.

***Example 1.57:*** Let Z be the ring of integers. Z[x] is the polynomial ring.

***Example 1.58:*** Let $Z_5$ be the field of characteristic five. $Z_5[x]$ is the polynomial ring.

***Example 1.59:*** Let $Z_6$ be the ring of integers modulo 6. $Z_6[x]$ is the polynomial ring.

***Example 1.60:*** $Z_2$ [x] is the polynomial ring; let $p(x) = x^2 + 1$ be in $Z_2[x]$; p(x) is reducible in $Z_2$ [x] for $p(x) = (x+1)^2$.

Consider $q(x) = x^2 + x + 1 \in Z_2[x]$; q(x) is irreducible as q(x) cannot be written as the product of polynomials of degree one in $Z_2[x]$.

Inview of this we will define the notion of reducible polynomial and irreducible polynomial.

Let $p(x) \in R[x]$ be a polynomial of degree n in R[x] with coefficients from R. We say p(x) is reducible in R[x] if p(x) = g(x) b(x) where g(x) and b(x) are polynomials in R(x) with deg g(x) < deg p(x) and deg b(x) < deg p(x).

If p(x) is not reducible that is p(x) = g(x) b(x) implies deg g(x) = deg p(x) or deg b(x) = deg p(x) then we define p(x) to be irreducible in R[x].

$p(x) = x^2 + 1 \in R[x]$; R reals is an irreducible polynomial in R[x].

$p(x) = x^2 + 1 \in Z_2[x]$ is reducible in $Z_2[x]$ as

p(x) = (x+1) (x+1).

$p(x) = x^2 - 4 \in R[x]$ is a reducible polynomial as $x^2 - 4 = (x+2)(x-2)$ where R is the set of reals.

$p(x) = x^2 + q$ is irreducible polynomial in R[x] where q is a positive value in R, R-reals.

We have seen examples of reducible and irreducible polynomials.

We first see the concept of reducibility or irreducibility is dependent on the ring or the field over which they are defined. The purpose of these irreducible polynomials is they generate a maximal ideal and this concept will be used in the construction of fields and finite fields.

Let Z[x] be the polynomial ring. $x^2 + 1$ is a irreducible polynomial in Z[x].

Let I be the ideal generated by $x^2 + 1$. Thus
I = {all polynomials of degree greater than or equal to two}

$$\frac{Z[x]}{I} = \{I, b+I, ax + I, (b + ax) + I \,/\, a, b \in Z\}$$
$$= \{a + bx + I \,/\, a, b \in Z\}.$$

Suppose we replace Z by $Z_2$ and consider the polynomial ring $Z_2[x]$; $x^2 + 1$ is a reducible polynomial in $Z_2[x]$.

Let I be the ideal generated by $x^2 + 1$.

Consider $\frac{Z_2[x]}{I}$ = {I, 1 + I, x + I, x + 1 + I}

I acts as the additive identity 1 + I acts as the multiplicative identity.

Clearly $\frac{Z_2[x]}{I}$ is not a field; for

$(1 + x + I)(1 + x + I) = (1 + x)^2 + I = x^2 + 1 + I = I$ is a zero

divisor so $\frac{Z_2[x]}{I}$ is only a ring.

Now if we consider $p(x) = x^2 + x + 1 \in Z_2[x]$, clearly p(x) is irreducible in $Z_2[x]$.

We see $\frac{Z_2[x]}{\langle x^2 + x + 1\rangle} = \frac{Z_2[x]}{I}$

= {I, 1 + I, x + I, 1 + x + I}

where I is the ideal generated by p(x) and I acts as the additive identity.

Now 1 + I is the multiplicative identity.

Further x + I is the inverse of 1 + x + I. To prove this we must show

$(x + I)(1 + x + I) = 1 + I.$

Now $(x + I)(1 + x + I) =$

$$\begin{aligned} x(1 + x) + I &= x^2 + x + I \\ &= (x^2 + x + 1) + 1 + I \\ &= 1 + I \text{ as } 1 + x + x^2 \in I. \end{aligned}$$

Thus $\frac{Z_2[x]}{I}$ is a field with four elements and is of characteristic two so $\frac{Z_2[x]}{I}$ is not a prime field of characteristic two as P = {I, 1+I} $\subseteq \frac{Z_2[x]}{I}$ is a subfield of $\frac{Z_2[x]}{I}$.

Now based on this we give the following facts interested reader can refer [2].

Let $Z_p[x]$ be a polynomial ring with coefficients from the field $Z_p$, p a prime. p (x) $\in Z_p[x]$ be a polynomial of degree n. Suppose p(x) is irreducible over $Z_p$ and I be the ideal generated by p (x).

Now $\frac{Z_p[x]}{I}$ is a field and this field has $p^n$ elements in it and its characteristic is p.

Thus this method gives one the way of constructing non prime fields of characteristic p.

Suppose $V = \frac{Z_2[x]}{I = \langle x^3 + x + 1 \rangle}$

$= \{I, 1 + I, x + I, x^2 + I, 1 + x + I, 1 + x^2 + I, x + x^2 + I, 1 + x + x^2 + I\}$.

$V = \frac{Z_2[x]}{I}$ is a field of characteristic two and V has $2^3$ elements and the irreducible polynomial $x^3 + x + 1$ is of degree 3 so V has $2^3$ elements in it. Suppose $p(x) = x^3 + 1$, clearly p(x) is reducible in $Z_2[x]$ for $x^3 + 1 = (x + 1)(x^2 + x + 1)$.

Consider $\frac{Z_2[x]}{\langle p(x) \rangle} = \frac{Z_2[x]}{I}$ = V = {I, 1 + I, x + I, $x^2$ + I, 1 + x + I, 1 + $x^2$ + I, x + $x^2$ + I, 1 + x + $x^2$ + I}. V is not a field for V has zero divisors. I is the additive identity or zero in V. Consider 1 + x + I, 1 + x + $x^2$ + I in V. (1+x+I) (1+ $x^2$+x) = I (using the fact I is an ideal of $Z_2$ [x] so (x+1) I = I and (1+x+$x^2$)I = I and I + I + I).

(1 + x) (1 + $x^2$ + x) + I = 1 + $x^2$ + x + x + $x^2$ + $x^3$ + I = 1 + $x^3$ + I (as 2x = 0 and $2x^2$ = 0)

= I as 1 + $x^3$ $\in$ I.

Thus V has zero divisors; so V is not a field only a ring as p(x) is a reducible polynomial.

Thus we wish to make a mention that the ideals generated by a single element will be known as principal ideals. Also we expect the reader to be familiar with the notion of minimal ideals, maximal ideals and prime ideals.

However we just give some examples of them.

***Example 1.61:*** Let $Z_{12}$ = {0, 1, 2, …, 11} be the ring of modulo integers.

Consider I = {0, 2, 4, 6, 8, 10} $\subseteq$ $Z_{12}$ to be a maximal ideal of $Z_{12}$. Take J = {0, 6} $\subseteq$ $Z_{12}$, J is a minimal ideal of $Z_{12}$. However K = {0, 3, 6, 9} $\subseteq$ $Z_{12}$ is also a maximal ideal. T = {0, 4, 8} $\subseteq$ $Z_{12}$ is also a minimal ideal of $Z_{12}$. Thus a ring can have many maximal ideals and many minimal ideals.

However we show a ring which has no minimal ideals.

***Example 1.62:*** Let Z be the ring of integers. Z has several maximal ideals, pZ is a maximal ideal where p is a prime.

For all primes p we obtain an infinite collection of maximal ideals.

However Z has no minimal ideal but it has infinite number of ideals which are neither maximal nor minimal. For nZ when n is a non prime is neither a maximal nor a minimal ideal.

2Z, 5Z, 3Z, 7Z, 11Z are maximal ideals.

When R is a commutative ring with unit and I a maximal ideal of R; then R / I is a field.

This property is used in the construction of finite fields of desired order.

We consider $Z_p[x]$; p a prime $Z_p[x]$, the polynomial ring. I be the ideal generated by an irreducible polynomial p(x); then

$\left(\frac{Z_p[x]}{I}\right)$ is a field.

If p(x) is of degree n then $Z_p$ [x] / I is a field of order $Z_p^n$.

Thus this method helps us in generating finite fields which are not prime.

We will give some examples.

***Example 1.63:*** Let $Z_5$ [x] be the polynomial ring.

Consider $p(x) = x^2 + x + 1$ in $Z_5[x]$. p(x) is irreducible in $Z_5[x]$. Let I be the ideal generated by p(x). Now

$Z_5[x]$ / I = {I, 1 + I, 2 + I, 3 + I, 4 + I, x + I, 2x + I, 3x + I, 4x + I, 1 + x + I, 2x + 3 + I, …, 4x + 4 + I}. It is easily verified $Z_5[x]$ / I is a field with 25 elements in it.

However $Z_5$ [x] / I is not a prime field for it has a subfield {I, 1 + I, 2 + I, 3 + I, 4 + I} $\cong Z_5$.

But characteristic of the field $Z_5$ [x] / I is five.

Now if we take $p(x) = x^2 + 2x + 2$, we see p(x) is a reducible polynomial in $Z_5[x]$. For $(x + 3)(x + 4) = x^2 + 2x + 2$. So –3, –4 or 1, 2 are the roots of p(x).

$Z_5[x] / I = \langle p(x) \rangle$ = {I, 1 + I, 2 + I, …, x + 3 + I, x + 4 + I} is not a field only a commutative ring with unit 1 + I.

For consider x + 3 + I and x + 4 + I in $Z_5[x]$ / I.

$$
\begin{aligned}
(x + 3 + I)(x + 4 + I) &= (x + 3)(x + 4) + I \\
&= x^2 + 2x + 2 + I \\
&= I.
\end{aligned}
$$

(I is the zero element in $Z_5[x]$ / I).

Hence $Z_5$ [x] / I cannot be a field.

Likewise we can create fields or rings of finite order.

***Notation:*** If $p(x) = a_0 + a_1x + \ldots + a_nx_n$ then p(x) has a representation in n + 1 tuples given by $(a_0, a_1, \ldots, a_n)$, $a_i$'s are coefficients using which the polynomial ring R[x] is defined and $p(x) \in R[x]$.

Now we just recall the definition of a primitive polynomial.

A polynomial $p(x) = a_0 + a_1x + \ldots + a_nx_n$ in R[x] is said to be a primitive polynomial (where $a_0, a_1, \ldots, a_n$ are integers) if the greatest common divisor of $a_0, a_1, \ldots, a_n$ is 1.

We will use this concept also in the construction of linear codes. Now we proceed onto recall the definition of vector space as we associate or define the algebraic code as a subspace of a vector space.

**DEFINITION 1.11:** *Let V be an additive abelian group. F a field. We say V is a vector space over the field F if the following conditions are true.*

*(1) For every $v \in V$ and $a \in F$ we have av and va are in V (av = va and conventionally we write as av a is called the scalar and v a vector).*

*(2) $a(v_1 + v_2) = av_1 + av_2$.*

*(3) $(a+b)v = av + bv$.*

*(4) $a(bv) = (ab)v$.*

*(5) $1.v = v$ for all $a, b \in F$ and $v_1, v_2, v \in V$ and $0.v = 0 \in V$.*

(Here 1 represents the unit element of F).

We give examples of vector spaces.

***Example 1.64:*** Let $V = Z_2 \times Z_2 \times Z_2 \times Z_2$ be a group under addition, V is a vector space over the field $Z_2$.

***Example 1.65:*** Let V = Q[x] be an additive abelian group. V is a vector space over the field Q.

***Example 1.66:*** Let $V = Z_2[x]$, polynomial group under addition. V is a vector space over $Z_2$.

We recall the definition of a vector subspace of a vector space.

**DEFINITION 1.12:** *Let V be a vector space over the field F. $W \subseteq V$ be a proper subset of V. If W itself is a vector space over the field F then we call W to be vector subspace of V over the field F.*

We will illustrate this situation by some examples.

***Example 1.67:*** Let $V = Z_2 \times Z_2 \times Z_2 \times Z_2$ be a vector space over the field $Z_2$. Consider $W = \{Z_2 \times \{0\} \times Z_2 \times \{0\}\} \subseteq V$; W is a vector subspace of V over $Z_2$.

***Example 1.68:*** Let V = Q[x] be the vector space over Q.

Consider

$$P = \left\{ \sum_{i=0}^{\infty} a_i x^{2i} \,\middle|\, a_i \in Q \right\} \subseteq V;$$

P is a vector subspace of V over Q.

***Example 1.69:*** Let $V = Q \times Q \times Q \times Q \times Q \times Q$ be a vector space over Q. Consider $W = Q \times \{0\} \times Q \times \{0\} \times Q \times Q \subseteq V$; W is a vector subspace of V over Q.

***Example 1.70:*** Let $V = R \times R \times R \times R \times R$ be a vector space over Q. Consider

$W = Q \times Q \times Q \times Q \times Q \subseteq R \times R \times R \times R \times R$; W is a vector subspace of V over Q.

Now the reader is expected to recall the definition of basis and dimension of a vector space [2, 16].

We will be using only finite dimensional vector spaces that too defined over the finite field. Infact we will be mainly using the field $Z_2$ or $Z_2^n$ of order $2^n$ and of characteristic two; denoted by $GF(2^n)$ as most of the messages transmitted are binary, it is sufficient to study over $Z_2$ or $GF(2^n)$.

Now we just recall the definition of Hamming distance and Hamming weight in vector spaces.

We say for any two vectors $v_1 = (x_1, \ldots, x_n)$ and $v_2 = (y_1, \ldots, y_n)$; $v_1$ and $v_2 \in V = F \times \ldots \times F$; n-times where F is a field and V is a vector space over F. The Hamming distance between $v_1$ and $v_2$ denoted by d ($v_1$, $v_2$) = number of places in which $v_1$ is different from $v_2$.

For instance if $v_1$ = (1 0 0 1 0 0 1) and $v_2$ = (0 1 1 1 0 0 1) then the Hamming distance between $v_1$ and $v_2$ denoted by d ($v_1$, $v_2$) = 3.

Now Hamming weight x of a vector in V is the distance of x from the zero vector. Thus Hamming weight of x is the number of non zero co ordinates in the vector x, thus d(x, 0) = w(x) = number of non zero coordinates in x.

Thus if x = (1 0 1 0 1 0 0 1 1 1) be a vector in V then w(x) = d (x, 0) = 6.

Now having defined Hamming weight and Hamming distance we now proceed onto define linear codes and illustrate them with examples in chapter two.

**Chapter Two**

# ALGEBRAIC LINEAR CODES AND THEIR PROPERTIES

In this chapter we recall the definition of algebraic linear codes and discuss the various properties associated with them. We give examples and define several types of codes.

Let $x = (x_1, x_2, \ldots, x_n)$ where $x_i \in F_q$, where $F_q$ is a finite field. (q, a power of a prime). In x the first k symbols are taken as message symbols and the remaining $n - k$ elements $x_{k+1}, x_{k+2}, \ldots, x_n$ are check symbols (or control symbols). $x = (x_1, \ldots, x_n)$ is defined as the code word and it can be denoted by $(x_1, \ldots, x_n)$ or $x_1\, x_2 \ldots x_n$ or $x_1, x_2, \ldots, x_n$.

We will now roughly indicate how messages go through the system starting from the sender (Information source).

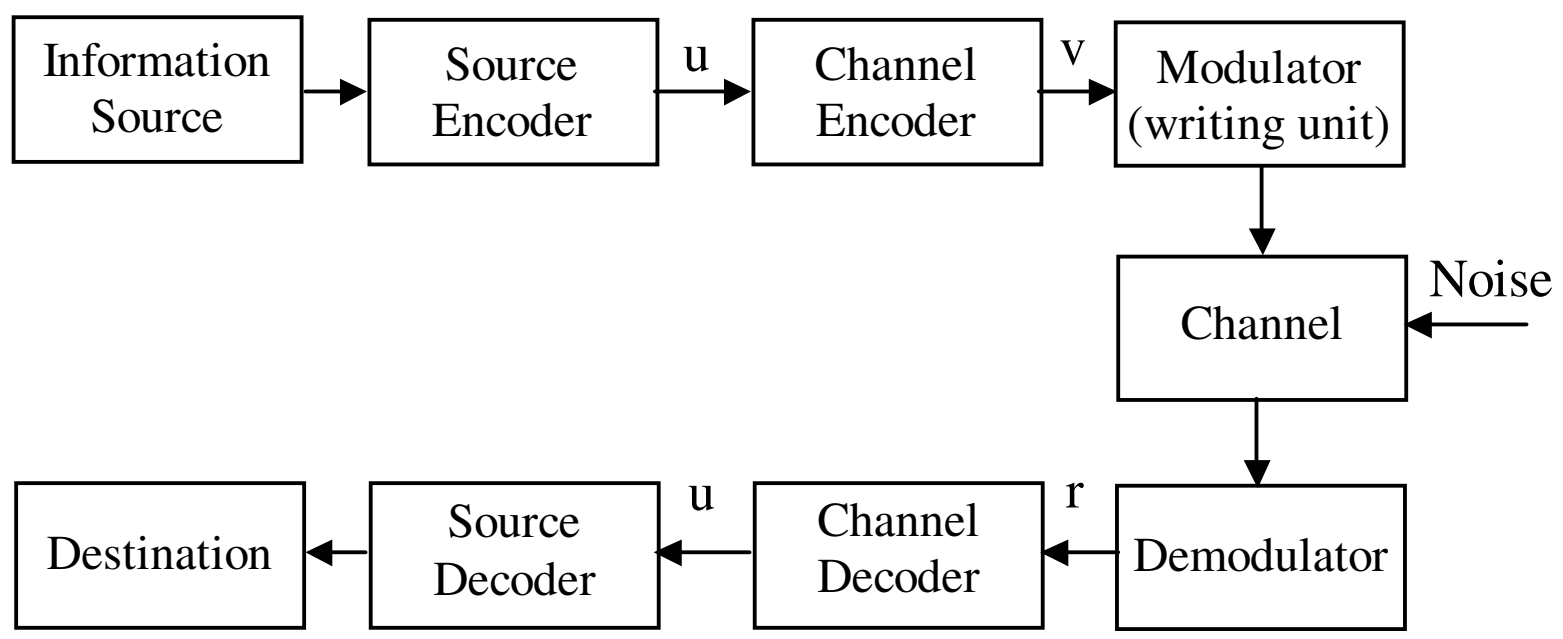

We shall consider senders with finite number of discrete signals (telegraph) in contrast to continuous sources (examples radio).

Further the signals emanating from the source cannot be transmitted directly by the channel. For instance, a binary channel cannot transmit words in the usual latin alphabet. Encoder performs the important work of data reduction and suitably transforms the message into usable form. Thus there is a difference between the source encoding and the channel encoding.

The former reduces the message to its essential parts and the latter adds redundant information to enable detection and correction of possible errors in the transmission.

Also the channel decoding and the source decoding are distinct for they invert the corresponding channel and source encoding besides detecting and correcting error.

The main aim of coding theory is to formulate techniques for transmitting message, free of errors, at a less cost and with a great speed. However in few cases the possibility of repeating messages is acceptable; in some cases repeating messages is impossible; like in case of taking picture of planets or Mars or Moon where there is a steady motion of the unmanned machine. In these cases repetition is impossible and also the cost is very high so high accuracy is expected with minimum or no error. As repeating the message is also time consuming besides being costly.

We want to find efficient algebraic methods to improve the realiability of the transmission of messages.

In this chapter we give only simple coding and decoding algorithms which can be easily understood by a beginner.

Binary symmetric channel is an illustration of a model for a transmission channel. Now we will proceed onto define a linear code algebraically.

Let $x = (x_1 \ldots x_n)$ be a code word with k message symbols and n–k check symbols or control symbols. We know the message symbols are from $F_q$, to determine the check symbols. We obtain the check symbols from the message symbols in such a way that the code words x satisfy the system of linear equations.

$Hx^t = (0)$, where H is the given $n - k \times n$ matrix with elements from $F_q$.

The standard form for H is $(A, I_{n-k})$ with A an $n - k \times k$ matrix and $I_{n-k}$ the $n - k \times n - k$ identity matrix.

The set of all n-dimensional vectors $x = (x_1, \ldots, x_n)$ satisfying $Hx^t = (0)$ over $F_q$ is called a linear (block) code C over $F_q$ of block length n.

The matrix H is called the parity check matrix of the linear (n,k) code C. If q = 2 then we call C a binary code. k/n is called transmission (or information) rate.

Since C under addition is a group we call C as a group code. Also C can be defined as the null space of the matrix H.

We will first illustrate this situation by some examples.

***Example 2.1:*** Let q = 2, n = 7 and k = 4 and C be a C(7, 4) code with entries from $Z_2$. The message $a_1$, $a_2$, $a_3$, $a_4$ is encoded as the code word x = $a_1\ a_2\ a_3\ a_4\ \ x_5\ x_6\ x_7$. Here the check symbols are $x_5\ x_6\ x_7$, such that for the given parity check matrix

$$H = \begin{pmatrix} 0 & 0 & 1 & 0 & 1 & 1 & 1 \\ 0 & 1 & 0 & 1 & 1 & 1 & 0 \\ 1 & 0 & 1 & 1 & 1 & 0 & 0 \end{pmatrix},$$

with the set of message symbols from $Z_2^4$ the code words are given by

$$\begin{matrix} 0 & 0 & 0 & 0 \\ 1 & 0 & 0 & 0 \\ 0 & 1 & 0 & 0 \\ 0 & 0 & 1 & 0 \\ 0 & 0 & 0 & 1 \end{matrix} \quad \begin{matrix} 1 & 1 & 0 & 0 \\ 1 & 0 & 1 & 0 \\ 1 & 0 & 0 & 1 \\ 0 & 1 & 1 & 0 \\ 0 & 1 & 0 & 1 \end{matrix} \quad \begin{matrix} 0 & 0 & 1 & 1 \\ 1 & 1 & 1 & 0 \\ 1 & 1 & 0 & 1 \\ 1 & 0 & 1 & 1 \end{matrix} \quad \begin{matrix} 0 & 1 & 1 & 1 \\ 1 & 1 & 1 & 1 \end{matrix}.$$

To find C. C = $\{x \in Z_2^7 / Hx^t = (0)\}$. Hence $|C| = 2^4 = 16$.

$$Hx^t = \begin{pmatrix} 0 & 0 & 1 & 0 & 1 & 1 & 1 \\ 0 & 1 & 0 & 1 & 1 & 1 & 0 \\ 1 & 0 & 1 & 1 & 1 & 0 & 0 \end{pmatrix} (x_1\ x_2\ x_3\ x_4\ a_5\ a_6\ a_7)^t = (0)$$

gives

$$C = \begin{matrix} \{0 & 0 & 0 & 0 & 0 & 0 & 0 & & 1 & 0 & 1 & 0 & 0 & 0 & 1 \\ 1 & 0 & 0 & 0 & 1 & 1 & 0 & & 1 & 0 & 0 & 1 & 0 & 1 & 1 \\ 0 & 1 & 0 & 0 & 0 & 1 & 1 & & 0 & 1 & 0 & 1 & 0 & 0 & 0 \\ 0 & 0 & 1 & 0 & 1 & 1 & 0 & & 1 & 1 & 1 & 0 & 0 & 1 & 0 \\ 0 & 0 & 0 & 1 & 1 & 0 & 1 & & 1 & 0 & 1 & 1 & 1 & 0 & 0 \\ 1 & 1 & 0 & 0 & 1 & 0 & 1 & & 1 & 1 & 0 & 1 & 0 & 0 & 0 \\ 0 & 1 & 1 & 0 & 1 & 0 & 0 & & 0 & 1 & 1 & 1 & 0 & 0 & 1 \\ 0 & 0 & 1 & 1 & 0 & 1 & 0 & & 1 & 1 & 1 & 1 & 1 & 1 & 1\} \end{matrix}$$

***Example 2.2:*** Take n = 7, k = 4 and q = 2. To construct the C(7, 4) code using the parity check matrix

$$H = \begin{pmatrix} 1 & 1 & 1 & 0 & 1 & 0 & 0 \\ 0 & 1 & 1 & 1 & 0 & 1 & 0 \\ 0 & 0 & 1 & 1 & 1 & 0 & 1 \end{pmatrix}.$$

$C = \{Hx^T = (0) \text{ where } x \in Z_2^4\}$.

C = {0 0 0 0 0 0 0, 1 0 0 0 1 0 1, 0 1 0 0 1 1 1, 0 0 1 0 1 1 0, 0 0 0 1 0 1 1, 1 1 0 0 0 1 0, 1 0 1 0 0 1 1, 1 0 0 1 1 1 0, 0 1 1 0 0 0 1, 0 0 1 1 1 0 1, 0 1 0 1 1 0 0, 1 1 1 0 1 0 0, 1 1 0 1 0 0 1, 1 0 1 1 0 0 0, 0 1 1 1 0 1 0, 1 1 1 1 1 1 1}.

We see the two codes given in examples 2.1 and 2.2 are C(7, 4) codes but they are different as their parity check matrices are different. Further both the codes have the same set of message symbols.

***Example 2.3:*** Let C be a (4, 2) code given by the parity check matrix

$$H = \begin{pmatrix} 1 & 1 & 1 & 0 \\ 0 & 1 & 0 & 1 \end{pmatrix}.$$

$C = \{x \in Z_2^4 / Hx^t = (0)\}$.

= {0 0 0 0, 1 0 1 0, 0 1 1 1, 1 1 0 1}.

Suppose we consider the (4, 2) code using another parity check matrix.

$$H_1 = \begin{pmatrix} 1 & 0 & 1 & 0 \\ 1 & 1 & 0 & 1 \end{pmatrix} \text{ we get}$$

$C = \{x \in Z_2^4 / H_1x^t = (0)\}$.

= {0 0 0 0, 1 0 1 1, 0 1 0 1, 1 1 1 0}.

We see both the codes are different though they have the same set of message symbols.

We now recall how the repetition code is constructed and the parity check matrix associated with it. If each code word of a code consists of only one message symbol $x_1 \in Z_2$ and n–1 check symbols $a_2 = a_3 = \ldots = a_n$ are all equal to $x_1$ ($x_1$ is repeated n–1 times). Thus we obtain a binary (n, 1) code with parity - check matrix;

$$H = \begin{pmatrix} 1 & 1 & 0 & \ldots & 0 \\ 1 & 0 & 1 & \ldots & 0 \\ \vdots & \vdots & \vdots & & \vdots \\ 1 & 0 & 0 & \ldots & 1 \end{pmatrix}.$$

There are only two code words in this code namely (0 0 … 0) and (1 1 … 1).

This code is used when it impossible and impracticable or too costly to send original message more than once, like transmission of information from space crafts or satellites where it is impossible to use ARQ protocols owing to time limitations. Moving space crafts which takes photos of heavenly bodies is an example where this code can be used.

***Example 2.4:*** Let C(5, 1) be a binary code obtained from the parity check matrix;

$$H = \begin{pmatrix} 1 & 1 & 0 & 0 & 0 \\ 1 & 0 & 1 & 0 & 0 \\ 1 & 0 & 0 & 1 & 0 \\ 1 & 0 & 0 & 0 & 1 \end{pmatrix}.$$ The two code words are 1 1 1 1 1 and 0 0 0 0 0.

Now we proceed onto describe the parity-check code.

Parity check code is a (n, n–1) code where we have n–1 message symbols and one check symbol. The parity check matrix is H = (1 1 … 1).

Each code word has only one check symbol and H has only even number of ones.

These codes are used in banking where the last digit of the account number, usually is a control digit.

***Example 2.5:*** Let H = (1 1 1 1 1 1) be the parity check code. C is a (6, 5) code and

| | | |
|---|---|---|
| C ={0 0 0 0 0 0 | 1 1 1 1 1 1 | 1 1 0 0 0 0 |
| 1 0 0 0 0 1 | 1 0 1 0 0 0 | 1 0 0 0 1 0 |
| 1 0 0 1 0 0 | 0 1 1 0 0 0 | 0 1 0 1 0 0 |
| 0 1 0 0 1 0 | 0 1 0 0 0 1 | 0 0 1 1 0 0 |
| 0 0 1 0 1 0 | 0 0 1 0 0 1 | 0 0 0 1 1 0 |
| 0 0 0 1 0 1 | 0 0 0 0 1 1 | 1 1 1 1 0 0 |
| 1 1 1 0 1 0 | 1 1 1 0 0 1 | 1 1 0 1 1 0 |
| 1 1 0 1 0 1 | 1 1 0 0 1 1 | 1 0 1 1 1 0 |
| 1 0 1 1 0 1 | 1 0 1 0 1 1 | 1 0 0 1 1 1 |
| 0 1 0 1 1 1 | 0 0 1 1 1 1 | 1 0 1 1 1 0 |
| 0 1 1 1 0 1 | 0 1 1 1 1 0} | |

is code associated with the parity check matrix H.

Now we will proceed onto describe the canonical generator matrix of a linear (n, k) code. Suppose $H = (A, I_{n-k})$ is the parity check matrix associated with the (n, k) code then the generator matrix $G = (I_k - A^t)$ is such that $GH^T = (0)$.

Further every code word $x = (x_1, \ldots, x_n) = (a_1, \ldots, a_k)\, G$.

We will now describe this situation by some examples.

***Example 2.6:*** Let H = $\begin{pmatrix} 0 & 1 & 1 & 0 & 1 & 0 & 0 \\ 1 & 0 & 0 & 1 & 0 & 1 & 0 \\ 0 & 1 & 1 & 1 & 0 & 0 & 1 \end{pmatrix}$ be a parity check matrix of a C(7, 4) code.

The generator matrix G associated with this (7, 4) code with the parity check matrix H is given by

$$G = \begin{pmatrix} 1 & 0 & 0 & 0 & 0 & 1 & 0 \\ 0 & 1 & 0 & 0 & 1 & 0 & 1 \\ 0 & 0 & 1 & 0 & 1 & 0 & 1 \\ 0 & 0 & 0 & 1 & 0 & 1 & 1 \end{pmatrix}.$$

Now using the message symbols from $Z_2^4$ we get the following code words generated by G.

C = {0 0 0 0 0 0 0, 1 0 0 0 0 1 0, 0 1 0 0 1 0 1, 0 0 1 0 1 0 1, 0 0 0 1 0 1 1, 1 1 0 0 1 1 1, 1 0 1 0 1 0 1, 1 0 0 1 0 0 1, 0 1 1 0 0 0 0, 0 1 0 1 1 1 0, 0 0 1 1 1 1 0, 1 1 1 0 0 1 0, 0 1 1 1 0 1 1, 1 1 0 1 1 0 0, 1 0 1 1 1 0 0, 1 1 1 1 0 0 1}.

***Example 2.7:*** Let H = $\begin{pmatrix} 1 & 0 & 1 & 1 & 0 & 0 & 0 \\ 0 & 1 & 1 & 0 & 1 & 0 & 0 \\ 1 & 1 & 1 & 0 & 0 & 1 & 0 \\ 0 & 1 & 0 & 0 & 0 & 0 & 1 \end{pmatrix}$ be parity check matrix associated with a C = C(7, 3) code.

The associated generator matrix

$$G = \begin{pmatrix} 1 & 0 & 0 & 1 & 0 & 1 & 0 \\ 0 & 1 & 0 & 0 & 1 & 1 & 1 \\ 0 & 0 & 1 & 1 & 1 & 1 & 0 \end{pmatrix}.$$

Now we generate the code,

C = {0 0 0 0 0 0 0, 1 0 0 1 0 1 0, 0 1 0 0 1 1 1, 0 0 1 1 1 1 0, 1 0 1 0 1 0 0, 0 1 1 1 0 0 1, 1 1 0 1 1 0 1, 1 1 1 0 0 1 1} is the code generated by G.

$$\text{Consider } GH^T = \begin{pmatrix} 1 & 0 & 0 & 1 & 0 & 1 & 0 \\ 0 & 1 & 0 & 0 & 1 & 1 & 1 \\ 0 & 0 & 1 & 1 & 1 & 1 & 0 \end{pmatrix} \begin{bmatrix} 1 & 0 & 1 & 0 \\ 0 & 1 & 1 & 1 \\ 1 & 1 & 1 & 0 \\ 1 & 0 & 0 & 0 \\ 0 & 1 & 0 & 0 \\ 0 & 0 & 1 & 0 \\ 0 & 0 & 0 & 1 \end{bmatrix}$$

$$= \begin{pmatrix} 0 & 0 & 0 & 0 \\ 0 & 0 & 0 & 0 \\ 0 & 0 & 0 & 0 \end{pmatrix}.$$

Having seen the generator matrix and parity check matrix of a code we now proceed onto analyse the nature of the generator matrix.

A generator matrix G for a linear code C is a $k \times n$ matrix for which the rows are a basis of C.

Now a natural question would be what is the purpose of the parity check matrix H.

We see the parity check matrix serves as the fastest means to detect errors. So error detection is done by the parity check matrix. Suppose y is the received code word then find $S(y) = Hy^T$, S (y) is defined as the syndrome. If S(y) = (0) we say no error and accept y as the correct received word.

If $S(y) \neq (0)$ we declare error has occurred. Thus the parity check matrix helps in detecting the error in the received word. So error detection is not a very difficult task as far as coding is concerned. However it is pertinent to mention that when $S(y) = 0$, y is a code word, it need not be the transmitted code word. In certain cases when the error pattern e is identical to a non zero code word, y is the sum of two code words which is a code word so $Hy^T = (0)$. These errors are not detectable. We accept them as a correct transmitted message. Thus we have $2^k - 1$ non zero code words which can lead to undetectable errors, so we have $2^k - 1$ undetectable error patterns. Now the real problem lies in correcting the error.

We will just describe the coset leader method which is used to correct errors. Once the error is detected, we know every code C is a subspace of the vector space $F_q^n$ where $F_q^n$ is defined over $F_q$.

The factor space $F_q^n / C$ consists of all cosets $a + C = \{a + x \mid x \in C\}$ for any $a \in F_q^n$.

Clearly each coset contains $q^k$ vectors as C has $2^k$ elements in it. Since a coset is either disjoint or identical we get a partition on $F_q^n$ so $F_q^n = C \cup (a^1 + C) \cup \ldots \cup (a^t + C)$ for $t = q^{n-k} - 1$.

If a vector y is received then y must be an element of one of these cosets say $a^i + C$. If the code word $x^1$ has been transmitted then the error vector $e = y - x^1 \in a^1 + C - x^1 = a^1 + C$. Thus we quote the decoding rule [16].

If a vector y is received then the possible error vectors e are the vectors in the coset containing y. The most likely error is the vector with minimum weight in the coset of y. Thus y is decoded as $= y - \bar{e}$. The vector of minimum weight in a coset is called the coset leader. If there are several coset leaders arbitrarily choose any one of them.

Let $a^{(1)}, a^{(2)}, \ldots, a^{(t)}$ be the coset leaders. We have the following table.

$$\begin{array}{llll} x^{(1)} = 0 & x^{(2)} & \ldots & x^{(q^k)} \} \text{ code words in C.} \\ a^{(1)} + x^{(1)} & a^{(1)} + x^{(2)} & \ldots & a^{(1)} + x^{(q^k)} \} \text{ other cosets} \\ \vdots & \vdots & & \vdots \\ a^{(t)} + x^{(t)} & a^{(t)} + x^{(2)} & \ldots & a^{(t)} + x^{(q^k)} \} \text{ other cosets} \end{array}$$

$\underbrace{\qquad\qquad}$
coset leaders

If a vector y is received then we have to find y in the table.

Let $y = a^{(i)} + x$, then the decoder decides that the error $\overline{e}$ is the coset leader $a^{(i)}$. Thus y is decoded as the code word $\overline{x} = y - \overline{e} = x^{(i)}$. The code word $\overline{x}$ occurs as the first element in the column of y. The coset of y can be found by evaluating the so called syndrome.

We will illustrate this situation by an example.

***Example 2.8:*** Let C be a binary linear (4, 2) code with the generator matrix

$$G = \begin{pmatrix} 1 & 0 & 1 & 0 \\ 0 & 1 & 1 & 1 \end{pmatrix} \text{ and parity check matrix}$$

$$H = \begin{pmatrix} 1 & 1 & 1 & 0 \\ 0 & 1 & 0 & 1 \end{pmatrix}.$$

The corresponding coset table is

| Message symbols | 0 0 | 10 | 0 1 | 1 1 |
|---|---|---|---|---|
| Code words | 0 0 0 0 | 1 0 1 0 | 0 1 1 1 | 1 1 0 1 |
| Other cosets | 1 0 0 0 | 0 0 1 0 | 1 1 1 1 | 0 1 0 1 |
| | 0 1 0 0 | 1 1 1 0 | 0 0 1 1 | 1 0 0 1 |
| | 0 0 0 1 | 1 0 1 1 | 0 1 1 0 | 1 1 0 0 |

Coset leaders

If y = 1 1 1 1 is received then

$$S(y) = \begin{pmatrix} 1 & 1 & 1 & 0 \\ 0 & 1 & 0 & 1 \end{pmatrix} \begin{bmatrix} 1 \\ 1 \\ 1 \\ 1 \end{bmatrix} = (1\ 0).$$

Thus error e = 1 0 0 0 and y is decoded as x = y – e = 0 1 1 1 and the corresponding message is 0 1.

Now we will proceed onto define cyclic codes.

We say a code word v in C (C a k dimensional subspace of $F_q^n$) is a cyclic code if $v = (v_1 \ldots v_n)$ is in C then $(v_n\ v_1 \ldots v_{n-1})$ is in C.

We generate cyclic codes using polynomial called the generator polynomial of a cyclic code.

If $g = g_0 + g_1x + \ldots + g_mx^m$
is a generator polynomial then $g / x^n - 1$ and deg $g = m < n$.

Let C be a linear (n, k) code with k = n–m defined by the generator matrix;

$$G = \begin{pmatrix} g_0 & g_1 & \cdots & g_m & 0 & \cdots & 0 \\ 0 & g_0 & \cdots & g_{m-1} & g_m & \cdots & 0 \\ \vdots & \vdots & & \vdots & \vdots & & \vdots \\ 0 & \cdots & \cdots & g_0 & g_1 & \cdots & g_m \end{pmatrix} = \begin{pmatrix} g \\ xg \\ \vdots \\ xg^{k-1} \end{pmatrix}.$$

Then C is cyclic. The rows of G are linearly independent and rank G = k, the dimension of C.

If $x^n - 1 = g_1 \ldots g_t$ is a complete factorization of $x^n - 1$ into irreducible polynomials over $F_q$ then the cyclic codes $(g_i)$ generated by polynomials $g_i$ are called maximal cyclic codes.

A maximal cyclic code is a maximal ideal in $\frac{F_q[x]}{\langle x^n - 1\rangle}$.

If g is the polynomial that generates the cyclic code C then $h = x^n - 1 / g$ is defined as the check polynomial of C.

Thus if $h = \Sigma\ h_i\ x_i$, $h_k \neq 0$ then the parity check matrix associated with the cyclic code C is given by

$$H = \begin{pmatrix} 0 & 0 & \ldots & 0 & h_k & \ldots & h_1 & h_0 \\ 0 & 0 & \ldots & h_k & h_{k-1} & \ldots & h_0 & 0 \\ . & . & . & . & . & . & . & . \\ h_k & h_{k-1} & \ldots & h_0 & 0 & \ldots & 0 & 0 \end{pmatrix}.$$

We will illustrate this situation by an example.

***Example 2.9:*** Let C = C (7, 4) be a code of length 7 with 4 message symbols and q = 2.

Suppose $g(x) = x^3 + x + 1$ be the generator polynomial of the cyclic code C. The check polynomial of the cyclic code C is $h = x^7 - 1 / g = x^4 + x^2 + x + 1$.

Now the generator matrix

$$G = \begin{pmatrix} 1 & 1 & 0 & 1 & 0 & 0 & 0 \\ 0 & 1 & 1 & 0 & 1 & 0 & 0 \\ 0 & 0 & 1 & 1 & 0 & 1 & 0 \\ 0 & 0 & 0 & 1 & 1 & 0 & 1 \end{pmatrix}$$

and the parity check matrix of this cyclic code is

$$H = \begin{pmatrix} 0 & 0 & 1 & 0 & 1 & 1 & 1 \\ 0 & 1 & 0 & 1 & 1 & 1 & 0 \\ 1 & 0 & 1 & 1 & 1 & 0 & 0 \end{pmatrix}.$$

Now we can generate the cyclic code using G. aG = x where a is the message from $Z_2^4$ and x the resulting code word of C. We have referred and information are from [1, 16].

Now we proceed onto just recall the definition of rank distance codes and give some of the properties related with them.

Error detection, error correction for these linear block codes can be found in [1, 6, 16, 18]. However the erasure techniques for these codes are very meagre hence we just describe the erasure techniques in these codes.

When the code symbols are from the Galois field $GF(2^n)$ of an arbitrary dimension the function of the modulator is to match the encoder output to the signals of the transmission channel.

The modulator accepts the binary encoded symbols and produces wave forms or signals appropriate to the physical transmission medium. At the receiving end of the communication link, the demodulator operates on the signals received from a separate transmission symbol interval or a set of elements in {0, 1}. The demodulator is designed to make a definite decision for each received symbol 0 or 1.

The definition of a channel includes the modulator, the demodulator and all intervening transmission equipment and media. Most of them are discrete memoryless channel. The assumption is made in that the output symbol at any instant of time depends statistically only on the input symbol at that time.

The coding system is described by the following figure.

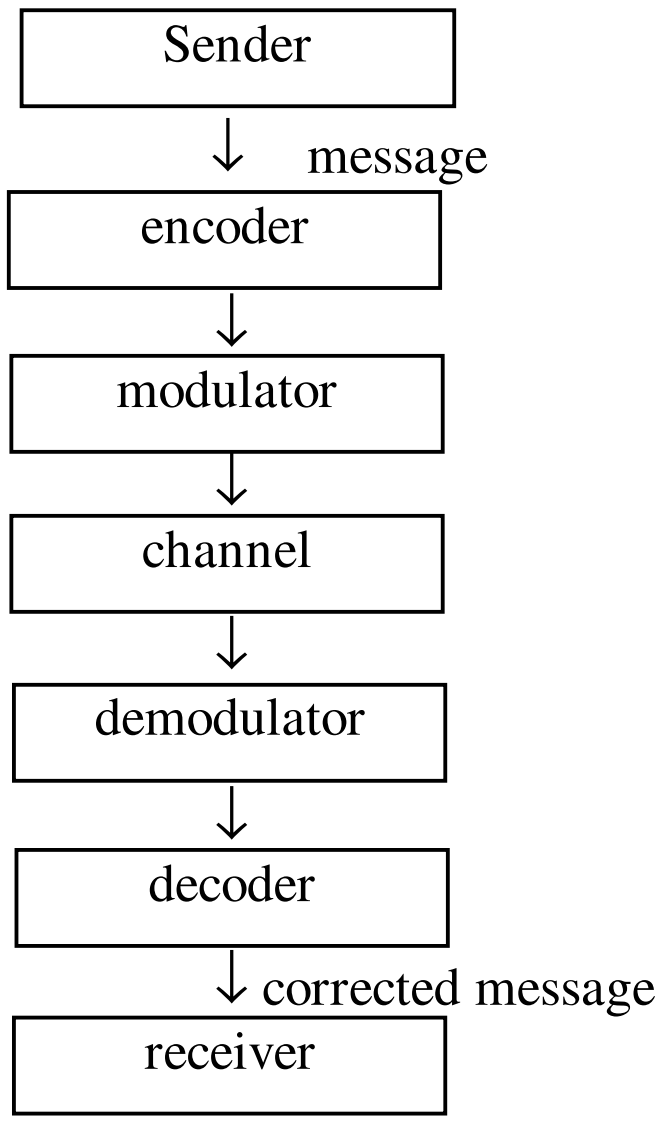

This simple model [20] is known as the binary symmetric channel and is described in the following figure.

In this model the modulator input x has value 0 or 1 and the demodulator output y has value 0 or 1. When 0 is transmitted and 1 is received the error probability is p. When 0 is transmitted and 0 is received the probability is 1 – p. Similarly when 1 is transmitted and 0 is received the probability is p and when 1 is transmitted and 1 is received the probability is 1 – p. This model can be described simply as a binary-input and binary-output model.

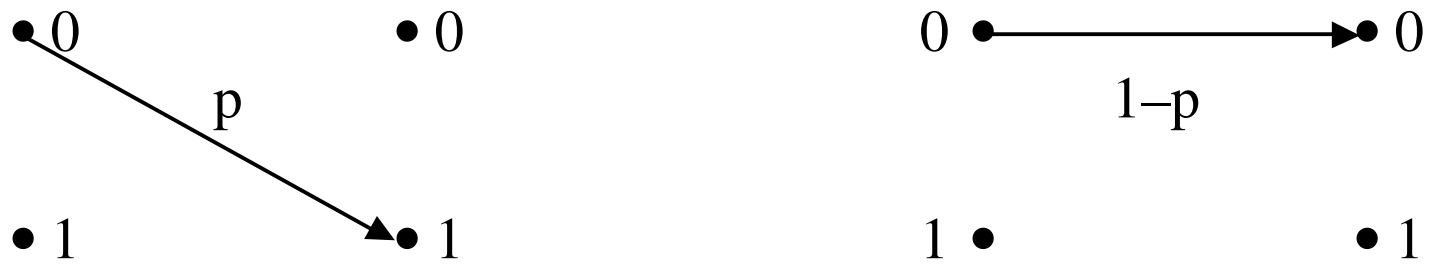

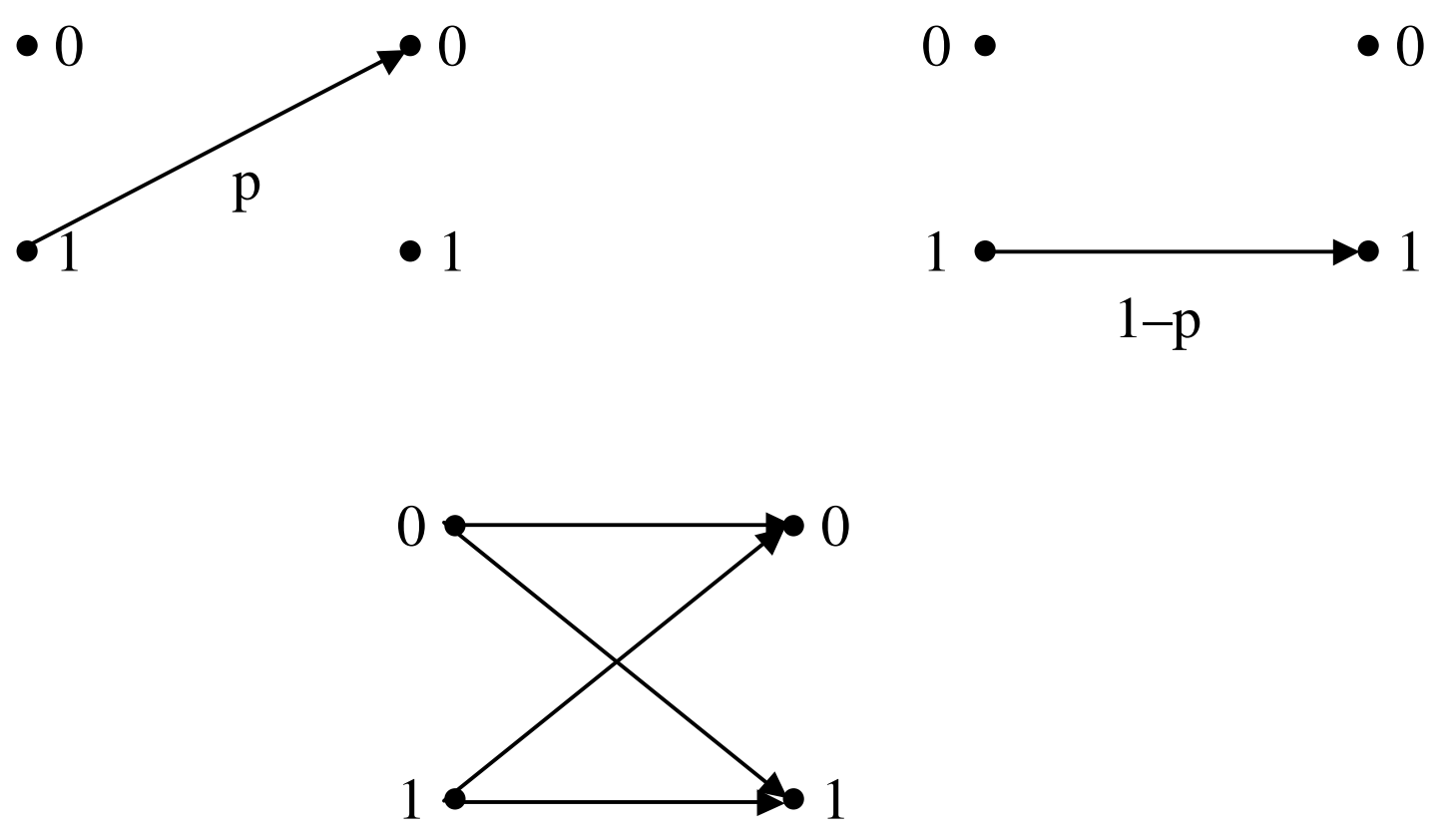

**Binary symmetric channel**

When erasure occurs [20] has studied the binary symmetric erasure channel. This channel model is depicted in the following figure which includes a symmetric transmission from either input symbol to an output symbol labeled ‘?’ to denote ambiguity. Now when an input symbol is sent we have the three possibilities.

- (i) The correct output is received.
- (ii) An erroneous output is received.
- (iii) The demodulator is unable to decide, the result is a blank space, that is ambiguous.

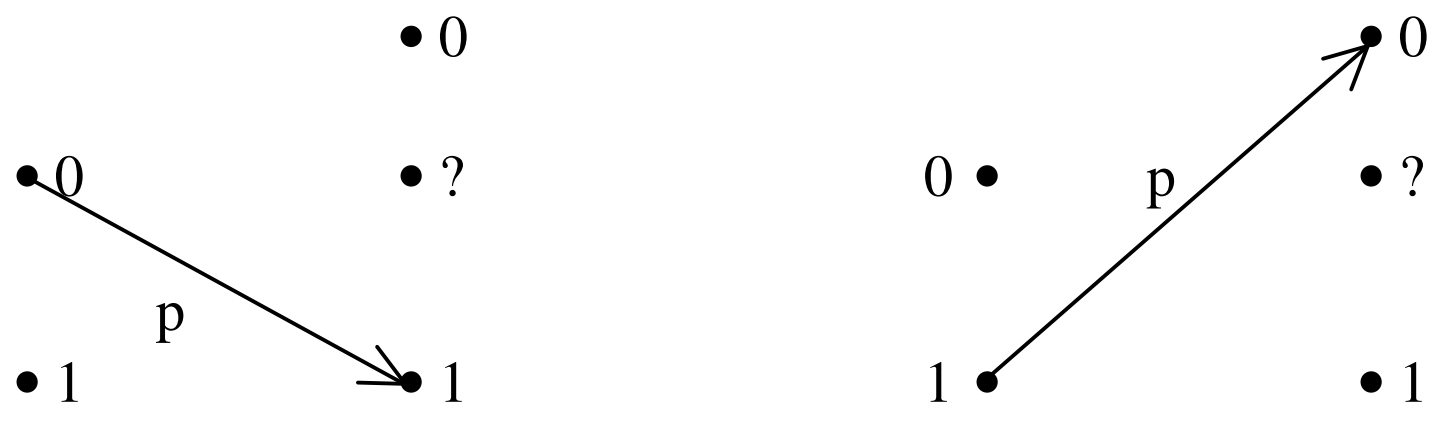

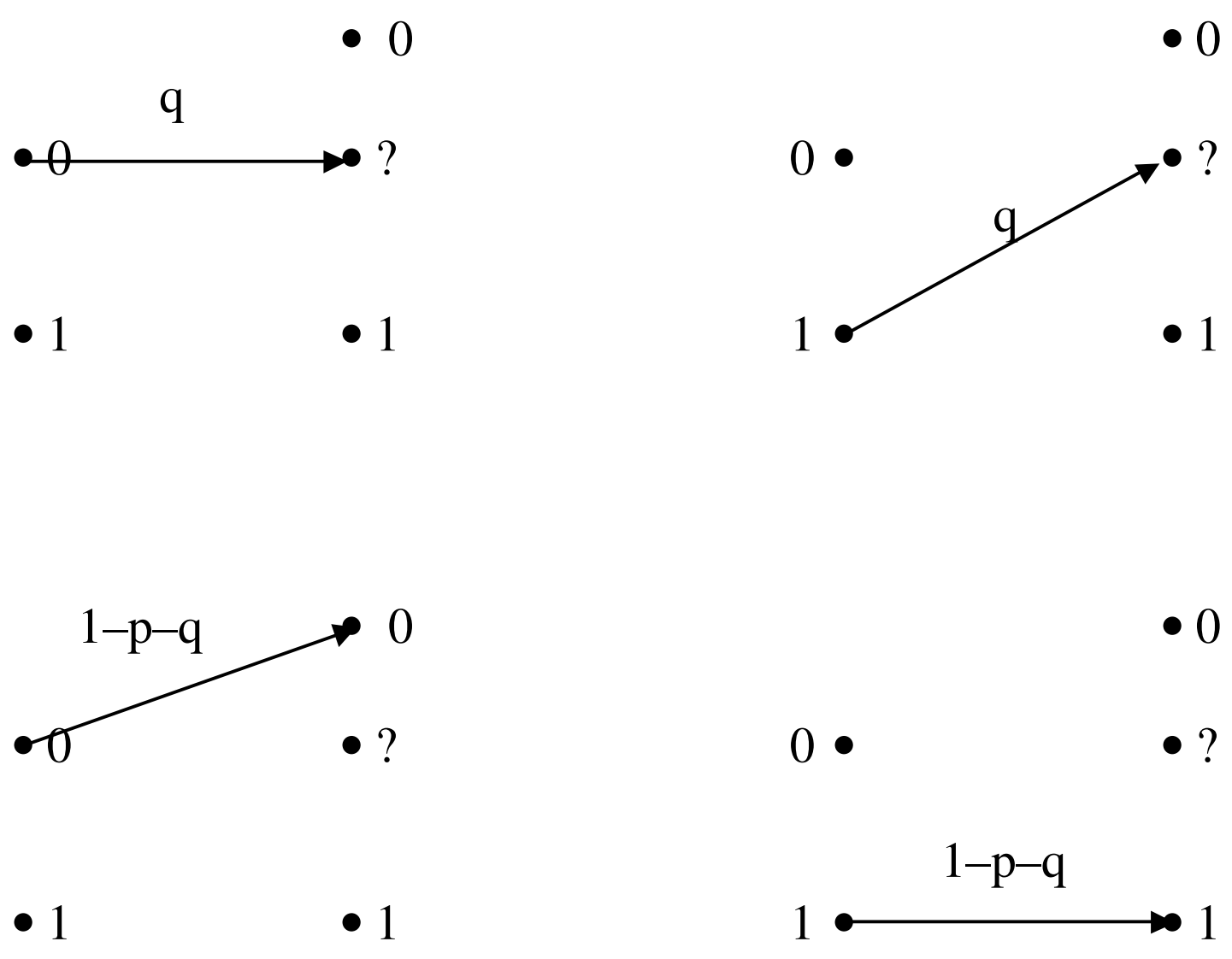

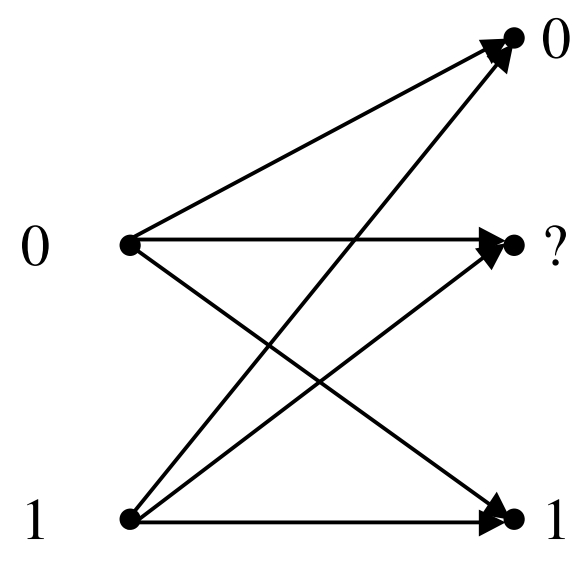

**Binary symmetric erasure channel**

Here when the received signal is very weak the demodulator does not give an output. That is the output corresponding to that particular input is erased and the blank space is left out. Here we assume when 0 is transmitted and the demodulator does not give an output and the probability of that is q when 0 is transmitted and 0 is received the probability is 1 – p – q. Similarly when 1 is transmitted and the demodulator does not give a output the probability is q. When 1 is transmitted and the modulator does not give an output the probability is p. When 1 is transmitted and 1 is received the probability is 1 – p – q.

The outputs that are erased by the demodulator are called erasures or blank spaces [6-8]. Thus when erasures are present in the received code word those coordinates received in the code word would be blank spaces.

The study of properties of erasures out weighs the study of properties of errors in a code as we are sure of the number of errors that has occurred by counting the number of blank spaces (erasures) and we are also aware of their locations as they are blank during the transmission process.

But if are to study only errors we may not completely be certain of the number of errors occurred during the transmission, for instance when a message is sent and the received message is also a code word different from the sent message, we may not be able to determine it as error but in case of erasures it may be blank.

Another advantage of erasure techniques over the study of error is that even a lay man can guess that the received message is an erroneous one. Also we can say whether the original message is retrievable or not. The study of erasures in case of Hamming metric has been widely studied by [21-2] and [6-8].

We would be defining the new notion of "blanks" which are not erasures. These 'blanks' will be known as "special blanks" and we will not be using the notion of erasure decoding we use only the error decoding technique. We will use this "special blank" notion in the last chapter of this book where we will be using them in concatenation of linear coding with Hamming metric defined on it.

Now refer Gabidulin for the notion of rank distance codes.

**DEFINITION 2.1:** *Let $X^n$ be a n-dimensional vector space over the field $GF(2^N)$. Let $u_1, u_2, \ldots, u_n$ be a fixed basis for $X^n$ over $GF(2^N)$. Then any element $x \in X^n$ can be represented as a n-tuple $(x_1, x_2, \ldots, x_n)$ where $x_i \in GF(2^N)$ $(1 \leq i \leq n)$.*

*$GF(2^N)$ is a vector space of dimension N over GF(2). Let $v_1, v_2, \ldots, v_N$ be a fixed basis for $GF(2^N)$ over GF(2). Then any element $x_i \in GF(2^N)$ can be uniquely represented in the form of a N-tuple $(m_{1i}, m_{2i}, \ldots, m_{Ni})$, $M_N^n$ denote the ensemble of all $(N \times n)$ matrices with elements from GF(2).*

*Consider the bijection $M : X^n \to \mathrm{M}_\mathrm{N}^\mathrm{n}$ defined by the following condition for any vector $x = (x_1, x_2, \ldots, x_n) \in X^n$; the associated matrix;*

$$M(x) = \begin{bmatrix} m_{11} & m_{12} & \ldots & m_{1n} \\ m_{21} & m_{22} & \ldots & m_{2n} \\ \vdots & \vdots & & \vdots \\ m_{N1} & m_{N2} & \ldots & m_{Nn} \end{bmatrix}$$

*where the $i^{th}$ column represents the $i^{th}$ coordinate $x_i$ of x over GF(2).*

*The rank of a vector $x \in X^n$ over GF(2) is defined as the rank of the matrix M(x) over GF(2). Let r(x) denote the rank of the vector $x \in X^n$ over GF(2). By the properties of the rank of a matrix the mapping $x \to r(x)$ defines a norm on $X^n$; called the rank norm.*

*Let $X^n$ be a vector space of dimension n over $GF(2^n)$ equipped with the rank norm. Clearly the rank norm induces a metric defined as the rank metric (rank distance) on $X^n$ and is denoted by $d_R$. For $x, y \in X^n$, the rank distance between x and y is $d_R(x, y) = r(x-y)$. A vector space $X^n$ over $GF(2^n)$ such that $n \leq N$ equipped with the rank metric $d_R$ is defined as a rank distance space. So if $X^n$ is a rank distance space, a linear (n, k)*

*rank distance code is a linear subspace of dimension k in the rank distance space $X^n$ and is denoted by C.*

*A generator matrix G of C is a $k \times n$ matrix with entries from $GF(2^N)$ whose rows form a basis for C.*

*Then a $(n - k) \times n$ matrix H with entries from $GF(2^N)$ such that $GH^T = (0)$ is called the parity check matrix of C where (0) denotes a $k \times n - k$ zero matrix.*

Suppose C is a (n, k) rank distance code with generator matrix G and parity check matrix H then C is the row space of G or the null space of H. We have minimum distance of C defined as $d = \min \{r(x - y) : x, y \in C; x \neq y\}$. C is a k - dimensional subspace of the rank distance space $X^n$, if $x, y \in C$ then $x - y \in C$.

Hence d the minimum distance,

$d = \min \{r(x) \mid x \in C; x \neq 0\}$. The notion of maximum rank and the erasure techniques would be studied in the following chapters.

## Chapter Three

# ERASURE DECODING OF MAXIMUM RANK DISTANCE CODES

This chapter has three sections. Section one is introductory in nature. In section two we describe the class of MRD codes as given by E.M. Gabidulin [9]. In section three we discuss the systematic guessing process or filling up of the blank spaces in case of erasures for MRD codes [17, 29, 32].

### 3.1 Introduction

Algebraic coding theory is required in communication systems to combat the errors that occur during transmission. In many communication systems, it is often convenient to represent the set of signals to be transmitted as a higher dimensional Galois field. There are many reasons to do so. One is that it makes it possible to visualize the signals by means of vectors, which in turn has the advantage of recognizing the relationship among various types of signals that is to be considered. Secondly, the length of the message will be very much reduced resulting an increase in the rate of transmission. Here the set of basic signals will be represented by a prime Galois field GF(p) and all the possible linear combinations of the basic signals will be represented by a higher dimensional Galois field, GF($p^n$).

When the code symbols are from a higher dimensional Galois field the function of the modulator is to match the encoder output to the transmission channel. A definition of a channel generally includes the modulator, the demodulator and all the intervening transmission equipment and media. In this model in certain situations, for example when the received signal is very weak the demodulator does not give an output. That is, the output corresponding to that particular input is erased and a black space of left out.

The outputs that are erased by the demodulator are called erasures or blank spaces. Hence the events in which the demodulator does not give a output when the evidence does not clearly indicate one signal as the most probable are called erasures. Hence when erasures are present in the received vector those coordinates in the received vector will be blank spaces. Erasure decoding in case of Hamming metric has been widely studied by W.W. Peterson [21-2], David Forney [6-8]. Here we obtain a method of erasure decoding for the class of MRD codes when the minimum distance is same as the length of the code. The error correcting capability of the MRD code depends on the minimum distance and greater the minimum distance greater the error correcting capability. In this chapter the minimum distance of the MRD code is equal to its length.

By making use of this systematic guessing process an erasure or blank space can be regarded as an error which can either be detected or corrected by making use of the decoding algorithm for MRD codes. We have proved that a MRD code of length n, dimension 1 and minimum distance $n = 2t + 1$ can correct atmost t erasures and detect more than t erasures. We have obtained the number of ways in which a particular erasure can be chosen during the guessing process and we have established that the result is unaffected by various choices for the erasures available during the guessing process at each stage.

Also the method of erasure decoding is illustrated through an example.

## 3. 2 Maximum Rank Distance Codes

Maximum Rank Distance (MRD) codes are a class of codes which are analogs of generalized Reed-Solomon codes [26]. MRD codes are codes of length $n \leq N$ defined over $GF(2^N)$ equipped with the rank metric.

Suppose $X^n$ is a n-dimensional vector space over the field $GF(2^N)$. Let $u_1, u_2, \ldots, u_n$ be a fixed basis for $X^n$ over $GF(2^N)$. Then any element $x \in X^n$ can be represented as an n-tuple $(x_1, x_2, \ldots, x_n)$ where $x_i \in GF(2^N)$.

$GF(2^N)$ is a vector space of dimension N over GF(2). Let $v_1, v_2, \ldots, v_N$ be a fixed basis for $GF(2^N)$ over GF(2). Then any element $x_i \in GF(2^N)$ can be uniquely represented in the form of a N-tuple $(a_{1i}, a_{2i}, \ldots, a_{Ni})$. Let $A_N^n$ denote the ensemble of all $(N \times n)$ matrices with elements from GF(2).

Consider the bijection $A : X^n \rightarrow A_N^n$ defined by the following rule:
For any vector $x = (x_1, x_2, \ldots, x_n) \in X^n$ the associated matrix

$$A(x) = \begin{bmatrix} a_{11} & a_{12} & \cdots & a_{1n} \\ a_{21} & a_{22} & \cdots & a_{2n} \\ \vdots & \vdots & & \vdots \\ a_{N1} & a_{N2} & \cdots & a_{Nn} \end{bmatrix} \qquad (3.2.1)$$

where the $i^{th}$ column represents the $i^{th}$ coordinate '$x_i$' of 'x' over GF(2). We recall some of the basic definition from [9, 29, 30, 32].

**DEFINITION 3.2.1:** *The rank of a vector $x \in X^n$ over GF(2) is defined as the rank of the matrix A(x) over GF(2). In other words the rank of a vector $x \in X^n$ is the maximum number of* *columns that are linearly independent over GF(2) in the associated matrix A(x) of the vector x. Let r(x) denote the rank of the vector $x \in X^n$ over GF(2). By the properties of the rank of a matrix the mapping $x \mapsto r(x)$ specifies a norm on $X^n$ and is called the rank norm.*

**DEFINITION 3.2.2:** *Let $X^n$ be a vector space of dimension n over $GF(2^N)$ equipped with rank norm. The rank norm induces a metric defined as the rank metric (rank distance) on $X^n$ and is denoted by $d_R$. For $x, y \in X^n$ the rank distance between x and y is $d_R(x, y) = r(x-y)$.*

**DEFINITION 3.2.3:** *A vector space $X^n$ over $GF(2^N)$ such that $n < N$ equipped with the rank metric $d_R$ is defined as a rank distance space.*

**DEFINITION 3.2.4:** *Let $X^n$ be the rank distance space. A linear (n, k) rank distance code is a linear subspace of dimension k in the rank distance space $X^n$ and is denoted by C.*

**DEFINITION 3.2.5:** *Let C be a linear (n, k) rank distance code. A generator matrix G of C is a $k \times n$ matrix with entries from $GF(2^N)$ whose rows form a basis for C.*

**DEFINITION 3.2.6:** *Let C be a linear (n, k) rank distance code with generator matrix G. Then a $(n-k) \times n$ matrix H with entries from $GF(2^N)$ such that $GH^T = (0)$ is called the parity check matrix of C where (0) denotes the $k \times (n-k)$ zero matrix.*

Suppose C is a linear (n, k) rank distance code with generator matrix G and parity check matrix H, then C can be thought of as :

1. the row space of G or
2. the null space of H.

**DEFINITION 3.2.7:** *Let C be a linear (n, k) rank distance code. The minimum distance of C is defined as d = min {r (x – y) : x, y*

*$\in$ C, x $\neq$ y}. Since C is linear space, d is also equal to min {r(x) : x $\in$ C, x $\neq$ 0}.*

A linear (n, k) rank distance code C with minimum distance d satisfies the following bound.

Singleton-style bound[12, 16]: A linear (n, k) rank distance code C with minimum distance d satisfies the inequality $d \leq n-k+1$.

**DEFINITION 3.2.8:** *Rank distance codes which attain equality in the singleton style bound are called Maximum Rank Distance codes (MRD codes).*

MRD codes are analogs of generalized Reed-Solomon codes and can be defined through generator and parity check matrices. A MRD code with length n = N can be defined as follows:

Let [i] = 2i, i = 0, ± 1, ± 2, …

Let $h_i \in GF(2^N)$, i = 1, 2, …, n be linearly independent over GF(2).

For a given design distance $d \leq n$ let us generate the matrix

$$H = \begin{bmatrix} h_1 & h_2 & \dots & h_n \\ h_1^{[1]} & h_2^{[1]} & \dots & h_n^{[1]} \\ \vdots & \vdots & & \vdots \\ h_1^{[d-2]} & h_2^{[d-2]} & \dots & h_n^{[d-2]} \end{bmatrix} \qquad (3.2.2)$$

The linear (n, k) rank distance code with parity check matrix H given in equation (3.2.2) is an MRD code of length n and minimum distance d. We denote a (n, k, d) MRD code as C[n, k].

The encoding and decoding algorithm for MRD codes are given by Gabidulin [9].

### 3.3 Erasure Decoding of MRD codes

It has long been recognized that there are advantages in allowing the demodulator not to guess at all on certain transmission (for example when the received signal is very weak) when the evidence does not clearly indicate one signal as the most probable: such events are called erasures. In the event of an erasure or blank space it is convenient to passes on the side information to the decoder that this guess is not completely reliable. By the guessing process or filling up of the blank space an erasure or blank space can be regarded as an error which can be either detected or corrected by making use of the decoding algorithm for MRD codes.

The guessing process or filling up of the blank spaces can be done systematically by making use of the decoding properties of the MRD codes and this procedure is described below.

The C[n, 1] be a (n, 1, n) MRD code defined over $GF(2^n)$, where $n = 2t + 1$. Since the minimum distance of this code is n the associated matrix of every code word in this code has n columns to be linearly independent over GF(2).

Suppose $x = (x_1, x_2, \ldots, x_n)$ is the codeword that is transmitted. Let $y = (y_1, y_2, \ldots, y_n)$ be the received vector where t coordinates are erasures. Let $*_1, *_2, \ldots *_t$ denote the t erasures. The t erasures can occur either in an array or can occur intermittently with the $y_i$'s or it can be randomly placed in between the $y_i$'s. Still the method adopted by us will not affect the result. Without loss of generality we can assume that $*_1, *_2, \ldots, *_t$ be the first t coordinates of y. Let $y_{t+1}, y_{t+2}, \ldots y_n$ be the rest of the n – t coordinates. The choosing of each of the $*_j$'s for $j = 1, 2, \ldots, t$ is detailed below.

Now choose $*_1$ ($*_1 \neq 0$) such that the set $\{*_1, y_{t+1}, y_{t+2}, \ldots, y_n\}$ is a linearly independent set over GF(2). Having chosen $*_1$, now choose $*_2$ ($*_2 \neq 0$) such that the set $(*_2, *_1, y_{t+1}, y_{t+2}, \ldots, y_n)$ is a linearly independent set over GF(2). Having chosen $*_1, *_2,$

…, $*_s$, s < t, $*_{s+1}$ is chosen such that $*_{s+1} \neq 0$ and the set $\{*_{s+1}, *_1, *_2, \ldots, *_s, y_{t+1}, y_{t+2}, \ldots, y_n\}$ is linearly independent over GF(2).

By the above guessing process or filling up of the blank spaces the t erasures $*_1, *_2, \ldots, *_t$ in the received vector y are filled by the t probable errors. Let the vector $y' = (y'_1, y'_2, \ldots, y'_n)$ be the vector obtained from y after converting the t erasures or t blank spaces into t errors by using the guessing process. Now we apply the decoding algorithm of the MRD codes for decoding $y'$ and obtain the correct message transmitted.

The following theorem gives the number of ways in which the erasures or blank spaces can be chosen by using the guessing process or filling up of the blank spaces.

**THEOREM 3.3.1:** *Let C[n,1] be a (n, 1, n) MRD code defined over $GF(2^n)$ where n=2t + 1. Let $x = (x_1, x_2, \ldots, x_n)$ be the transmitted code word. Let $y = (y_1, y_2, \ldots, y_n)$ be the received vector where t coordinates are erasures or blank spaces. Let $*_1, *_2, \ldots, *_t$ denote the t erasures or t blank spaces. Then $*_{s+1}$, where $s + 1 \leq t$ can be chosen in $2^{n-t+2}\,[2^{t-s} - 1]$ ways by using the guessing process.*

***Proof:*** Let C[n, 1] be a MRD (n, 1, n) code defined over $GF(2^n)$ where n = 2t + 1. Suppose $x = (x_1, x_2, \ldots, x_n)$ is the codeword that is transmitted. Let $y = (y_1, y_2, \ldots, y_n)$ be the received vector where t coordinates are erasures or blank spaces.

Let $*_1, *_2, \ldots, *_t$ denote the t erasures or blank spaces. Without loss of generality let $*_1, *_2, \ldots, *_t$ be the first coordinates of the vector y. Then when applying the guessing process $*_1$ is chosen such that $*_1 \neq 0$ and the set $\{*_1, y_{t+1}, y_{t+2}, \ldots, y_n\}$ is a linearly independent set over GF(2). That is, $*_1$ is chosen from the set $B_1 = GF(2^n) \setminus \{0, y_{t+1}, \ldots, y_n, y_{t+1} + y_{t+2}, \ldots, y_{n-1} + y_n, y_{t+1} + y_{t+2} + y_{t+3}, \ldots, y_{n-2} + y_{n-1} + y_n, \ldots, y_{t+1} + y_{t+2} + \ldots + y_{n-1} + y_n\}$.

Hence $*_1$ can be chosen in $|B_1|$ ways, where $|B_1|$ denotes the number of elements in the set $B_1$. Clearly

$$
\begin{aligned}
|B_1| &= 2^n - (1 + n - {}_tC_1 + n - {}_tC_2 + n - {}_tC_3 \\
&\quad + \ldots + n - {}_tC_{n-t-1} + n - {}_tC_{n-1}) \\
&= 2^n - (1+1)^{n-1} \qquad \text{(by Binomial theorem)} \\
&= 2^n - 2^{n-t} \\
&= 2^{n-t}\ [2^t - 1].
\end{aligned}
$$

Having chosen $*_1$, $*_2$ can be chosen such that $*_2 \neq 0$ and the set $\{*_2, *_1, y_{t+1}, y_{t+2}, \ldots, y_n\}$ is linearly independent over GF(2). That is, $*_2$ is chosen from the set

$B_2 = GF(2^n) \setminus \{0, *_1, y_{t+1}, \ldots, y_n, *_1 + y_{t+1}, y_{t+1} + y_{t+2}, \ldots, y_{n-1} + y_n, *_1 + y_{t+1} + y_{t+2}, y_{t+1} + y_{t+2} + y_{t+3}, \ldots, y_{n-2} + y_{n-1} + y_n, \ldots, *_1 + y_{t+1} + y_{t+2} + \ldots + y_{n-1} + y_n\}$. Hence $*_2$ can be chosen in $|B_2|$ ways. Clearly $|B_2|$ is given by;

$$
\begin{aligned}
|B_2| &= 2^n - (1 + n - {}_{(t+1)}C_1 + n - {}_{(t+1)}C_2 \\
&\quad + \ldots + n - {}_{(t+1)}C_{n-t} + n - {}_{(t+1)}C_{n-t+1}) \\
&= 2^n - (1+1)^{n-t+1} \qquad \text{(By binomial theorem)} \\
&= 2^n - 2^{n-t+1} \\
&= 2^{n-t+1}\ [2^{t-1} - 1].
\end{aligned}
$$

Having chosen $*_1, *_2, \ldots, *_s$, $s < t$, $*s+1$ is chosen such that $*_{s+1} \neq 0$ and $\{*_{s+1}, *_1, *_s, \ldots, *_s, y_{t+1}, y_{t+2}, \ldots, y_n\}$ is linearly independent set over GF(2). Hence $*_{s+1}$ can be chosen in $|B_{s+1}|$ ways. Clearly $|B_{s+1}|$ is

$$
\begin{aligned}
|B_{s+1}| &= 2^n - (1 + n - {}_{(t+s)}C_1 + n - {}_{(t+s)}C_2 + n - {}_{(t+s)}C_3 \\
&\quad + \ldots + n - {}_{(t+s)}C_{n-t+s-1} + n - {}_{(t+s)}C_{n-t+s}) \\
&= 2^n - (1+1)^{n-t+s} \\
&= 2^n - 2^{n-t+s} \\
&= 2^{n-t+s}\ [2^{t-s} - 1].
\end{aligned}
$$

Hence the theorem.

**THEOREM 3.3.2:** *Let C[n, 1] be a (n, 1, n) MRD code with minimum distance n = 2t + 1. Then C[n, 1] corrects atmost t ensures and detects more than t erasures.*

***Proof:*** Let C [n, 1] be a (n, 1, n) MRD code with minimum distance n = 2t + 1. Suppose x = $(x_1, x_2, \ldots, x_n)$ is the codeword that is transmitted. Let y = $(y_1, y_2, \ldots, y_n)$ be the received vector where t coordinates are erasures or blank spaces. Let y′ = $(y'_1, y'_2, \ldots, y'_n)$ be the vector obtained from y after converting the t erasures into t errors by using the guessing process or filling up of the blank spaces. Therefore, $d_R$ (x, y′) < t, by [9, 29, 30] we know that a MRD code with n = 2t + 1 can correct atmost t errors and detect more than t errors. Hence C[n, 1] can correct atmost t erasures and detect more than t erasures.

**THEOREM 3.3.3:** *Let C[n, 1] be a (n, 1, n) MRD code defined over GF($2^n$) where n = 2t + 1. Let x = $(x_1, x_2, \ldots, x_n)$ be the transmitted codeword. Let y = $(y_1, y_2, \ldots, y_n)$ be the received vector where t coordinates are erasures. Let $*_1, *_2, \ldots, *_t$ denote the t erasures. The erasures $*_{s+1}$, where $s + 1 \leq t$ can be chosen in $2^{n-t+s}$ $[2^{t-s}-1]$ ways by using the guessing process. Then the various choices of $*_j$'s does not affect the erasure correcting capability of the MRD code.*

***Proof :*** Let C[n, 1] be a (n, 1, n) MRD code defined over GF($2^n$) where n = 2t + 1. Suppose x = $(x_1, x_2, \ldots, x_n)$ be the transmitted codeword. Let y = $(y_1, y_2, \ldots, y_n)$ be the received vector where t coordinates are erasures or blank spaces. Let $*_1$, $*_2, \ldots, *_t$ denote the t erasures or t blank spaces. Let y′ = $(y'_1, y'_2, \ldots, y'_n)$ be the vector obtained from y after applying the guessing process or filling up of the blank spaces to obtain from y after applying the guessing process or filling up of the blank spaces to the received vector y. Let y″ = $(y''_1, y''_2, \ldots, y''_n)$ be another vector obtained from y after applying the guessing process or filling up of the blank spaces to the received vector y and be such that y′ ≠ y″, then $d_R$(y′, y″) ≤ t. By the guessing process the t erasures in the received vector y are converted into t errors. Hence $d_R$(x, y′) ≤ t and $d_R$(x, y″) ≤ t. We know MRD code with n = 2t + 1 can correct atmost t errors and detect more than t errors. Since $d_R(x, y') \le t$ and $d_R(x, y'') \le t$ both the vector $y''$ and $y'$ are both decoded as the vector x, which is the transmitted message. Hence the theorem.

The guessing process and the above results are illustrated by the following example.

***Example 3.3.1:*** Consider $GF(2^3)$ where $GF(2^3) = \{0, 1, \alpha, \alpha^2, \ldots, \alpha^6\}$ and $\alpha$ is the root of the primitive polynomial $p(x) = x^3 + x + 1$. Consider the (3, 1, 3) MRD code having parity check matrix H.

$$H = \begin{bmatrix} 1 & \alpha & \alpha^2 \\ 1 & \alpha^2 & \alpha^4 \end{bmatrix} \tag{3.3.1}$$

obtained by choosing $h_1 = 1$, $h_2 = \alpha$, $h_3 = \alpha^2$ in the matrix (3.3.1).

Let $m_1 = \alpha^5$ be the message. Let $x = (\alpha^5, a_1, a_2)$. The $Hx^T = (0)$ gives $a_1 = \alpha^6$ and $a_2 = \alpha^2$. The codeword corresponding to the message $m_1 = \alpha^5$ is $c = (\alpha^5, \alpha^6, \alpha^2)$. Suppose c is the transmitted vector and let $y = (\alpha^5, *, \alpha^2)$ be the received vector with one erasure or one blank space. By applying the guessing process we choose $* \neq 0$ and such that the set $\{*, \alpha^5, \alpha^2\}$ is linearly independent over GF(2). Then there are four choices for * namely $* = 1$, $* = \alpha$, $* = \alpha^4$ and $* = \alpha^6$. We see that the choice of * does not change the result and we get the corrected message for all the four choices of *.

Choose * = 1.

Now the vector $y'$ obtained after filling up of the blank space in the received vector y is $y' = (\alpha^5, 1, \alpha^2)$. By applying the decoding algorithm for MRD codes to the vector $y'$ we have

1. Syndrome = $(\alpha^3, \alpha^4) = (s_0, s_1)$.

2. $\Delta(z) = \alpha^2 z + z^2$.

3. The non zero root of $\Delta(z)$ is $E_1 = \alpha^2$.

4. Solving $s_p = \sum_{j=1}^{1} E_j x_j^{[p]}$, p = 0 we get $x_1 = \alpha$.

5. Solving $x_1 = \sum_{j=1}^{3} Y_{1j} h_j$ gives $Y_{11} = 0$, $Y_{12} = 1$, $Y_{13} = 0$ and hence the matrix Y = (0 1 0).

6. Error vector $e = E_1 Y = (0, \alpha^2, 0)$.

7. $y' + e = (\alpha^5, \alpha^6, \alpha^2)$ is the required codeword.

Choose * = $\alpha$.

Now the vector $y'$ obtained after filling up of the blank space in the received vector y is $y' = (\alpha^5, \alpha, \alpha^2)$. By applying the decoding algorithm for MRD codes to the vector $y'$ we have

1 Syndrome = $(\alpha^6, 1) = (s_0, s_1)$.

2. $\Delta(z) = \alpha^5 z + z^2$.

3. The non zero root of $\Delta(z)$ is $E_1 = \alpha^5$.

4. Solving $s_p = \sum_{j=1}^{1} E_j x_j^{[p]}$, p = 0 we get $x_1 = \alpha$.

5. Solving $x_1 = \sum_{j=1}^{3} Y_{1j} h_j$ gives $Y_{11} = 0$, $Y_{12} = 1$, $Y_{13} = 0$ and hence the matrix Y = (0 1 0).

6. Error vector $e = E_1 Y = (0, \alpha^5, 0)$.

7. $y' + e = (\alpha^5, \alpha^6, \alpha^2)$ is the required codeword.

Choose $* = \alpha^4$.

Now the vector $y'$ obtained after filling up of the blank space in the received vector y is $y' = (\alpha^5, \alpha^4, \alpha^2)$. By applying the decoding algorithm for MRD codes to the vector $y'$ we have

Syndrome = (0, 0).

This indicates that $y'$ is the correct message.

Let $m_2 = \alpha^3$ be the message. Let $x = (\alpha^3, a_1, a_2)$. The $Hx^T = (0)$ gives $a_1 = \alpha^4$ and $a_2 = 1$. The codeword corresponding to the message $m_1 = \alpha^3$ is $c = (\alpha^3, \alpha^4, 1)$. Suppose c is the transmitted vector and let $y = (\alpha^3, *_1, *_2)$ be the received vector with two erasures or blank spaces. By applying the guessing process we choose $*_1 \neq 0$ and such that the set $\{*_1, \alpha^3\}$ is linearly independent over GF(2). Since the set $\{\alpha^6, \alpha^3\}$ is linearly independent over GF(2) choose $*_1 = \alpha^6$. Again since the set $\{\alpha^2, \alpha^6, \alpha^3\}$ is a linearly independent set over GF(2) choose $*_2 = \alpha^2$. Now the vector $y'$ obtained after filling up of the blank spaces in the received vector y is $y' = \{\alpha^3, \alpha^6, \alpha^2\}$.

By applying the decoding algorithm for MRD codes to the vector $y'$ we have

1. Syndrome = $(\alpha^4, \alpha^5) = (s_0, s_1)$.

2. $\Delta(z) = \alpha^3 z + z^2$.

3. The nonzero root of $\Delta(z)$ is $E_1 = \alpha^3$.

4. Solving $s_p = \sum_{j=1}^{1} E_j x_j^{[p]}$, $p = 0$ we get $x_1 = \alpha$.

5. Solving $x_1 = \sum_{j=1}^{3} Y_{1j} h_j$ gives $Y_{11} = 0$, $Y_{12} = 1$, $Y_{13} = 0$ and hence the matrix Y = (0 1 0).

6. Error vector $e = E_1 Y = (0, \alpha^3, 0)$.

8. $y' + e = (\alpha^5, \alpha^6, \alpha^2)$ is the required codeword.

Choose $* = \alpha^6$.

Now the vector $y'$ obtained after filling up of the blank space in the received vector y is $y' = (\alpha^5, \alpha^6, \alpha^2)$. By applying the decoding algorithm for MRD codes to the vector $y'$ we have

Syndrome = (0, 0).

This indicates that $y'$ is the correct message.

Let $m_2 = \alpha^3$ be the message. Let $x = (\alpha^3, a_1, a_2)$. The $Hx^T = (0)$ gives $a_1 = \alpha^4$ and $a_2 = 1$ . The codeword corresponding to the message $m_1 = \alpha^3$ is $c = (\alpha^3, \alpha^4, 1)$. Suppose c is the transmitted vector and let $y = (\alpha^3, *_1, *_2)$ be the received vector with two erasures or two blank space. By applying the guessing process we choose $*_1 \neq 0$ and such that the set $\{*_1, \alpha^3\}$ is linearly independent over GF(2). Since the set $\{\alpha^6, \alpha^3\}$ is linearly independent over GF(2) choose $*_1 = \alpha^6$. Again since the set $\{\alpha^2, \alpha^6, \alpha^3\}$ is a linearly independent set over GF(2) choose $*_2 = \alpha^2$. Now the vector $y'$ obtained after filling up of the blank spaces in the received vector y is $y' = \{\alpha^3, \alpha^6, \alpha^2\}$.

By applying the decoding algorithm for MRD codes to the vector $y'$ we have

1. Syndrome $= (\alpha^2, \alpha^2) = (s_0, s_1)$.

2. $\Delta(z) = \alpha^2 z + z^2$.

3. The nonzero root of $\Delta$ (z) is $E_1 = \alpha^2$.

4. Solving $s_p = \sum_{j=1}^{1} E_j x_j^{[p]}$, p = 0 we get $x_1 = 1$.

5. Solving $x_1 = \sum_{j=1}^{3} Y_{1j} h_j$ gives $Y_{11} = 1$, $Y_{12} = 0$, $Y_{13} = 0$ and hence the matrix Y = (1 0 0).

Error detected.

The erasure decoding method declares that more than one erasure has been detected.

Hence the (3, 1, 3) MRD code corrects atmost one erasure and detects more than one erasure in the received vector. We also note that the result is independent of the choice of *.

The study of erasure decoding is better than the study of errors. Erasure decoding in case of Hamming codes has been widely studied in literature. In this chapter we have started the study of erasure decoding of the class of MRD codes. We give a guessing process by which the erasures or blank spaces in the received vector can be converted into errors. We have proved that the MRD code of length n, dimension one and minimum distance n = 2t + 1 can correct atmost t erasures and detect more than t erasures. We have obtained the number of ways in which a particular erasure can be chosen during the guessing process and we have established that the result is unaffected by the various choices of erasures or blank spaces that are available in the guessing process or filling up of the blank spaces at each stage.

It is pertinent to mention that these results are true for any (n, k, d) MRD code defined over $GF(2^n)$ where d = 2t + 1.

Chapter Four

# MRD CODES – SOME PROPERTIES AND A DECODING TECHNIQUE

## 4.1 Introduction

Minimum distance is a one of the chief parameters which determines the error-correcting capability of a code. In fact, the maximum number of errors corrected by a code is proportional to its minimum distance. For a linear code of length n, dimension k and minimum distance d, the upper bound for the minimum distance is n – k + 1. A code which has the minimum distance d = n – k + 1 is defined as the Maximum Distance Separable codes. The Reed-Solomon codes are important Maximum Distance Separable codes. Analogues to Maximum Distance Separable codes, Gabidulin in [9] defines the class of Maximum Rank Distance (MRD) codes. An $[n,k,d]_{q^N}$ Rank Distance code whose minimum distance d is equal to n – k + 1 is called an MRD code, n ≤ N. As the error-correcting capability of a code is proportional to its minimum distance, codes with larger minimum distance is preferred for error correction. It makes it justifiable to study the characteristics or properties enjoyed by such class of codes so as to enunciate their uses in the communication channels.

This chapter has four sections. Section one is introductory in nature. Section 2 introduces for the first time, a combined error-erasure decoding technique to the class of $[n,k,d]_{q^N}$ MRD codes, $n \leq N$. A code is said to be invertible if, knowing only the parity-check symbols of a codeword, the corresponding information symbols can be determined through an inversion process. The invertible property facilitates the data recovery process in error-control schemes. A comprehensive study on the invertible property for the class of $[n,k,d]_{q^N}$ q-Cyclic RD codes is carried out in section 3. Section 3 also presents the systematic encoding and the shortening technique for the class of $[n,k,d]_{q^N}$ q-Cyclic RD codes. Section 4 carries out a study on the class of Rank Distance codes having complementary duals. It is proved that the class of $[n,k,d]_{2^n}$ MRD codes generated by the trace-orthogonal-generator matrices are LCD codes. Further, description to the (noiseless and noisy) 2-user F-Adder Channel and coding for the noiseless 2-user F-Adder Channel via the class of $[n,k,d]_{2^n}$ MRD codes having complementary duals are presented [25]. A coding problem for the noisy 2-user F-Adder Channel is explained.

## 4.2 Error-Erasure Decoding Technique to MRD Codes

This section introduces a combined error-erasure decoding technique for the first time to the class of $[n,k,d]_{q^N}$ MRD codes, $n \leq N$. The combined error-erasure decoding technique to the class of MRD codes presented later in this section enables the decoder to correct all combinations of r rank-errors and s erasures in an erroneously received vector as long as $2r + s < d$, where d is the minimum-rank distance of the MRD codes.

Let $\Gamma$ denote an $[n,k,d]_{q^N}$ MRD code, $n \leq N$ with the parity-check matrix H.

$$H = \begin{bmatrix} h_1 & h_2 & \cdots & h_n \\ h_1^{[1]} & h_2^{[1]} & \cdots & h_n^{[1]} \\ \vdots & \vdots & \ddots & \vdots \\ h_1^{[d-2]} & h_2^{[d-2]} & \cdots & h_n^{[d-2]} \end{bmatrix} \qquad (4.1)$$

where $h_1, h_2, \ldots, h_n \in GF(q^N)$ are linearly independent over GF (q).

Before presenting the error-erasure decoding technique, in what follows are the descriptions to the deletion of m coordinates of $x \in \Gamma$ and the deletion of m columns of H that are required in the decoding technique, where m < n.

Let $x = (x_1, x_2, \ldots, x_n) \in \Gamma$. Define the deletion of the $l_1^{th}$ coordinate $x_{l_1}$ of x as the (n–1)-tuple $(x_1, x_2, \ldots, x_{l_1-1}, x_{l_1+1}, \ldots, x_n)$ and denote it by $x^{(1)}$. Similarly, the deletion of m coordinates, say $x_{l_1}, x_{l_2}, \ldots, x_{l_m}$ of x is defined as the (n–m)-tuple $(x_1, x_2, \ldots, x_{l_1-1}, x_{l_1+1}, \ldots, x_{l_2-1}, x_{l_2+1}, \ldots, x_{l_m-1}, x_{l_m+1}, \ldots, x_n)$ and is denoted by $x^{(m)}$.

Generate a matrix, say $H^{(1)}$ of order $(d–2) \times (n–1)$ from H by performing the following row operations so as to delete the $l_1^{th}$ column of H.

For each j = 1, 2, …, d – 2, multiply the $j^{th}$ row-vector of H by $\dfrac{h_{l_1}^{[j]}}{h_{l_1}^{[j-1]}}$ and subtract the $(j+1)^{th}$ row-vector from the $j^{th}$ row-vector and let the deletion of the $l_1^{th}$ coordinate of the resultant vector, namely ( $h_1^{[j-1]}\dfrac{h_{l_1}^{[j]}}{h_{l_1}^{[j-1]}} - h_1^{[j]},\ h_2^{[j-1]}\dfrac{h_{l_1}^{[j]}}{h_{l_1}^{[j-1]}} - h_2^{[j]}, \ldots,$

$h_{l_1-1}^{[j-1]}\dfrac{h_{l_1}^{[j]}}{h_{l_1}^{[j-1]}} - h_{l_1-1}^{[j]},\ h_{l_1+1}^{[j-1]}\dfrac{h_{l_1}^{[j]}}{h_{l_1}^{[j-1]}} - h_{l_1+1}^{[j]}, \ldots,\ h_n^{[j-1]}\dfrac{h_{l_1}^{[j]}}{h_{l_1}^{[j-1]}} - h_n^{[j]}$ )

be the j[th] row-vector in $H^{(1)}$.

The resultant matrix $H^{(1)}$ is given by

$$\left[\begin{matrix} h_1^{[0]}\dfrac{h_{l_1}^{[1]}}{h_{l_1}^{[0]}} - h_1^{[1]} & \cdots & h_{l_1-1}^{[0]}\dfrac{h_{l_1}^{[1]}}{h_{l_1}^{[0]}} - h_{l_1-1}^{[1]} \\ h_1^{[1]}\dfrac{h_{l_1}^{[2]}}{h_{l_1}^{[1]}} - h_1^{[2]} & \cdots & h_{l_1-1}^{[1]}\dfrac{h_{l_1}^{[2]}}{h_{l_1}^{[1]}} - h_{l_1-1}^{[2]} \\ \vdots & \ddots & \vdots \\ h_1^{[d-3]}\dfrac{h_{l_1}^{[d-2]}}{h_{l_1}^{[d-3]}} - h_1^{[d-2]} & \cdots & h_{l_1-1}^{[d-3]}\dfrac{h_{l_1}^{[d-2]}}{h_{l_1}^{[d-3]}} - h_{l_1-1}^{[d-2]} \end{matrix}\right.$$

$$\left.\begin{matrix} h_{l_1+1}^{[0]}\dfrac{h_{l_1}^{[1]}}{h_{l_1}^{[0]}} - h_{l_1+1}^{[1]} & \cdots & h_n^{[0]}\dfrac{h_{l_1}^{[1]}}{h_{l_1}^{[0]}} - h_n^{[1]} \\ h_{l_1+1}^{[1]}\dfrac{h_{l_1}^{[2]}}{h_{l_1}^{[0]}} - h_{l_1+1}^{[2]} & \cdots & h_n^{[1]}\dfrac{h_{l_1}^{[2]}}{h_{l_1}^{[0]}} - h_n^{[2]} \\ \vdots & \ddots & \vdots \\ h_{l_1+1}^{[d-3]}\dfrac{h_{l_1}^{[d-2]}}{h_{l_1}^{[d-3]}} - h_{l_1+1}^{[d-2]} & \cdots & h_n^{[d-3]}\dfrac{h_{l_1}^{[d-2]}}{h_{l_1}^{[1]}} - h_n^{[d-2]} \end{matrix}\right].$$

Let $h_i' = \dfrac{h_{l_1}^{[j]}}{h_{l_1}^{[0]}} - h_i^{[1]}$ for each $i \neq l_1$.

Then $h_i^{[j-1]} \dfrac{h_{l_1}^{[j]}}{h_{l_1}^{[j-1]}} - h_i^{[j]} = h_i'^{[j-1]}$ for each $i \neq l_1$ and $j = 1, 2, \ldots, d-2$.

The $(d-2) \times (n-1)$ matrix $H^{(1)}$ now takes the following form

$$H^{(1)} = \begin{bmatrix} h_1'^{[0]} & h_1'^{[0]} & \cdots & h_{l_1-1}'^{[0]} & h_{l_1+1}'^{[0]} & \cdots & h_n'^{[0]} \\ h_1'^{[1]} & h_2'^{[1]} & \cdots & h_{l_1-1}'^{[1]} & h_{l_1+1}'^{[1]} & \cdots & h_n'^{[1]} \\ \vdots & \vdots & \ddots & \vdots & \vdots & \ddots & \vdots \\ h_1'^{[d-3]} & h_2'^{[d-3]} & \cdots & h_{l_1-1}'^{[d-3]} & h_{l_1+1}'^{[d-3]} & \cdots & h_n'^{[d-3]} \end{bmatrix}$$

where $h_i' = h_i^{[0]} \dfrac{h_{l_1}^{[j]}}{h_{l_1}^{[0]}} - h_i^{[1]}$ for each $i \neq l_1$.

Call the matrix $H^{(1)}$ thus obtained as the deletion of $l_1^{th}$ column of H. Proceeding with a similar row operations in $H^{(1)}$ so as to delete a column of $H^{(1)}$ that is corresponding to the $l_2^{th}$ column of H, one obtains a $(d-3) \times (n-2)$ matrix, say $H^{(2)}$; call it as the deletion of $l_1^{th}$ and $l_2^{th}$ columns of H.

In general, proceeding with a similar row operations in $H^{(m-1)}$ so as to delete a column of $H^{(m-1)}$ that is corresponding to the $l_m^{th}$ column of H, one obtains a $(d-1-m) \times (n-m)$ matrix say $H^{(m)}$, called as the deletion of $l_1^{th}$, $l_2^{th}$, …, $l_m^{th}$ columns of H, where $H^{(m-1)}$ is the $(d-m) \times (n-m+1)$ matrix and is the deletion of $l_1^{th}$, $l_2^{th}$, …, $l_{m-1}^{th}$ columns of H.

An important relationship between $x^{(1)}$ and $H^{(1)}$ is proved in the following lemma.

***Lemma 4.2.1:*** Let Γ be an $[n,k,d]_{q^N}$ MRD code with the parity-check matrix H as defined in (4.1). Let $x^{(1)}$ be the deletion of the $l_1^{th}$ coordinate of $x \in \Gamma$. Let $H^{(1)}$ be the deletion of the $l_1^{th}$ column of H. Then $x^{(1)} \, H^{(1)^T} = (0)$.

***Proof:*** Let $x = (x_1, x_2, \ldots, x_n) \in \Gamma$. Then $x^{(1)} = (x_1, x_2, \ldots, x_{l_1-1}, x_{l_1+1}, \ldots, x_n)$ is the deletion of the $l_1^{th}$ coordinate $x_{l_1}$ of x. Let $y = x^{(1)} H^{(1)^T} = (y_1, y_2, \ldots, y_{d-2})$.

Then $y = x^{(1)} H^{(1)^T}$

$$= x^{(1)} \begin{bmatrix} h_1'^{[0]} & h_1'^{[0]} & \ldots & h_{l_1-1}'^{[0]} & h_{l_1+1}'^{[0]} & \ldots & h_n'^{[0]} \\ h_1'^{[1]} & h_2'^{[1]} & \ldots & h_{l_1-1}'^{[1]} & h_{l_1+1}'^{[1]} & \ldots & h_n'^{[1]} \\ \vdots & \vdots & \ddots & \vdots & \vdots & \ddots & \vdots \\ h_1'^{[d-3]} & h_2'^{[d-3]} & \ldots & h_{l_1-1}'^{[d-3]} & h_{l_1+1}'^{[d-3]} & \ldots & h_n'^{[d-3]} \end{bmatrix}^T$$

where $h_i'^{[0]} = h_i^{[0]} \dfrac{h_{l_1}^{[1]}}{h_{l_1}^{[0]}} - h_i^{[1]}$ for each $i \neq l_1$.

One needs to show that

$$x^{(1)} \left( h_1^{[j-1]} \frac{h_{l_1}^{[j]}}{h_{l_1}^{[j-1]}} - h_1^{[j]},\ h_2^{[j-1]} \frac{h_{l_1}^{[j]}}{h_{l_1}^{[j-1]}} - h_2^{[j]}, \ldots, \right.$$

$$\left. h_{l_{1-1}}^{[j-1]} \frac{h_{l_1}^{[j]}}{h_{l_1}^{[j-1]}} - h_{l_{1-1}}^{[j]},\ h_{l_1+1}^{[j-1]} \frac{h_{l_1}^{[j]}}{h_{l_1}^{[j-1]}} - h_{l_1+1}^{[j]}, \ldots, h_n^{[j-1]} \frac{h_{l_1}^{[j]}}{h_{l_1}^{[j-1]}} - h_n^{[j]} \right)^T$$

= 0, for each j = 1, 2, …, d – 2.

Now $y_j = x^{(1)} ( h_1'^{[j-1]}, h_2'^{[j-1]}, \ldots, h_{l_1-1}'^{[j-1]}, h_{l_1+1}'^{[j-1]}, \ldots, h_n^{[j-1]})^T$

$$= x_1 h_1'^{[j-1]} + x_2 h_2'^{[j-1]} + \ldots$$

$$+ x_{l_{1-1}} h_{l_1-1}'^{[j-1]} + x_{l_{1+1}} h_{l_1+1}'^{[j-1]}, \ldots, x_n h_n'^{[j-1]}$$

$$= \begin{cases} x_1 h_1^{[j-1]}\left(\dfrac{h_{l_1}^{[j]}}{h_{l_1}^{[j-1]}}\right) + \ldots + x_{l_{1}-1} h_{l_{1}-1}^{[j-1]}\left(\dfrac{h_{l_1}^{[j]}}{h_{l_1}^{[j-1]}}\right) + x_{l_{1}+1} h_{l_{1}+1}^{[j-1]}\left(\dfrac{h_{l_1}^{[j]}}{h_{l_1}^{[j-1]}}\right) + \ldots + \\ x_n h_n^{[j-1]}\left(\dfrac{h_{l_1}^{[j]}}{h_{l_1}^{[j-1]}}\right) - (x_{l_1} h_1^{[j]} + \ldots + x_{l_1 - 1} h_{l_1-1}^{[j]} + x_{l_1+1} h_{l_1+1}^{[j]} + \ldots + x_n h_n^{[j]} \end{cases}$$

$$= -x_1\ h_{l_1}^{[j-1]} \frac{h_{l_1}^{[j]}}{h_{l_1}^{[j-1]}} - (-x_1\ h_{l_1}^{[j]}) \qquad (\because\ xH^T = (0))$$

$$= -x_1\ h_{l_1}^{[j]} + (x_1\ h_{l_1}^{[j]}) = 0$$

i.e. $y_j = 0$, for each j = 1, 2, …, d–2.

Thus $x^{(1)}\ H^{(1)^T} = (0)$. Hence the lemma.

The above lemma is true for the general case also. That is, if $x^{(m)}$ is the deletion of m coordinates, say $x_{l_1}, x_{l_2}, ..., x_{l_m}$ of $x \in \Gamma$ and $H^{(m)}$ is the deletion of $l_1^{th}$, $l_2^{th}$, …, $l_m^{th}$ columns of H, then $x^{(m)}\ H^{(m)^T} = (0)$.

Thus the proof of the following lemma follows immediately.

***Lemma 4.2.2:*** Let Γ be an $[n,k,d]_{q^N}$ MRD code with the parity-check matrix H as defined in (4.1). Let $x^{(m)}$ be the deletion of m coordinates $x_{l_1}, x_{l_2}, ..., x_{l_m}$ of $x \in \Gamma$ and $H^{(m)}$ be the deletion of $l_1^{th}$, $l_2^{th}$, …, $l_m^{th}$ columns of H, m < n.

Then $x^{(m)}\ H^{(m)^T} = 0$.

The above lemma gives a relationship between $x^{(m)}$ and $H^{(m)}$. This relation plays a crucial role in the combined error-erasure decoding technique to the class $[n,k,d]_{q^N}$ MRD codes presented in what follows.

## Error-Erasure Decoding Technique:

Consider the parity-check matrix H of an $[n,k,d]_{q^N}$ MRD code, $n \le N$.

$$H = \begin{bmatrix} h_1 & h_2 & \cdots & h_n \\ h_1^{[1]} & h_2^{[1]} & \cdots & h_n^{[1]} \\ \vdots & \vdots & \ddots & \vdots \\ h_1^{[d-2]} & h_2^{[d-2]} & \cdots & h_n^{[d-2]} \end{bmatrix}$$

where $h_i \in GF(q^N)$, i=1, 2, …, N are linearly independent over GF(q).

Let $x = (x_1, x_2, \ldots, x_n) \in \Gamma$ be a codeword transmitted over a noisy channel. Because of the channel noise, the receiver may not receive the transmitted codeword x. Let $y = \overline{y} \dotplus e$ be the received vector, where $e = (e_1, e_2, \ldots, e_n)$ is an error-vector and $\overline{y}$ denotes the codeword x with erasures; the details of the notations '$\doteq$' and '$\dotplus$' are given in the next line. Since the erroneously received vector y also has erasures and erasures are nothing but blank spaces, '$\doteq$' and '$\dotplus$' are so used to represent the received vector in terms of the error-vector and erasure-vector.

Assume that the received vector $y = (y_1, y_2, \ldots, y_n)$ has $m \le r$ rank-errors and $t \le s$ erasures such that $2r + s < d$, where d is the minimum-rank distance of the $[n,k,d]_{q^N}$ MRD code and m is the rank of the error-vector e. Without loss of generality, assume that the received vector y has erasures in the first t coordinates; i.e., $\overline{y} = (*_1, *_2, \ldots, *_t, x_{t+1}, \ldots, x_n)$ with $*_1, *_2, \ldots, *_t$ representing erasures. One should note that the t coordinates of e that are corresponding to the erasures in the received vector y are zeros; i.e., $e_1, e_2, \ldots, e_t$ are all zero. Therefore, e = $(\underbrace{0, 0, \ldots, 0}_{t\ terms}, e_{t+1}, e_{t+2}, \ldots, e_n)$. Then the received vector takes the form $y \doteq (*_1, *_2, \ldots, *_t, x_{t+1}, \ldots, x_n) \dot{+} (0, 0, \ldots, 0, e_{t+1}, e_{t+2}, \ldots, e_n)$.

Let $f = (f_1, f_2, \ldots, f_t, \underbrace{0, 0, \ldots, 0}_{n-t\ terms})$ be the erasure-vector.

Let $\overline{y}' = (\overline{y}'_1, \overline{y}'_2, \ldots, \overline{y}'_n)$ be such that $\overline{y}'_i = \begin{cases} f_i & \text{if } x_i \text{ is erased} \\ x_i & \text{otherwise} \end{cases}$,

where the unknowns $f_1, f_2, \ldots, f_t$ $(t \leq s)$ are to be determined.

Then $\overline{y}'$ = $(f_1, f_2, \ldots, f_t, x_{t+1}, \ldots, x_n)$. Having replaced the erasures with the unknown $f_1, f_2, \ldots, f_t$, the received vector y is now expressed as $\overline{y}'$ + e; i.e., y = $\overline{y}'$ + e. The syndrome of y is given by

$$
\begin{aligned}
S &= y H^T \\
&= (\overline{y}' + e)H^T \\
&= (f_1, f_2, \ldots, f_t, x_{t+1}, \ldots, x_n)H^T \\
&\quad + (0, 0, \ldots, 0, e_{t+1}, \ldots, e_n)H^T \\
&= (f_1, f_2, \ldots, f_t, x_{t+1}, e_{t+1}, \ldots, x_n + e_n)\, H^T. \qquad (4.2)
\end{aligned}
$$

The decoder's problem is to first determine the error-vector $e = (0, 0, \ldots, 0, e_{t+1}, \ldots, e_n)$ and then the erasure-vector $f = (f_1, f_2, \ldots, f_t, 0, \ldots, 0)$ on the basis of the syndrome vector $S = (s_0, s_1, \ldots, s_{d-2})$, where $s_i$ denotes the $i^{th}$ coordinates of S for each i = 0, 1, … d–2.

Equating each component in the right-hand side of (4.2),

$$f_1 h_1 + \ldots + f_t h_t + (x_{t+1} + e_{t+1})\, h_{t+1} + \ldots + (x_n + e_n)\, h_n = 0$$

$$f_1\ h_1^{[1]} + \ldots + f_t\ h_t^{[1]} + (x_{t+1} + e_{t+1})\ h_{t+1}^{[1]}$$
$$+ \ldots + (x_n + e_n)\ h_n^{[1]} = 0$$
$$\vdots$$
$$f_1\ h_1^{[d-2]} + \ldots + f_t\ h_t^{[d-2]} + (x_{t+1} + e_{t+1})\ h_{t+1}^{[d-2]}$$
$$+ \ldots + (x_n + e_n)\ h_n^{[d-2]} = 0$$

Let $\left(\theta_j^{(0)}\right)$ denote the equation

$f_1\ h_1^{[j]} + \ldots + f_t\ h_t^{[j]} + (x_{t+1} + e_{t+1})\ h_{t+1}^{[j]} + \ldots + (x_n + e_n)\ h_n^{[j]} = 0$ for each j = 0, 1, …, d–2 and $\theta^{(d-1)} = \left(\left(\theta_0^{(0)}\right), \left(\theta_1^{(0)}\right), \ldots, \left(\theta_{d-2}^{(0)}\right)\right)$ represent the above system of d–1 equations $\left(\theta_0^{(0)}\right)$, $\left(\theta_1^{(0)}\right)$, …, $\left(\theta_{d-2}^{(0)}\right)$.

In the above system $\theta^{(d-1)}$ of d–1 equations, one needs to eliminate $f_1, f_2, \ldots, f_t$. For each i = 1, 2,…, t, the elimination of $f_i$, results in a system, say $\theta^{(d-1-i)}$, of d–1–i equations, say $\left(\theta_0^{(i)}\right)$, $\left(\theta_1^{(i)}\right)$, …, $\left(\theta_{d-2-i}^{(i)}\right)$. Let $\theta^{(d-1-i)} = (\left(\theta_0^{(i)}\right), \left(\theta_1^{(i)}\right), \ldots, \left(\theta_{d-2-i}^{(i)}\right))$.

Therefore, $\theta^{(d-1-t)} = \left(\left(\theta_0^{(t)}\right), \left(\theta_1^{(t)}\right), \ldots, \left(\theta_{d-2-t}^{(t)}\right)\right)$ would represent the system $\theta^{(d-1-t)}$ of d–1–t equations obtained after the elimination $f_1, f_2, \ldots, f_t$ from the system $\theta^{(d-1)}$.

The recursive procedure for the elimination of $f_1, f_2, \ldots, f_t$ from the system $\theta^{(d-1)} = \left(\theta_0^{(0)}\right), \left(\theta_1^{(0)}\right), \ldots, \left(\theta_{d-2}^{(0)}\right)$ is given as follows.

For each i = 1, 2, …, t, perform the operation.

Define

$$Z_{(i-1)} = \frac{\text{Coefficient of } f_i \text{ in equation } (\theta^{(i-1)}) \text{ of system } \theta^{(d-i)}}{\text{Coefficient of } f_i \text{ in equation } (\theta_0^{(i-1)}) \text{ of system } \theta^{(d-i)}}.$$

Then (to eliminate $f_i$ from the system $\theta^{(d-i)}$ of d–i equations) multiply the equation $\left(\theta_j^{(i-1)}\right)$ by $Z_{(i-1)}^{[j]}$ and subtract the equation $\left(\theta_{j+1}^{(i-1)}\right)$ from the equation $\left(\theta_j^{(i-1)}\right)$ for each j = 0, 1, 2, …, d–2–i. Then one obtains the system $\theta^{(d-1-i)} = \left(\left(\theta_0^{(i)}\right), \left(\theta_1^{(i)}\right), \ldots, \left(\theta_{d-2-i}^{(i)}\right)\right)$ of d–1–i equations after the elimination of $f_i$.

The reduced system $\theta^{(d-1-t)}$ of d–1–t equations $\left(\theta_0^{(t)}\right)$, $\left(\theta_1^{(t)}\right)$, …, $\left(\theta_{d-2-t}^{(t)}\right)$ after the elimination of $f_1, f_2, \ldots, f_t$ from the system $\theta^{(d-1)} = \left(\left(\theta_0^{(0)}\right), \left(\theta_1^{(0)}\right), \ldots, \left(\theta_{d-2}^{(0)}\right)\right)$ is of the following form.

$$s_0' = (x_{t+1} + e_{t+1})\, h_{t+1}' + (x_{t+2} + e_{t+2})\, h_{t+2}' + \ldots + (x_n + e_n)\, h_n'$$

$$s_1' = (x_{t+1} + e_{t+1})\, h_{t+1}'^{[1]} + (x_{t+2} + e_{t+2})\, h_{t+2}'^{[1]} + \ldots + (x_n + e_n)\, h_n'^{[1]}$$

$$\vdots$$

and $$s_{d-2-t}' = (x_{t+1} + e_{t+1})\, h_{t+1}'^{[d-2-t]} + (x_{t+2} + e_{t+2})\, h_{t+2}'^{[d-2-t]} + \ldots + (x_n + e_n)\, h_n'^{[d-2-t]}.$$

The above system of equations can be rewritten as,

$$(s_0', \ldots, s_{d-2-t}')$$

$$= (x_{t+1} + e_{t+1}, \ldots, x_n + e_n) \begin{bmatrix} h_1' & h_2' & \cdots & h_n' \\ h_{t+1}'^{[1]} & h_{t+2}'^{[1]} & \cdots & h_n'^{[1]} \\ \vdots & \vdots & \ddots & \vdots \\ h_{t+1}'^{[d-2-t]} & h_{t+2}'^{[d-2-t]} & \cdots & h_n'^{[d-2-t]} \end{bmatrix}^T$$

$$= (x_{t+1} + e_{t+1}, \ldots, x_n + e_n)\, H^{(t)^T}$$

$$= (x_{t+1}, \ldots, x_n)\, H^{(t)^T} + (e_{t+1}, \ldots, e_n)\, H^{(t)^T}$$

$= x^{(t)} H^{(t)^T} + e^{(t)} H^{(t)^T}$

$= (0) + e^{(t)} H^{(t)^T}$ (by Lemma 4.2.2)

$= e^{(t)} H^{(t)^T}$

i.e., $(s'_0, \ldots, s'_{d-2-t}) = e^{(t)} H^{(t)^T}$ (4.3)

where

$e^{(t)} = (e_{t+1}, e_{t+2}, \ldots, e_n)$ is the deletion of the first t coordinates of e and

$$H^{(t)} = \begin{bmatrix} h'_{t+1} & h'_{t+2} & \cdots & h'_n \\ h'^{[1]}_{t+1} & h'^{[1]}_{t+2} & \cdots & h'^{[1]}_n \\ \vdots & \vdots & \ddots & \vdots \\ h'^{[d-2-t]}_{t+1} & h'^{[d-2-t]}_{t+2} & \cdots & h'^{[d-2-t]}_n \end{bmatrix}$$

is the $(d-1-t) \times (n-t)$ matrix, which is the deletion of the first t columns of H.

If $s'_i = 0$ for each i, then it is concluded that the received vector y contains no errors. Then solving the system $\theta^{(d-1)}$ of $d-1$ equations, one can determine the values for the unknowns $f_1, f_2, \ldots, f_t$. If $s'_i \neq 0$ for some i, then the received vector y contains errors. To find the error-vector

$e = (0, 0, \ldots, 0, e_{t+1}, e_{t+2}, \ldots, e_n)$, one is to continue the following procedure with the known syndrome values $s'_i$, $i = 0, 1, \ldots, d-2-t$. Actually, one needs to find the error vector $e^{(t)} = (e_{t+1}, e_{t+2}, \ldots, e_n)$.

Since the rank norm of the error-vector e is assumed to be m, $e^{(t)}$ can be written in the following form:

$e^{(t)} = EY = (E_1, E_2, \ldots, E_m)Y,$ (4.4)

where ($E_1$, $E_2$, …, $E_m$) are linearly independent over GF(q), and $Y = \left[Y_{ij}\right]_{i,j=1,t+1}^{m,n}$ is an $m \times (n - t)$ matrix of rank m with entries from GF(q).

Equation (4.3) becomes,

$$(s'_0, \ldots, s'_{d-2-t}) = EYH^{(t)^T} = EX, \qquad (4.5)$$

where the transposed matrix of $X = YH^{(t)^T}$ has the form

$$X^T = \begin{bmatrix} x'_1 & x'_2 & \cdots & x'_m \\ x'^{[1]}_1 & x'^{[1]}_2 & \cdots & x'^{[1]}_m \\ \vdots & \vdots & \ddots & \vdots \\ x'^{[d-2-t]}_1 & x'^{[d-2-t]}_2 & \cdots & x'^{[d-2-t]}_m \end{bmatrix} \qquad (4.6)$$

and

$$x'_p = \sum_{j=t+1}^{n} Y_{pj} h'_j, \; p = 1, 2, \ldots, m. \qquad (4.7)$$

are linearly independent over GF(q).

Equation (4.5) is equivalent to the following system of equations in the unknowns $E_1$, $E_2$, …, $E_m$, $x'_1, x'_2, \ldots, x'_m$,

$$s'_p = \sum_{i=1}^{m} E_i x'^{[p]}_i, \; p = 0, 1, \ldots, d-2-t. \qquad (4.8)$$

If the solution of the system (4.8) has been found, then from (4.7) and (4.4), one can determine the matrix Y and the error-vector $e^{(t)}$ respectively. Note that the above system (4.8) has many solutions for a specified m; for $m \le \left\lfloor \frac{d-1-t}{2} \right\rfloor$. However, all solutions lead to the same vector $e^{(t)}$. Thus, the decoding problem reduces to finding the solution of the system (4.8).

Let $S'(z) = \sum_{j=0}^{d-2-t} s'_j z^{[j]}$ be termed as the syndrome polynomial.

Set $F_0(z) = z^{[d-1-t]}$, $F_1(z)=S'(z)$ and employ Euclid's division algorithm until reaching a remainder polynomial $F_{m+1}(z)$ such that $\deg(F_m(z)) \geq q^{(d-1-t)/2}$ and $\deg(F_{m+1}(z)) < q^{(d-1-t)/2}$. Then it is concluded that the received vector y has m rank-errors; that is, $r[e^{(t)}; q] = m$.

Let $\Delta(z) = \sum_{p=0}^{m} \Delta_p z^{[p]}$, $\Delta_m = 1$ be a polynomial whose roots are all possible linear combinations of $E_1, E_2, \ldots, E_m$ with coefficients from GF(q).

Using the coefficients of the remainder polynomial $F_{m+1}(z)$, the coefficients $\Delta_0, \Delta_1, \ldots, \Delta_m$ of the polynomial $\Delta(z)$ can be determined recursively as follows:

Let j be such that $s_i' = 0$ and $s_j' \neq 0$ for $i < j$.

Then $\Delta_0 = f_j / s_j'$,

$$\Delta_p = (f_{j+p} - \sum_{i=0}^{p-i} \Delta_{i-p+j-1}^{[i]} / s_j'^{[p]}, \quad p = 1, 2, \ldots, m,$$

where $f_j$ is the coefficient of $F_{m+1}(z)$ for degree [j] and for $j + p \geq m$, set $f_{j+p} = 0$.

Determine the roots $E_1, E_2, \ldots, E_m$ of $\Delta(z)$ that are linearly independent over GF(q). Methods for determining the roots of $\Delta(z)$ are described in [25]. After determining the roots $E_1, E_2, \ldots, E_m$ of $\Delta(z)$, consider the following truncated system;

$$s_p' = \sum_{i=1}^{m} E_i x_i'^{[p]}, \quad p = 0, 1, \ldots, m-1 \qquad (4.9)$$

Solving the system (4.9), one obtains $x_j'^{[p]}$ and hence the error-vector $e^{(t)}$. Then, by substituting $e^{(t)}$ in the system $\theta^{(d-1)}$ of $d - 1$ equations, one can determine the values for the unknowns $f_1, f_2, \ldots, f_t$. Hence $x = y - e$, which is the actually transmitted codeword.

Summary of the error-erasure decoding technique:

**Step 1:** Compute the syndrome values ($s'_0, \ldots, s'_{d-2-t}$) for the vector $y = (f_1, f_2, \ldots, f_t, x_{t+1}, \ldots, x_n) + (0, 0, \ldots, 0, e_{t+1}, \ldots, e_n)$ and the corresponding syndrome polynomial

$$S'(z) = \sum_{j=0}^{d-2-t} s'_j z^{[j]} .$$

**Step 2:** Set $F_0(z) = z^{[d-1-t]}$, $F_1(z) = S'(z)$ and employ Euclid's division algorithm until reaching a $F_{m+1}(z)$ such that deg $(F_m(z)) \geq q^{(d-1-t)/2}$ and deg $(F_{m+1}(z)) < q^{(d-1-t)/2}$. Let $\Delta(z) = \sum_{p=0}^{m} \Delta_p z^{[p]}$, $\Delta_m = 1$ be such that its roots are all possible linear combinations of $E_1, E_2, \ldots, E_m$ over GF(q).

**Step 3:** Determine the coefficients of $\Delta(z)$ as follows:
Let j be such that $s'_i = 0$ and $s'_j \neq 0$ for $i < j$.

$$\Delta_0 = f_j / s'_j,$$

$$\Delta_p = (f_{j+p}, \sum_{i=0}^{p-i} \Delta_{i-p+j-1}^{[i]} / s_j'^{[p]}, \; p = 1, 2, \ldots, m.$$

where $f_j$ is the coefficient of $F_{m+1}(z)$ for degree [j] and $f_{j+p} = 0$ for $j + p \geq m$.

**Step 4:** Compute the roots $E_1, E_2, \ldots, E_m$ of $\Delta(z)$ that are linearly independent over GF(q). Then, by substituting $E_1, E_2, \ldots, E_m$ in system (4.9), one can determine $x_j^{[p]}$ and hence the error-vector $e^{(t)}$.

**Step 5:** Substituting $e^{(t)}$ in system $\theta^{(d-1)}$, one can obtain the values for the unknown $f_1, f_2, \ldots, f_t$.

**Step 6:** The transmitted codeword is then obtained as

$$x = y - e.$$

Using the combined error-erasure decoding technique to the class of $[n,k,d]_{q^N}$ MRD codes presented above, one can correct any combination of $m \leq r$ rank-errors and $t \leq s$ erasures in an erroneously received vector as long as $2r + s < d$, where d is the minimum-rank distance of the MRD codes.

A detailed description of the combined error-erasure decoding technique to the class of $[n,k,d]_{2^N}$ MRD codes presented above is demonstrated through the following example which applies the combined error-erasure decoding technique to the $[7,1,7]_{2^8}$ MRD code for the correction of 2 errors and 2 erasures in an erroneously received vector.

***Example 4.2.1***: Let $\Gamma = [7,1,7]_{2^8}$ be the MRD code defined over $GF(2^8) = \{0, 1, \alpha, \ldots, \alpha^{2^8-2}\}$ with the parity-check matrix H:

$$H = \begin{bmatrix} 1 & \alpha & \alpha^2 & \alpha^3 & \alpha^4 & \alpha^5 & \alpha^6 \\ 1 & \alpha^2 & \alpha^4 & \alpha^6 & \alpha^8 & \alpha^{10} & \alpha^{12} \\ 1 & \alpha^4 & \alpha^8 & \alpha^{12} & \alpha^{16} & \alpha^{20} & \alpha^{24} \\ 1 & \alpha^8 & \alpha^{16} & \alpha^{24} & \alpha^{32} & \alpha^{40} & \alpha^{48} \\ 1 & \alpha^{16} & \alpha^{32} & \alpha^{48} & \alpha^{64} & \alpha^{80} & \alpha^{96} \\ 1 & \alpha^{32} & \alpha^{64} & \alpha^{96} & \alpha^{128} & \alpha^{160} & \alpha^{192} \end{bmatrix}$$

where $\alpha$ is a root of the primitive polynomial $x^8 + x^6 + x^5 + x + 1$ over GF(2).

Let $x = (0, 0, 0, 0, 0, 0, 0) \in \Gamma$ be the transmitted codeword and $y = \alpha^{31}, \alpha^{147}, 0, 0, *_5, 0, *_7)$ be the received vector, where '$*_5$' and '$*_7$' represent erasures. Note that the received vector y has two erasures.

Let $f = (0, 0, 0, 0, f_5, 0, f_7)$. Then replacing the erasures with unknowns $f_5$ and $f_7$, one has the vector $y' = (\alpha^{31}, \alpha^{147}, 0, 0, f_5, 0, f_7)$ where the unknowns $f_5$ and $f_7$ are to be determined. Let e =

$(e_1, e_2, e_3, e_4, 0, e_6, 0)$ be the error-vector to be determined. Then $x = y' - e$ would be the transmitted code word. The syndrome of $y' - e$ is

$$S = (y' - e)H^T$$
$$= (\alpha^{31}, \alpha^{147}, 0, 0, f_5, 0, f_7)H^T - (e_1, e_2, e_3, e_4, 0, e_6, 0)H^T \quad (4.10)$$

Equating each component on the right-hand side of (4.10) to 0,

$$\alpha^4 f_5 + \alpha^6 f_7 + 1 = e_1 + \alpha e_2 + \alpha^2 e_3 + \alpha^3 e_4 + \alpha^5 e_6,$$
$$\alpha^8 f_5 + \alpha^{12} f_7 + \alpha^{147} = e_1 + \alpha^2 e_2 + \alpha^4 e_3 + \alpha^6 e_4 + \alpha^{10} e_6,$$
$$\alpha^{16} f_5 + \alpha^{24} f_7 + \alpha^{108} = e_1 + \alpha^4 e_2 + \alpha^8 e_3 + \alpha^{12} e_4 + \alpha^{20} e_6,$$
$$\alpha^{32} f_5 + \alpha^{48} f_7 + \alpha^{113} = e_1 + \alpha^8 e_2 + \alpha^{16} e_3 + \alpha^{24} e_4 + \alpha^{40} e_6,$$
$$\alpha^{64} f_5 + \alpha^{96} f_7 + \alpha^{41} = e_1 + \alpha^{16} e_2 + \alpha^{32} e_3 + \alpha^{48} e_4 + \alpha^{80} e_6,$$

and

$$\alpha^{128} f_5 + \alpha^{192} f_7 + \alpha^{62} = e_1 + \alpha^{32} e_2 + \alpha^{64} e_3 + \alpha^{96} e_4 + \alpha^{160} e_6.$$

Let $\theta^{(6)} = \left(\left(\theta_0^{(0)}\right), \left(\theta_1^{(0)}\right), \left(\theta_2^{(0)}\right), \left(\theta_3^{(0)}\right), \left(\theta_4^{(0)}\right), \left(\theta_5^{(0)}\right)\right)$ denote the above system of 6 equations.

Here

$$Z_{(0)} = \frac{\text{Coefficient of } f_5 \text{ in equation } (\theta_1^{(0)}) \text{ of system } \theta^{(6)}}{\text{Coefficient of } f_5 \text{ in equation } (\theta_0^{(0)}) \text{ of system } \theta^{(6)}} = \alpha^4.$$

The eliminate $f_5$ from the above system $\theta^{(6)}$ of 6 equations, multiply the equation $\left(\theta_j^{(0)}\right)$ by $Z_{(0)}^{[j]}$ and subtract the equation $\left(\theta_{j+1}^{(0)}\right)$ from the equation $\left(\theta_j^{(0)}\right)$ for each j = 0, 1, 2, 3, 4. Then one obtains the following reduced system $\theta^{(5)} = \left(\left(\theta_0^{(1)}\right), \left(\theta_1^{(1)}\right), \left(\theta_2^{(1)}\right), \left(\theta_3^{(1)}\right), \left(\theta_4^{(1)}\right)\right)$ of 5 equations after the elimination of $f_5$.

$$\alpha^{149} f_7 + \alpha^{78} = \alpha^{23}e_1 + \alpha^{52} e_2 + \alpha^{143} e_3 + \alpha^{203} e_4 + \alpha^{206} e_6,$$
$$\alpha^{43} f_7 + \alpha^{105} = \alpha^{46}e_1 + \alpha^{104} e_2 + \alpha^{31} e_3 + \alpha^{151} e_4 + \alpha^{157} e_6,$$
$$\alpha^{86} f_7 + \alpha^{104} = \alpha^{92}e_1 + \alpha^{208} e_2 + \alpha^{62} e_3 + \alpha^{47} e_4 + \alpha^{59} e_6,$$
$$\alpha^{172} f_7 + \alpha^{16} = \alpha^{184}e_1 + \alpha^{161} e_2 + \alpha^{124} e_3 + \alpha^{94} e_4 + \alpha^{118} e_6,$$

and

$$\alpha^{89} f_7 + \alpha^{240} = \alpha^{113}e_1 + \alpha^{67} e_2 + \alpha^{248} e_3 + \alpha^{188} e_4 + \alpha^{236} e_6.$$

Now

$$Z_{(1)} = \frac{\text{Coefficient of } f_7 \text{ in equation } (\theta_1^{(1)}) \text{ of system } \theta^{(5)}}{\text{Coefficient of } f_7 \text{ in equation } (\theta_0^{(1)}) \text{ of system } \theta^{(5)}} = \alpha^{149}.$$

To eliminate $f_7$ from the above system $\theta^{(5)}$ of 5 equations, multiply the equation $\left(\theta_j^{(1)}\right)$ by $Z_{(1)}^{[j]}$ and subtract the equation $\left(\theta_{j+1}^{(1)}\right)$ from the equation $\left(\theta_j^{(1)}\right)$ for each j = 0, 1, 2, 3. Then one obtains the following reduced system $\theta^{(4)} = \left(\left(\theta_0^{(1)}\right), \left(\theta_1^{(1)}\right), \left(\theta_2^{(1)}\right), \left(\theta_3^{(1)}\right)\right)$ of 4 equations after the elimination of $f_7$.

$$\alpha^{95} = \alpha^{197}e_1 + \alpha^{71} e_2 + \alpha^{131} e_3 + \alpha^{196} e_4 + \alpha^{10} e_6,$$
$$\alpha^{68} = \alpha^{139}e_1 + \alpha^{142} e_2 + \alpha^{7} e_3 + \alpha^{137} e_4 + \alpha^{20} e_6,$$
$$\alpha^{44} = \alpha^{23}e_1 + \alpha^{29} e_2 + \alpha^{14} e_3 + \alpha^{19} e_4 + \alpha^{40} e_6,$$
and $$\alpha^{48} = \alpha^{46}e_1 + \alpha^{58} e_2 + \alpha^{28} e_3 + \alpha^{38} e_4 + \alpha^{80} e_6.$$

The above system $\theta^{(4)}$ of 4 equations can be written as
$(\alpha^{95}, \alpha^{68}, \alpha^{44}, \alpha^{48})$

$$= (e_1, e_2, e_3, e_4, e_6) \begin{bmatrix} \alpha^{197} & \alpha^{71} & \alpha^{131} & \alpha^{196} & \alpha^{10} \\ \alpha^{139} & \alpha^{142} & \alpha^{7} & \alpha^{137} & \alpha^{20} \\ \alpha^{23} & \alpha^{29} & \alpha^{14} & \alpha^{19} & \alpha^{40} \\ \alpha^{46} & \alpha^{58} & \alpha^{28} & \alpha^{38} & \alpha^{80} \end{bmatrix}^T \quad (4.11)$$

Let $(s'_0, s'_1, s'_2, s'_3) = (\alpha^{95}, \alpha^{68}, \alpha^{44}, \alpha^{48})$
and $e^{(2)} = (e_1, e_2, e_3, e_4, e_6)$

Here $S'(z) = \alpha^{95}z^{[0]} + \alpha^{68} z^{[1]} + \alpha^{44} z^{[2]} + \alpha^{48} z^{[3]}$.

Let $F_0(z) = z^{[4]}$ and

$F_1(z) = \alpha^{95}z + \alpha^{68} z^2 + \alpha^{44} z^4 + \alpha^{48} z^8$

Dividing $F_0(z)$ on the right by $F_1(z)$.

$F_0(z) = (\alpha^{207}z^8 + \alpha^{203} z^4 + \alpha^{227} z^2 + \alpha^{254} z + \alpha^{199})$
$* F_1(z) + F_2(z)$.
where $F_2(z) = \alpha^{103}z^4 + \alpha^{136} z^2 + \alpha^{39} z$.

Dividing $F_1(z)$ on the right by $F_2(z)$.

$F_1(z) = (\alpha^{200}z^4 + \alpha^{233} z^2 + \alpha^{136} z) * F_2(z) + F_3(z)$,
where $F_3(z) = \alpha^{48}z^2 + \alpha^{72} z$.

Since $\deg(F_2(z)) = 4 \geq 2^2$ and $\deg(F_3(z)) = 2 < 2^2$, it follows the received vector has $m = 2$ rank-errors i.e., $r[e^{(2)}; 2] = 2$.

Since $e^{(2)}$ is of rank norm $m = 2$,

$$e^{(2)} = (E_1, E_2)\begin{bmatrix} Y_{11} & Y_{12} & Y_{13} & Y_{14} & Y_{16} \\ Y_{21} & Y_{22} & Y_{23} & Y_{24} & Y_{26} \end{bmatrix} \qquad (4.12)$$

$$= \left(\bar{Y}_1, \bar{Y}_2, \bar{Y}_3, \bar{Y}_4, \bar{Y}_6\right) \text{ (say)}$$

where $E_1$ and $E_2$ are linearly independent over GF(2), and $Y = [Y_{ij}]$ is the $2 \times 5$ matrix of rank 2 with entries from GF(2).

Equation (4.11) becomes,
$(\alpha^{95}, \alpha^{68}, \alpha^{44}, \alpha^{48})$

$$= \left(\bar{Y}_1, \bar{Y}_2, \bar{Y}_3, \bar{Y}_4, \bar{Y}_6\right) \begin{bmatrix} \alpha^{197} & \alpha^{71} & \alpha^{131} & \alpha^{196} & \alpha^{10} \\ \alpha^{139} & \alpha^{142} & \alpha^{7} & \alpha^{137} & \alpha^{20} \\ \alpha^{23} & \alpha^{29} & \alpha^{14} & \alpha^{19} & \alpha^{40} \\ \alpha^{46} & \alpha^{58} & \alpha^{28} & \alpha^{38} & \alpha^{80} \end{bmatrix}^T \qquad (4.13)$$

Here $\Delta$ (z) $=\Delta_0 Z^{[0]} + \Delta_1 Z^{[1]} + \Delta_2 Z^{[2]}$ with $\Delta_0 = \alpha^{232}$, $\Delta_1 = \alpha^{160}$ and $\Delta_2 = \alpha^{160}$.

The roots of $\Delta$ (z) are $\alpha^{31}$, $\alpha^{147}$, $\alpha^{149}$ and 0. Take $E_1 = \alpha^{31}$, and $E_2 = \alpha^{147}$.

$$\alpha^{95} = \begin{cases} E_1(\alpha^{197}y_{11} + \alpha^{71}y_{12} + \alpha^{131}y_{13} + \alpha^{196}y_{14} + \alpha^{10}y_{16}) + \\ E_2(\alpha^{197}y_{21} + \alpha^{71}y_{22} + \alpha^{131}y_{23} + \alpha^{196}y_{24} + \alpha^{10}y_{26}) \end{cases}$$

and

$$\alpha^{68} = \begin{cases} E_1(\alpha^{139}y_{11} + \alpha^{142}y_{12} + \alpha^{7}y_{13} + \alpha^{137}y_{14} + \alpha^{20}y_{16}) + \\ E_2(\alpha^{139}y_{21} + \alpha^{142}y_{22} + \alpha^{7}y_{23} + \alpha^{137}y_{24} + \alpha^{20}y_{26}) \end{cases}.$$

Solving the above equation,

$$Y = \begin{bmatrix} 1 & 0 & 0 & 0 & 0 \\ 0 & 1 & 0 & 0 & 0 \end{bmatrix}.$$

Therefore, e = ($\alpha^{31}$, $\alpha^{147}$, 0, 0, 0, 0, 0). By substituting e = ($\alpha^{31}$, $\alpha^{147}$, 0, 0, 0, 0, 0) in the system $\theta^{(6)}$ of 6 equations, one obtains $f_5 = 0$ and $f_7 = 0$.

Hence $x = y' - e$

$= (\alpha^{31}, \alpha^{147}, 0, 0, 0, 0, 0) - (\alpha^{31}, \alpha^{147}, 0, 0, 0, 0, 0)$

$= (0, 0, 0, 0, 0, 0, 0)$,

which is the actually transmitted codeword.

## 4.3 Invertible q-Cyclic RD Codes

There are two categories of techniques for controlling transmission errors in data transmission systems: the Forward-Error Control (FEC) scheme and the Automatic-Repeat – Request (ARQ) scheme. In an FEC system, an error-correcting code is used, when the receiver detects the presence of errors in a received vector, it attempts to determine the error locations and then corrects the errors. If the exact locations of errors are determined, the received vector will be correctly decoded; if the receiver fails to determine the exact locations of errors, the received vector will be decoded incorrectly and erroneous data will be delivered to be destination. In an ARQ system, a code with good error-detecting capability is used. At the receiver, the syndrome of the received vector is computed. If the syndrome is zero, the received vector is assumed to be error-free and is accepted by the receiver. At the same time, the receiver notifies the transmitter, via a return channel, that the transmitted codeword has been successfully received. If the syndrome is not zero, errors are detected in the received vector. Then the transmitter is instructed, through the return channel, to retransmit the same codeword. Retransmission continues until the codeword is successful received.

The throughput efficiency (or throughput) is a measure of performance of an ARQ system. Throughout efficiency is defined as the ratio of the average number of information symbols successfully accepted by the receiver per unit of time to the total number of symbols that could be transmitted per unit of time. The retransmission of an erroneously received vector continues until the received vector is successfully received by the receiver. Though the retransmission requests provide a powerful means of improving reliability performance at the cost of a reduction in throughput, the frequency of retransmission must be reduced to improve the throughout efficiency. Both systems have their own limitations and drawbacks. To improve the throughput efficiency, a hybrid - ARQ scheme which is a combination of both FEC and ARQ, came into existence [25].

There are two types of hybrid - ARQ schemes : type I hybrid – ARQ scheme [25] and type II hybrid – ARQ scheme [25]. This section deals only the type II hybrid-ARQ scheme, the description of it is given in the following.

The type II hybrid - ARQ scheme is devised based on the concept that the parity - check symbols for error correction are sent to the receiver only when they are needed [25]. Two linear codes are used in this type of scheme; one is a high-rate (n, k) code $C_0$, which is designed for error detection only, the other is a half-rate invertible (2k, k) code $C_1$, which is designed for simultaneous error correction and error detection.

When a message, say u of k information symbols is ready for transmission, it is encoded into a codeword v = (f(u), u) of n symbols based on the error-detecting code $C_0$, where f(u) denotes the n–k parity-check symbols. The codeword v = (f(u), u) is then transmitted. At the same time, the transmitter computes the k parity-check symbols, denoted by q(u), based on the message u and the half-rate invertible code $C_1$. Clearly the 2k-tuple (q(u), u) is a buffer of the transmitter for later use.

The invertible property facilities the data recovery process. This section carries out a comprehensive study on the invertible property enjoyed by the class of $[n,k,d]_{q^n}$ q-Cyclic RD codes.

This section is divided into three subsections. In subsection 1, the systematic encoding of the class of $[n,k,d]_{q^n}$ q-cyclic RD codes analogues to the class of cyclic codes is given. Making use of the systematic encoding, the subsection 2 gives the shortening technique to the class of $[n,k,d]_{q^n}$ q-cyclic RD codes. Subsection 3 studies the invertible property enjoyed by the class of $[n,k,d]_{q^n}$ q-cyclic MRD codes for both the cases $n - k \geq k$ and $n - k < k$.

Systematic Encoding of q-Cyclic RD Codes

The generator and parity-check matrices for linear codes greatly simplify encoding at the transmitter and error detection at the receiver. The problem of recovering the message block from a codeword can be greatly simplified through the use of systematic encoding. As in the case of ordinary cyclic codes, systematic encoding for the class of $[n,k,d]_{q^n}$ q-Cyclic RD codes can be done either by means of a generator polynomial or by means of a generator matrix in systematic (or standard) form.

**(a) Systematic encoding through generator polynomial**

Let G(z) = $\sum_{i=0}^{n-k} G_i z^{[i]}$ be the given generator polynomial of an $[n,k,d]_{q^n}$ q-Cyclic RD code, say C.

Let H(z) = $\sum_{i=0}^{k} H_i z^{[i]}$ be such that $z^{[n]} - z = H(z) * G(z)$. If one assumes that $g_{n-1}$, $g_{n-2}$, …, $g_{n-k}$ are information symbols, then s/he can determine the parity-check symbols $g_{n-k-1}$, $g_{n-k-2}$, …, $g_0$ as follows.

Let $u(z) = g_{n-1}z^{[n-1]} + g_{n-2}z^{[n-2]} + \ldots + g_{n-k}z^{[n-k]}$ be a message polynomial to be encoded.

Divide u(z) on the right by G(z);

$u(z) = q(z) * G(z) + f(z)$, deg (f (z)) < [n–k].

where q(z) is the quotient polynomial and f(z) is the remainder polynomial.

Then the coefficients $g_{n-i}$ for degree $g_{n-i}$ for degree [n–i], i = k+1, k+2, …, n of the remainder polynomial f (z) will be the parity-check symbols. Then g(z) = u(z) – f(z) is the systematically encoded code polynomial corresponding to the message polynomial u(z), where f(z) is called as the parity-check polynomial.

**(b) Systematic encoding through generator matrix**

The systematic encoding of the $[n,k,d]_{q^n}$ q-cyclic RD code C can also be carried out through generator and parity-check matrices in systematic form. This is accomplished in the following.

The generator matrix and the parity-check matrix corresponding to G(z) and H(z) of C in systematic form can be obtained as follows:

Dividing $z^{[n-k+i]}$ on the right by G(z) for each i = 0, 1, …, k–1, one gets

$z^{[n-k+i]} = q_i(z) * G(Z) + f_i(z)$, where $q_i(z)$ is the quotient polynomial and $f_i(z) = f_{i0}z^{[0]} + f_{i1}z^{[1]} + \ldots + f_{i,\, n-k-1}z^{[n-k-1]}$ is the remainder polynomial.

Since G(z) is a right divisor of $z^{[n-k-i]} - f_i(z)$, i = 0, 1, …, k–1, they are code polynomials of C.

Arranging these k code polynomials as rows of a k × n matrix, one obtains the matrix say G:

$$G = \begin{bmatrix} f_{00} & f_{01} & \cdots & f_{0,n-k-1} & 1 & 0 & 0 & \cdots & 0 \\ f_{01} & f_{11} & \cdots & f_{1,n-k-1} & 0 & 1 & 0 & \cdots & 0 \\ \vdots & \vdots & \ddots & \vdots & \vdots & \vdots & \vdots & \ddots & 0 \\ f_{k-1,0} & f_{k-1,1} & \cdots & f_{k-1,n-k-1} & 0 & 0 & 0 & \cdots & 1 \end{bmatrix}$$

which is the generator matrix of C in systematic form. Recall that any set of k linearly independent vectors can be used as the rows of the generator matrix to form a k-dimensional linear code.

The corresponding parity-check matrix, say H in systematic form is given by

$$H = \begin{bmatrix} 1 & 0 & 0 & \ldots & 0 & -f_{00} & -f_{10} & \ldots & -f_{k-1,0} \\ 0 & 1 & 0 & \ldots & 0 & -f_{01} & -f_{11} & & -f_{k-1,1} \\ \vdots & \vdots & \ddots & & \vdots & \vdots & \vdots & \ddots & \vdots \\ 0 & 0 & 0 & \ldots & 1 & -f_{0,n-k-1} & -f_{1,n-k-1} & \ldots & -f_{k-1,n-k-1} \end{bmatrix}$$

Note that $GH^T = (0)$. If $u = (g_{n-k}, g_{n-k+1}, \ldots, g_{n-1})$ is a message vector, then $v = uG = (g_0, g_1, \ldots, g_{n-1})$ is the systematically encoded codeword corresponding to the message vector u and one can note that $H^T = (0)$, where $g_0, g_1, \ldots, g_{n-k-1}$ are the corresponding parity-check symbols.

The systematic encoding for the class of $[n,k,d]_{q^n}$ q-Cyclic RD codes given above is demonstrated through the following example.

***Example 4.3.1.1:*** Consider the $[5,3]_{2^5}$ 2-Cyclic RD code generated by $G(z) = \alpha^{24} z + \alpha^3 z^2 + \alpha^2 z^4$, where $\alpha$ is a root of the primitive polynomial $x^5 + x^2 + 1$ over GF(2).

Dividing $z^{[2]}$ on the right by G[z], one gets

$z^4 = \alpha^{29} z * G(z) + \alpha z^2 + \alpha^{22} z$ with $f_1(z) = \alpha z^2 + \alpha^{22} z$.

Dividing $z^{[3]}$ on the right by G[z], one gets

$z^8 = (\alpha^{27} z^2 + z) * G(z) + \alpha^7 z^2 + \alpha^{24} z$ with $f_2(z) = \alpha^7 z^2 + \alpha^{24} z$.

Dividing $z^{[4]}$ on the right by G[z], one gets

$Z^{16} = (\alpha^{23} z^4 + z^2 + \alpha^{12} z) * G(z) + \alpha^{20} z^2 + \alpha^5 z$ with $f_3(z) = \alpha^{20} z^2 + \alpha^5 z$.

Then the generator matrix G and the parity-check matrix H for the $[5,3]_{2^5}$ 5-Cyclic RD code in systematic form are respectively given by

$$G = \begin{bmatrix} \alpha^{22} & \alpha & 1 & 0 & 0 \\ \alpha^{24} & \alpha^{7} & 0 & 1 & 0 \\ \alpha^{5} & \alpha^{20} & 0 & 0 & 1 \end{bmatrix}$$

and

$$H = \begin{bmatrix} 1 & 0 & \alpha^{22} & \alpha^{24} & \alpha^{5} \\ 0 & 1 & \alpha & \alpha^{7} & \alpha^{20} \end{bmatrix}.$$

Let $u_1(z) = \alpha^5 z^{[4]} + z^{[3]} + \alpha^{23} z^{[2]}$ be the message polynomial to be encoded. Dividing $u_1(z)$ on the right by $G(z)$, one gets.

$$u_1(z) = (\alpha^{28}z^4 + \alpha^{12}z^2 + \alpha^6 z) * G(z) + \alpha^{17} z^2 + \alpha^{30} z.$$

Then $g_1(z) = \alpha^{30}z + \alpha^{17} z^{[1]} + \alpha^{23} z^{[2]} + z^{[3]} + \alpha^5 z^{[4]}$ is the systematically encoded code polynomial corresponding to the message polynomial $u_1(z)$. Note also that

$$(\alpha^{30}, \alpha^{17}, \alpha^{23}, 1, \alpha^5) H^T = (0).$$

Similarly, $g_2(z) = \alpha^{13}z + \alpha^{11}z^{[1]} + \alpha^9 z^{[2]} + 0z^{[4]}$ is the code polynomial corresponding to the message polynomial $u_2(z) = 0z^{[4]} + \alpha^{21}z^{[3]} + \alpha^9 z^{[2]}$ and that $(\alpha^{13}, \alpha^{11}, \alpha^9, \alpha^{21}, 0) H^T = (0)$.

The next subsection gives the shortening technique to the class of $[n,k,d]_{q^n}$ q-Cyclic RD codes. In shortening of an $[n,k,d]_{q^n}$ q-Cyclic RD code, each codeword is shortened (in length) by t information symbols, resulting in a linear code of length $n - t$, dimension $k - t$ with the same error-correcting capability as the original code, $t < k$.

### 4.3.2 Shortened q-Cyclic RD codes

In many applications, there are external constraints such as puncturing, extending, shortening, lengthening, expurgating, or augmenting [25], which are unrelated to error-control but determine the allowed length of the error-control code. In system design, if a code of suitable length or suitable number of information symbols can not be found, it may be possible to shorten a code to meet the requirements. In extending, a code is extended by adding an additional redundant coordinate. Thus an (n, k) code becomes an (n + 1, k) code. In shortening, a code is shortened by deleting a message coordinate from the encoding process. Thus an (n, k) code becomes an (n – 1, k – 1) code. This section gives shortening technique to the class of $[n,k,d]_{q^n}$ q-Cyclic RD codes. This shortening technique to the class of $[n,k,d]_{q^n}$ q-Cyclic RD codes enables the subsection 4.3.3, to study the invertible property for the class of $[n,k,d]_{q^n}$ q-Cyclic MRD codes when $n - k < k$.

Let C denote an $[n,k,d]_{q^n}$ q-Cyclic RD code.

Let $G(z) = \sum_{i=0}^{n-k} G_i z^{[i]}$ be a generator polynomial of C.

Let $H(z) = \sum_{i=0}^{k} H_i z^{[i]}$ be such that $z^{[n]} - z = H(z) * G(z)$. Consider the set of codewords of C for which the $t < k$ leading higher-order information symbols are identical to zero. There are $q^{n(k-t)}$ such codewords in C and these $q^{n(k-t)}$ codewords in fact form a linear subcode of C. If the t zero information symbols are deleted from each of these codewords, one obtains a set of $q^{n(k-t)}$ vectors of length n–t. These set of $q^{n(k-t)}$ shortened vectors forms an $[n-t,k-t]_{q^n}$ linear code. Call it as the shortened $[n-t,k-t]_{q^n}$ q-Cyclic RD code and denote it by $C_{(t)}$. The shortened $[n-t,k-t]_{q^n}$ q-Cyclic RD code $C_{(t)}$ has the same error-correcting capability as C. Also the encoding and decoding for the shortened $[n-t, k-t]_{q^n}$ q-Cyclic RD code $C_{(t)}$ can be accomplished by the generator polynomial G(z) and the parity-check polynomial H(z) in the same way as C. This is so because the deleted t zero information symbols do not affect the parity-check and syndrome computations.

The following example shortens the $[5,3]_{2^5}$ 2-Cyclic RD code by t = 1, 2 information symbols and describes encoding of the shortened codes.

***Example 4.3.2.1:*** Consider the $[5,3]_{2^5}$ 2-Cyclic RD code C generated by $G(z) = \alpha^{24}z + \alpha^3 z^2 + \alpha^2 z^4$, where $\alpha$ is a root of the primitive polynomial $x^5 + x^2 + 1$ over GF(2).

The generator matrix G and the parity-check matrix H in systematic form are respectively given by

$$G = \begin{bmatrix} \alpha^{22} & \alpha & 1 & 0 & 0 \\ \alpha^{24} & \alpha^7 & 0 & 1 & 0 \\ \alpha^5 & \alpha^{20} & 0 & 0 & 1 \end{bmatrix}$$

and

$$H = \begin{bmatrix} 1 & 0 & \alpha^{22} & \alpha^{24} & \alpha^5 \\ 0 & 1 & \alpha & \alpha^7 & \alpha^{20} \end{bmatrix}.$$

Shortening the code C by t = 1 leading zero information symbol, one obtains the shortened $[4,2]_{2^5}$ 2-Cyclic RD code $C_{(1)}$. Take $(\alpha^{24}, \alpha^3, \alpha^2, 0, 0) \in C$. Then the shortened codeword $(\alpha^{24}, \alpha^3, \alpha^2, 0)$ (by one information symbol) is in $C_{(1)}$.

Similarly shortening the code C by t = 2 leading zero information symbols, one obtains the shortened $[3,1]_{2^5}$ 2-Cyclic RD code $C_{(2)}$ and the shortened codeword $(\alpha^{24}, \alpha^3, \alpha^2)$, shortened by two zero information symbols from $(\alpha^{24}, \alpha^3, \alpha^2, 0, 0)$, is in $C_{(2)}$.

Having given the shortening technique to the class on $[n,k,d]_{q^n}$ q-Cyclic RD codes, the next subsection studies the invertible property enjoyed by the class of $[n,k,d]_{q^n}$ q-Cyclic MRD codes, for the cases: $n - k \geq k$ and $n - k < k$.

### 4.3.3 Invertible Property of q-Cyclic MRD Codes

If one is able to find the k message symbols in a codeword of an (n, k) F-ary linear code only with the knowledge of $n - k$ parity-check symbols through an inversion process, the code is said to be invertible. This section explores the invertible property for the class of $[n,k,d]_{q^n}$ q-Cyclic MRD codes and infers some interesting results.

It is observed that, the case when $n - k \geq k$, the $[n,k,d]_{q^n}$ q-Cyclic MRD codes are invertible and for the case when $n - k < k$, the shortened $[n-t_0,k-t_0]_{q^n}$ q-Cyclic MRD codes are invertible; where $2k - n \leq t_0 < k$.

The detailed discussion of invertible property for the class of $[n,k,d]_{q^n}$ q-Cyclic MRD codes is carried out in two cases: Case (i): $n - k \geq k$ and case (ii): $n - k < k$.

**Case (i) $n - k \geq k$**

Let C denote an $[n,k,d]_{q^n}$ q-Cyclic MRD code of length n, dimension k and minimum-rank distance d generated by $G(z) = \sum_{i=0}^{n-1} g_i z^{[i]}$ be a systematically encoded code polynomial. Then as described in subsection 4.3.1, in systematic encoding, the k leading high-order coefficients $g_{n-k}$, $g_{n-(k-1)}$, …, $g_{n-1}$ are identical to the information symbols, the $n - k$ low-order coefficients $g_0$, $g_1$, …, $g_{n-k-1}$ are the corresponding parity-check symbols.

Let $u(z) = u_{n-1} z^{[n-1]} + u_{n-2} z^{[n-2]} + \ldots + u_{n-k} z^{[n-k]}$ be a message polynomial to be encoded. Divide u(z) on the right by G(z).

$$u(z) = q(z) * G(z) + f(z), \deg (f(z)) < [n - k] \qquad (4.14)$$

where q(z) and f(z) are respectively the quotient and the remainder polynomials.

Then g(z) = u(z) – f(z) is the code polynomial corresponding to the message polynomial u(z).

Now one raises the question: Does there exists a one-to-one correspondence between the information symbols and parity-check symbols in a codeword of the $[n,k,d]_{q^n}$ q-Cyclic MRD code C so that one can talk about the invertible property of $[n,k,d]_{q^n}$ q-Cyclic MRD codes for n – k ≥ k?. The following theorem answers this question.

**THEOREM 4.3.3.1:** *If C denotes an $[n,k,d]_{q^n}$ q-Cyclic MRD code such that n – k ≥ k, then no two codewords of C will have same parity-check symbols.*

***Proof:*** Consider the generator polynomial G(z) of C in the form: $G(z) = z^{[n-k]} + G_{n-k-1}z^{[n-k+1]} + \ldots + G_0z^{[0]}$.

Let $u^{(1)}(z) = \sum_{i=1}^{k} u_i^{(1)} z^{[n-i]}$ and

$$u^{(2)}(z) = \sum_{i=1}^{k} u_i^{(2)} z^{[n-i]}$$

be two distinct message polynomials to be encoded.

Dividing $u^{(1)}(z)$ and $u^{(2)}(z)$ on the right by G(z).

$u^{(1)}(z) = q_1(z) * G(z) + f_1(z), \deg (f_1(z)) < [n–k]$ (4.15)

and

$u^{(2)}(z) = q_2(z) * G(z) + f_2(z), \quad \deg(f_2(z)) < [n-k]. \quad (4.16)$

The code polynomials corresponding to $u^{(1)}(z)$ and $u^{(2)}(z)$ respectively are,

$$g_1(z) = u^{(1)}(z) - f_1(z)$$

$$\text{and } g_2(z) = u^{(2)}(z) - f_2(z).$$

Suppose that $f_1(z) = f_2(z)$.

Subtracting (4.16) from (4.15),

$u^{(1)}(z) - u^{(2)}(z) = (q_1(z) - q_2(z)) * G(z)$

i.e., $z^{[n-k]} * (u'^{(1)}(z) - u'^{(2)}(z)) = (q_1(z) - q_2(z)) * G(z)$,

where $u^{(i)}(z) = z^{[n-k]} * u'^{(i)}(z)$ for each $i = 1, 2$.

The above equation shows that G(z) is a right divisor of $z^{[n-k]} * (u'^{(1)}(z) - u'^{(2)}(z))$. Since $z^{[n-k]}$ is relatively prime to G(z), G(z) must be a right divisor of $u'^{(1)}(z) - u'^{(2)}(z)$. However this is impossible, because $u'^{(1)}(z) \neq u'^{(2)}(z)$ and degree of $u'^{(1)}(z) - u'^{(2)}(z)$ is less than or equal to $[k - 1]$ but the degree of G(z) is $[n-k]$ which is strictly greater than $[k-1]$, since $n - k \geq k$.

Thus $f_1(z) \neq f_2(z)$. Hence the theorem.

Since the remainder f(z) resulting from dividing u(z) by G(z) is unique, the Theorem 4.3.3.1 implies that there exists an one-to-one correspondence between a message polynomial u(z) and its parity-check polynomial f(z). Therefore knowing only the parity-check polynomial f(z), the message polynomial u(z) can be determined uniquely.

In what follows, it is shown how to recover the message polynomial u(z) from its parity-check polynomial f(z). Now one has the equation (4.14):

$$u(z) = q(z) * G(z) + f(z), \ deg(f(z)) < [n–k].$$

Consider $z^{[k]} * f(z)$ :

$$z^{[k]} * f(z)= z^{[k]} * (–q(z) * G(z) + u(z))$$

$$= z^{[k]} * (–q(z) * G(z)) + z^{[k]} * u(z)$$

$$= z^{[k]} * (–q(z) * G(z)) + (u'(z) * z^{[k]}) * z^{[n–k]}$$

$$= z^{[k]} * (–q(z) * G(z)) + u'(z) * z^{[n]}$$

$$= z^{[k]} * (–q(z) * G(z)) + u'(z) * (z^{[n]}–z) + u'(z) * z$$

$$= (z^{[k]} * (–q(z))) * G(z) + u'(z) * (H(z)*G(z))+ u'(z)$$

(since G(z) is a right divisor of $z^{[n]} – z$).

$$= (z^{[k]} * (–q(z)) + u'(z) * H(z)) * G(z) + u'(z), \quad (4.17)$$

where $u'(z) = u_{n-1}^{[k]} z^{[k-1]} + u_{n-2}^{[k]} z^{[k-2]} + \ldots + u_{n-k}^{[k]} z^{[0]}$.

See that from (4.17), the message polynomial u(z) is nothing but $z^{[n–k]} * u'(z)$;

i.e., $u(z) = z^{[n–k]} * u'(z)$

$$= z^{[n–k]} *( u_{n-1}^{[k]} z^{[k-1]} + u_{n-2}^{[k]} z^{[k-2]} + \ldots + u_{n-k}^{[k]} z^{[0]})$$

$$= u_{n-1}^{[k][n-k]} z^{[n-1]} + u_{n-2}^{[k][n-k]} z^{[n-2]} + \ldots + u_{n-k}^{[k][n-k]} z^{[n-k]}$$

$$= u_{n-1}^{[n]} z^{[n-1]} + u_{n-2}^{[n]} z^{[n-2]} + \ldots + u_{n-k}^{[n]} z^{[n-k]}$$

$$= u_{n-1} z^{[n-1]} + u_{n-2} z^{[n-2]} + \ldots + u_{n-k} z^{[n-k]},$$

where $u'(z)$ is nothing but the remainder obtained when $z^{[k]} * f(z)$ is divided on the right by $G(z)$.

Thus the class of $[n,k,d]_{q^n}$ q-Cyclic MRD codes for $n - k \geq k$ are invertible.

An example is given in the following which demonstrates the invertible property of the $[5,2,4]_{2^5}$ 2-Cyclic MRD code.

***Example 4.3.3.1 :*** Consider the $[5,2,4]_{2^5}$ 2-Cyclic MRD code generated by $G(z) = z^8 + \alpha^{10} z^4 + \alpha^{17} z^2 + \alpha^{13} z$, where $\alpha$ is a root of the primitive polynomial $x^5 + x^2 + 1$ over GF(2).

The generator matrix G and the parity-check matrix H for the $[5,2,4]_{2^5}$ 2-Cyclic MRD code in systematic form are respectively given by

$$G = \begin{bmatrix} \alpha^{13} & \alpha^{17} & \alpha^{10} & 1 & 0 \\ \alpha^{2} & \alpha^{14} & \alpha^{9} & 0 & 1 \end{bmatrix}$$

and

$$H = \begin{bmatrix} 1 & 0 & 0 & \alpha^{13} & \alpha^{2} \\ 0 & 1 & 0 & \alpha^{17} & \alpha^{14} \\ 0 & 0 & 1 & \alpha^{10} & \alpha^{9} \end{bmatrix}.$$

Let $u_1(z) = z^{16} + \alpha z^8$ be the message polynomial to be encoded.

Dividing $u_1(z)$ on the right by $G(z)$, one obtains

$$z^{16} + \alpha z^8 = (z^2 + \alpha^{12} z) * G(z) + f_1(z),$$

where $f_1(z) = \alpha^{14} z^4 + \alpha^{24} z^2 + \alpha^{25} z$.

Dividing $z^{[2]} * f_1(z)$ on the right by G(z), one obtains

$$\alpha^{25}z^{16} + \alpha^3z^8 + \alpha^7z^4 = (\alpha^{25}z^2 + \alpha^{22}z) + G(z) + f_1'(z),$$

where $f_1'(z) = z^2 + \alpha^4z$. Then the message polynomial $u_1(z)$ is $z^{[3]} * (z^2 + \alpha^4z) = z^{16} + \alpha z^8$.

Thus the message polynomial $u_1(z) = z^{16} + \alpha z^8$ is obtained only with the knowledge of the parity-check polynomial $f_1(z) = \alpha^{14}z^4 + \alpha^{24}z^2 + \alpha^{25}z$.

Similarly, the message polynomial $u_2(z) = z^8$ can be obtained from the parity-check polynomial
$f_2(z) = \alpha^{10}z^4 + \alpha^{17}z^2 + \alpha^{13}z$, as the remainder obtained by dividing $z^{[2]} * f_2(z)$ on the right by G(z) is $f_2'(z) = z$ and thus $u_2(z) = z^{[3]} * f_2'(z)$.

The case when n–k < k, the q-Cyclic MRD codes are not invertible, but the class of shortened $[n-t_0, k-t_0]_{q^n}$ q-Cyclic MRD codes, $2k-n \le t_0 < k$ has the invertible property which is proved in what follows.

**Case (ii) n – k < k**

The case when n – k < k, the class of $[n,k,d]_{q^n}$ q-Cyclic MRD codes are not invertible. But the codes obtained by shortening the $[n,k,d]_{q^n}$ q-Cyclic MRD codes by removing the $t_0$ leading zero information symbols are invertible, where $2k - n \le t_0 < k$; i.e., the shortened $[n-t_0, k-t_0]_{q^n}$ q-Cyclic MRD codes are invertible.

Following example shows the face that an $[n,k,d]_{q^n}$ q-Cyclic MRD code with n – k < k is not invertible.

***Example 4.3.3.2:*** Consider GF $(3^3)$ = $\{0, 1, a, \ldots, \alpha^{3^3-2}\}$, where $\alpha$ is a root of the primitive polynomial $x^3 + 2x + 1$ over GF(3). Let C = $[3,2,2]_{3^3}$ be the 3-Cyclic MRD code with the generator polynomial $G(z) = z^3 + \alpha^{21} z$.

The generator matrix G and the parity-check matrix H for the $[3,2,2]_{3^3}$ 3-Cyclic MRD code in systematic form are respectively given by

$$G = \begin{bmatrix} \alpha^8 & 1 & 0 \\ \alpha^6 & 0 & 1 \end{bmatrix} \text{ and } H = [1 \; -\alpha^8 \; -\alpha^6].$$

Let $u(z) = \alpha^{22}z^9 + \alpha^7z^3$ be a message polynomial.

Dividing u(z) on the right by G(z), one obtains

$$\alpha^{22}z^9 + \alpha^7z^3 = \alpha^{22}z^3 * G(z) + f(z),$$

where $f(z) \equiv 0$, the zero polynomial.

Then dividing $z^{[2]} * f(z)$ on the right by G(z), one obtains the remainder polynomial $u'(z) \equiv 0$ so that $z^{[1]} * u'(z) \equiv 0$; but $u(z) = \alpha^{22}z^9 + \alpha^7z^3$, which is not the zero polynomial.

Thus from the parity-check polynomial f(z), the message polynomial u(z) is not retrieved. This is because the parity block is same for the codewords (0, 0, 0) and $(0, \alpha^7, \alpha^{22})$.

Thus an $[n,k,d]_{q^n}$ q-Cyclic MRD code for which $n - k < k$ need not be invertible.

Using the shortening technique to the class of $[n,k,d]_{q^n}$ q-Cyclic MRD codes discussed in subsection 4.3.2, one can obtain a class of invertible shortened $[n-t_0, k-t_0]_{q^n}$ q-Cyclic MRD codes as follows, $2k - n \le t_0 < k$.

Let C denote an $[n,k,d]_{q^n}$ q-Cyclic MRD code generated by G(z) = $\sum_{i=0}^{n-k} G_i z^{[i]}$ such that n – k < k. Let H(z) = $\sum_{i=0}^{k} H_i z^{[i]}$ be such that $z^{[n]} - z = H(z) * G(z)$. Let $t_0$ be such that $2k - n \le t_0 < k$.

Consider the shortened $[n - t_0, k - t_0]_{q^n}$ q-Cyclic MRD code $C_{(t_0)}$ of length $n - t_0$ and dimension $k - t_0$. The encoding scheme, similar to the $[n,k,d]_{q^n}$ q-Cyclic MRD code C, for the shortened $[n - t_0, k - t_0]_{q^n}$ q-Cyclic MRD code $C_{(t_0)}$ is given below.

Let u(z) = $u_{n-t_0-1} z^{[(n-t_0)-1]} + u_{n-t_0-2} z^{[(n-t_0)-2]} + \ldots + u_{n-k} z^{[n-k]}$ be a message polynomial to be encoded. Divide u(z) on the right by G(z):

$$u(z) = q(z) * G(z) + f(z), \deg (f(z)) < [n - k], \qquad (4.18)$$

where q(z) and f(z) are respectively the quotient and the remainder polynomials.

Then g(z) = u(z) – f(z) is the code polynomial corresponding to the message polynomial u(z).

Now one needs to ensure that no two codewords in the shortened $[n - t_0, k - t_0]_{q^n}$ q-Cyclic MRD code $C_{(t_0)}$ have same parity block. This is proved in the following theorem. Though the proof is similar to the case (i), here it is given for the sake of completeness.

**THEOREM 4.3.3.2:** *Let C be an $[n,k,d]_{q^n}$ q-Cyclic MRD code such that n – k < k. Let $C_{(t_0)}$ be the shortened $[n - t_0, k - t_0]_{q^n}$ q-Cyclic MRD code, $2k - n \le t_0 < k$. Then no two codewords of $C_{(t_0)}$ will have same parity-check symbols.*

***Proof:*** Consider the generator polynomial G(z) of C in the form: $G(z) = z^{[n-k]} + G_{n-k-1} z^{[n-k+1]} + \ldots + G_0 z^{[0]}$.

Let $u^{(1)}(z) = \sum_{i=1}^{k-t_0} u_i^{(1)} z^{[n-t_0-i]}$ and

$$u^{(2)}(z) = \sum_{i=1}^{k-t_0} u_i^{(2)} z^{[n-t_0-i]}$$

be two distinct message polynomials to be encoded.

Dividing $u^{(1)}(z)$ and $u^{(2)}(z)$ on the right by G(z),

$$u^{(1)}(z) = q_1(z) * G(z) + f_1(z), \quad \deg(f_1(z)) < [n-k] \qquad [4.19]$$

$$\text{and } u^{(2)}(z) = q_2(z) * G(z) + f_2(z), \quad \deg(f_2(z)) < [n-k] \qquad [4.20]$$

The code polynomials corresponding to $u^{(1)}(z)$ and $u^{(2)}(z)$ respectively are,

$$g_1(z) = u^{(1)}(z) - f_1(z)$$

and $g_2(z) = u^{(2)}(z) - f_2(z)$.

Suppose that $f_1(z) = f_2(z)$.

Subtracting (4.20) from (4.19),

$$u^{(1)}(z) - u^{(2)}(z) = (q_1(z) - q_2(z)) * G(z)$$

i.e., $z^{[n-k]} * (u'^{(1)}(z) - u'^{(2)}(z)) = (q_1(z) - q_2(z)) * G(z)$,

where $u'^{(1)}$, $u'^{(2)}$ are such $u^{(1)}(z) = z^{[n-k]} * u'^{(1)}(z)$ and

$u^{(2)}(z) = z^{[n-k]} * u'^{(2)}(z)$.

This shows that G(z) is a right divisor of
$z^{[n-k]} * (u'^{(1)}(z) - u'^{(2)}(z))$.

Since $z^{[n-k]}$ is relatively prime to G(z), G(z) is a right divisor of $u'^{(1)}(z) - u'^{(2)}(z)$.

However this is impossible, because $u'^{(1)}(z) - u'^{(2)}(z) \neq 0$ and its degree is less than $[k-t_0-1]$ but the degree of G(z) is [n – k] which is strictly great than $[k-t_0-1]$.

Thus $f_1(z) \neq f_2(z)$. Hence the theorem.

Since the remainder f(z) resulting from dividing u(z) by G(z) is unique, the Theorem 4.3.3.2 implies that there exists an one-to-one correspondence between a message polynomial u(z) and its parity-check polynomial f(z). Therefore knowing only the parity-check polynomial f(z) alone, the message polynomial u(z) can be determined uniquely.

In what follows, it is shown how to recover the message polynomial u(z) from its parity-check polynomial f(z). Now one has the equation (4.18):

$$u(z) = q(z) * G(z) + f(z), \deg(f(z)) < [n-k].$$

Consider $z^{[k]} * f(z)$ :

$$\begin{aligned}
z^{[k]} * f(z) &= z^{[k]} * (-q(z) * G(z) + u(z)) \\
&= z^{[k]} * (-q(z) * G(z)) + z^{[k]} * u(z) \\
&= z^{[k]} * (-q(z) * G(z)) + (u'(z) * z^{[k]}) * z^{[n-k]} \\
&= z^{[k]} * (-q(z) * G(z)) + u'(z) * z^{[n]} \\
&= z^{[k]} * (-q(z) * G(z)) + u'(z) * (z^{[n]} - z) + u'(z) * z \\
&= (z^{[k]} * (-q(z))) * G(z) + u'(z) * (H(z)*G(z)) + u'(z)
\end{aligned}$$

(since G(z) is a right divisor of $z^{[n]} - z$).

$$= (z^{[k]} * (-q(z)) + u'(z) * H(z)) * G(z) + u'(z), \qquad (4.21)$$

where $u'(z) = u^{[k]}_{n-t_0-1} z^{[k-t_0-1]} + u^{[k]}_{n-t_0-2} z^{[k-t_0-2]} + \ldots + u^{[k]}_{n-k} z^{[0]}$.

See that from (4.21), the message polynomial u(z) is nothing but $z^{[n-k]} * u'(z)$;

i.e., $u(z) = z^{[n-k]} * u'(z)$

$$= z^{[n-k]} *( u^{[k]}_{n-t_0-1} z^{[k-t_0-1]} + u^{[k]}_{n-t_0-2} z^{[k-t_0-2]} + \ldots + u^{[k]}_{n-k} z^{[0]})$$

$$= u^{[k][n-k]}_{n-t_0-1} z^{[n-t_0-1]} + u^{[k][n-k]}_{n-t_0-2} z^{[n-t_0-2]} + \ldots + u^{[k][n-k]}_{n-k} z^{[n-k]}$$

$$= u^{[n]}_{n-t_0-1} z^{[n-t_0-1]} + u^{[n]}_{n-t_0-2} z^{[n-t_0-2]} + \ldots + u^{[n]}_{n-k} z^{[n-k]}$$

$$= u_{n-t_0-1} z^{[n-t_0-1]} + u_{n-t_0-2} z^{[n-t_0-2]} + \ldots + u_{n-k} z^{[n-k]},$$

where $u'(z)$ is nothing but the remainder obtained when $z^{[k]} * f(z)$ is divided on the right by G(z).

Having carried out a comprehensive study on a characteristic of RD codes, namely the invertible property of the class of RD codes, the next section deals yet another characteristic enjoyed by the class of RD codes: Rank Distance codes having complementary duals.

## 4.4 Rank Distance Codes With Complementary Duals

It is known that an (n, k) F-ary linear code is just a k-dimensional subspace of the n-dimensional vector space $F^n$ of n-tuples with coordinates in the finite field F. Recall that the vectors u and v in $F^n$ are said to be orthogonal if $\langle u, v \rangle = 0$ and that, if $\Gamma$ is an (n, k) F-ary linear code, the dual code $\Gamma^{\perp}$ is the (n, n–k) F-ary linear code consisting of all vectors $v \in F^n$ that are orthogonal to every vector u in $\Gamma$. The class of Linear codes with Complementary Duals (LCD codes) is defined by J.L. Massey [19] in 1992. A F-ary linear code $\Gamma$ is called an LCD code if $\Gamma \cap \Gamma^{\perp} = \{0\}$. It is immediate that $\Gamma$ is an LCD code just when the following occurs

$$F^n = \Gamma \oplus \Gamma^{\perp} ;$$

that is, when $F^n$ is the direct sum of $\Gamma$ and $\Gamma^{\perp}$.

The following theorem due to J.L. Massey [19] gives an algebraic characterization to the class of LCD codes.

**THEOREM 4.4.1:** *[19] Let $\Gamma$ be an (n, k) F-ary linear code with generator matrix G. Then $\Gamma$ is an LCD code if and only if the $k \times k$ matrix $GG^T$ is non-singular. Further, if $\Gamma$ is in LCD code, then $\Pi_\Gamma = G^T (GG^T)^{-1} G$ is the orthogonal projector from $F^n$ onto $\Gamma$.*

Recall that the trace function tr : $GF(2^n) \to GF(2)$ is defined by $tr(a) = \sum_{j=0}^{n-1} a^{2j}$. Let $B = \{b_1, b_2, \ldots, b_n\}$ be an ordered basis of $GF(2^n)$ over GF(2) and $\{\overline{b}_1, \overline{b}_2, ..., \overline{b}_n\}$ be the dual basis of B with respect to the trace function, i.e.,

$tr(b_i\overline{b}_j) = \delta_{ij}$ for $1 \le i, j \le n$, where $\delta_{ij}$ is the Kronecker's symbol. Then B is called trace-orthogonal if $b_i = \overline{b}_i$, $1 \le i \le n$. The dual basis of any basis of GF $(2^n)$ over GF (2) is determined uniquely [20, Theorem 4.1.1]. It is known that a trace-orthogonal basis of $GF(2^n)$ over GF(2) exists for any positive integer n [23].

This section carries out a brief study on the class of $[n,k,d]_{2^n}$ Rank Distance codes having complementary duals.

Consider the generator matrix of an $[n,k,d]_{2^n}$ MRD code:

$$G = \begin{bmatrix} \alpha_1 & \alpha_2 & \cdots & \alpha_n \\ \alpha_1^{[1]} & \alpha_2^{[1]} & \cdots & \alpha_n^{[1]} \\ \vdots & \vdots & \ddots & \vdots \\ \alpha_1^{[k-1]} & \alpha_2^{[k-1]} & \cdots & \alpha_n^{[k-1]} \end{bmatrix} \qquad (4.22)$$

where $\alpha_1, \alpha_2, \ldots, \alpha_n$ are linearly independent over GF(2).

Represent the generator matrix G defined as in (4.22) by $G = \left[\alpha_j^{2i}\right]_{i,j=0,1}^{k-1,n}$ where $\alpha_1, \alpha_2, \ldots, \alpha_n$ are linearly independent GF(2). This book calls the matrix $G = \left[\alpha_j^{2i}\right]_{i,j=0,1}^{k-1,n}$ with the first row entries $\alpha_1, \alpha_2, \ldots, \alpha_n$ being a trace-orthogonal basis in $GF(2^n)$ as the trace-orthogonal-generator matrix. Clearly, G generates an $[n,k,d]_{2^n}$ MRD code.

This section is divided into three subsections. The subsection 1 proves that the class of $[n,k,d]_{2^n}$ MRD codes generated by the trace-orthogonal-generator matrices are LCD codes. Description to the (noiseless and noisy) 2-user F-Adder Channel is given in subsection 2. The final subsection gives the coding for the noiseless 2-user F-Adder Channel via the class of $[n,k,d]_{2^n}$ MRD codes having complementary duals and describes a coding problem with the noisy 2-user F-Adder Channel.

### 4.4.1 MRD Codes with Complementary Duals

It is observed that an $[n,k,d]_{q^n}$ RD code, $n \le N$ with $d < n - k + 1$ need not be an LCD code but it is interesting to see that the class of $[n,k,d]_{2^n}$ MRD codes generated by the trace-orthogonal-generator matrices are LCD codes.

The following counter example shows that an $[n,k,d]_{q^n}$ RD code, $n \le N$ with $d < n - k + 1$ need not be an LCD code.

***Example 4.4.1.1:*** Let GF $(2^3) = \{0, 1, \alpha, \alpha^2, \alpha^3, \alpha^4, \alpha^5, \alpha^6\}$, where $\alpha$ is a root of the primitive polynomial $x^3 + x + 1$ over GF(2).

Let $\Gamma = [3,2,1]_{2^3}$ be the 2-Cyclic RD code with the generator matrix G corresponding to the generator polynomial
$G(Z) = z^2 + \alpha z$ is given as

$$G = \begin{bmatrix} \alpha & 1 & 0 \\ 0 & \alpha^2 & 1 \end{bmatrix}.$$

Clearly $|GG^T| = (0)$. Thus the $[3,2,1]_{2^3}$ 2-Cyclic RD code is not an LCD code.

The following theorem proves that an $[n,k,d]_{2^n}$ MRD code generated by a trace-orthogonal-generator matrix is an LCD code.

**THEOREM 4.4.1.1:** *An $[n,k,d]_{2^n}$ MRD code generated by $G = \left[\alpha_j^{2i}\right]_{i,j=0,1}^{k-1,n}$ with $\{\alpha_1, \alpha_2, \ldots, \alpha_n\}$ being a trace-orthogonal basis in GF($2^n$) is an LCD code.*

Let $\Gamma$ denote an $[n,k,d]_{2^n}$ MRD code generated by

$$G = \begin{bmatrix} \alpha_1 & \alpha_2 & \cdots & \alpha_n \\ \alpha_1^{[1]} & \alpha_2^{[1]} & \cdots & \alpha_n^{[1]} \\ \vdots & \vdots & \ddots & \vdots \\ \alpha_1^{[k-1]} & \alpha_2^{[k-1]} & \cdots & \alpha_n^{[k-1]} \end{bmatrix},$$

where $\{\alpha_1, \alpha_2, \ldots, \alpha_n\}$ is a trace-orthogonal basis in GF($2^n$).

In order to prove that $\Gamma$ is an LCD code, one has to prove that the $k \times k$ matrix $GG^T$ is non-singular.

Since $\{\alpha_1, \alpha_2, \ldots, \alpha_n\}$ is a trace-orthogonal basis in $GF(2^n)$, the k row vectors $(\alpha_1^{2i}, \alpha_2^{2i}, \ldots, \alpha_n^{2i})$, i = 0, 1, …, k–1 are orthonormal vectors. It follows that $GG^T = I$, where I denotes the identity matrix. Thus, $GG^T$ is non-singular. Therefore, $\Gamma$ is an LCD code.

Further, the orthogonal projector $\prod_\Gamma$ from $[GF(2^n)]^n$ onto $\Gamma$ defined by $r\prod_\Gamma = rG^T(GG^T)^{-1}G$ for each $r \in [GF(2^n)]^n$ exists.

In the above theorem, it is proved that the class of $[n,k,d]_{2^n}$ MRD codes generated by the generator matrices of the form $G = \left[\alpha_j^{2^i}\right]_{i,j=0,1}^{k-1,n}$ with $\{\alpha_1, \alpha_2, \ldots, \alpha_n\}$ being a trace-orthogonal basis in $GF(2^n)$ are LCD codes. But there exists $[n,k,d]_{q^N}$ MRD codes, $n \leq N$ generated by generator matrices of the form $G = \left[\alpha_j^{q^i}\right]_{i,j=0,1}^{k-1,n}$ with the first row entries $\alpha_1, \alpha_2, \ldots, \alpha_n$ are linearly independent over GF(q) are LCD codes. This is evident from the following example.

***Example 4.4.1.2:*** Let $\Gamma = [3,1,3]_{3^4}$ be the MRD code with the generator matrix G :

$$G = [\alpha^4 \ \alpha^{65} \ 1],$$

where $\alpha$ is a root of the primitive polynomial $x^4 + x + 2$ over GF(3). Clearly $|GG^T| \neq (0)$.

The class of $[n,k,d]_{2^n}$ MRD codes having complementary duals is effective at coding over the noiseless 2-user F-Adder Channel, which can be seen in the next subsection.

### 4.4.2 The 2-user F-Adder Channel

Given a finite field F, the F-Adder Channel is described as the channel whose inputs are elements of F and the output is the sum (over F) of the inputs [25]. This section describes 2-user F-Adder Channel as the F-Adder Channel shared by 2 users. Following are the descriptions to the 2-user F-Adder Channel, for both the noiseless and noisy cases.

**Case (i) Noiseless 2-user F-Adder Channel**

A pictorial representation of the noiseless 2-user F-Adder Channel is depicted in the following figure. In this noiseless communication channel, the two users of the 2-user F-Adder Channel transmit two n-tuples, say $\gamma_1$ and $\gamma_2$ respectively from the F-ary linear codes $\Gamma_1$ and $\Gamma_2$. Then the received vector, say r is the componentwise sum $\gamma_1 + \gamma_2$ over the finite field F.

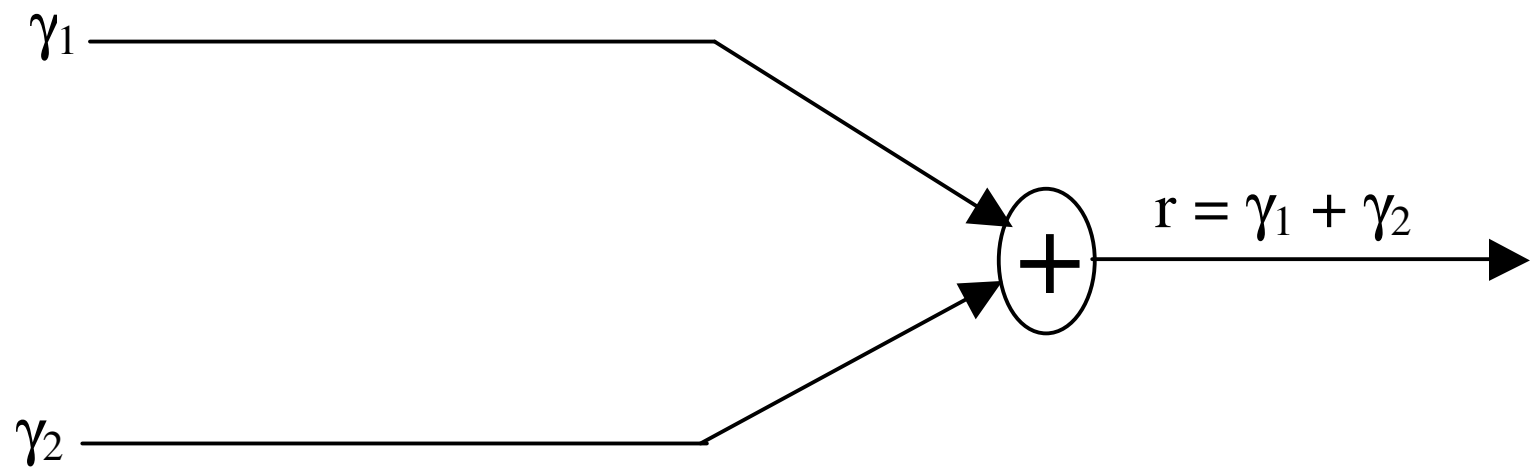

**Noiseless 2-user F-Adder Channel**

In this noiseless channel, the problem for the receiver is to decode the received vector $r = \gamma_1 + \gamma_2$ into the codewords $\gamma_1$ and $\gamma_2$ originally transmitted. The next subsection provides a solution to this problem through the class of $[n,k,d]_{2^n}$ MRD codes having complementary duals.

**Case (ii) Noisy 2-user F-Adder Channel**

Consider the pictorial representation of the noisy 2-user F-Adder Channel depicted in the following figure in which the two users are attempting to transmit two

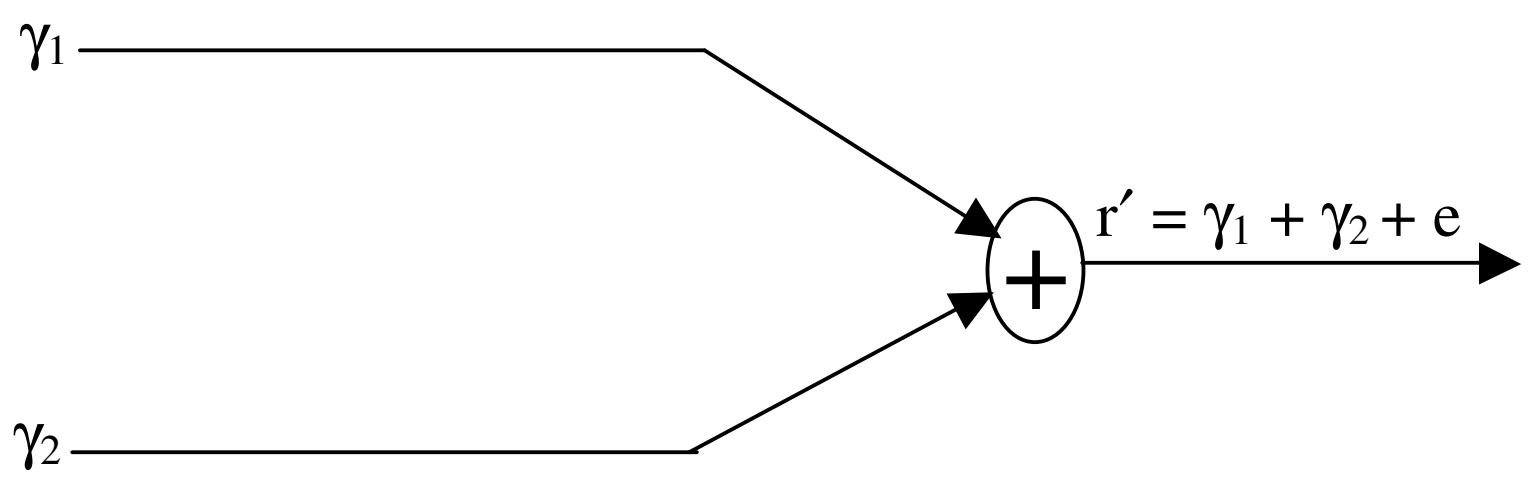

**Noisy 2-user F-Adder Channel**

codewords, say $\gamma_1$ and $\gamma_2$ respectively from the F-ary linear codes $\Gamma_1$ and $\Gamma_2$. Then, in this noisy channel, the received vector, say $r'$ is the componentwise sum $\gamma_1 + \gamma_2 + e$ over F, where e is an error-vector. The problem for the receiver is to decode the received vector $r' = \gamma_1 + \gamma_2 + e$ into the transmitted codewords $\gamma_1$ and $\gamma_2$.

**4.4.3 Coding for the 2-user F-Adder Channel**

As described in the previous subsection, the 2-user F-Adder Channel is a F-Adder Channel shared by 2 users. This subsection describes how the class of $[n,k,d]_{2^n}$ MRD codes having complementary duals can be effectively used over the noiseless 2-user $GF(2^n)$ - Adder Channel. For the case when the 2-user F-Adder Channel is noisy, the coding problem is described. It is observed that the class of LCD codes is not suitable for coding over the noisy 2-user F-Adder Channel.

**Case (i) Noiseless 2-user F-Adder Channel**

Consider the pictorial representation of the noiseless 2-user F-Adder Channel given in figure for F = GF $(2^n)$. Let Γ be an $[n,k]_{2^n}$ MRD code generated by a trace-orthogonal-generator matrix and $\Gamma^{\perp} = [n, n-k]_{2^n}$ denote its dual code. Assume that the two users of the noiseless 2-user GF$(2^n)$-Adder Channel transmit the codewords γ and β from Γ and $\Gamma^{\perp}$ respectively. Then, in this noiseless 2-user GF$(2^n)$ - Adder Channel, the received vector r is the componentwise sum $\gamma + \beta$ over GF$(2^n)$.

Since Γ is an LCD code, the orthogonal projector $\prod_{\Gamma} = GT\ (GG^T)^{-1}\ G$ defined from $[GF(2^n)]^n$ onto Γ exists. To receive the codewords γ and β, the receiver simply applies the orthogonal projector $\prod_{\Gamma}$ on r which gives $r\prod_{\Gamma} = \gamma$. The codeword β is then obtained by subtracting γ from r = γ + β. Thus the codewords γ and β transmitted are retrieved from the received vector r successfully.

The following example describes the coding for the noiseless 2-user GF$(3^4)$-Adder Channel via $[3,1,3]_{3^4}$ MRD code.

***Example 4.4.3.1:*** Let $\Gamma = [3,1,3]_{3^4}$ be the 3-Cyclic MRD code with the generator matrix G and parity-check matrix H.

$$G = [\alpha^4\ \alpha^{65}\ 1] \text{ and } H = \begin{bmatrix} 1 & \alpha & \alpha^2 \\ 1 & \alpha^3 & \alpha^6 \end{bmatrix},$$

where α is a root of the primitive polynomial $x^4 + x + 2$ over GF(3). Note that the generator matrix for the dual code $\Gamma^{\perp}$ is H.

Since Γ being an LCD code, the orthogonal projector $\prod_{\Gamma}$ defined from $[GF\ (3^4)]^4$ onto Γ is given by

$$\prod_\Gamma = \begin{bmatrix} \alpha^{41} & \alpha^{22} & \alpha^{37} \\ \alpha^{22} & \alpha^{3} & \alpha^{18} \\ \alpha^{37} & \alpha^{18} & \alpha^{33} \end{bmatrix}.$$

Suppose the codewords $\gamma_1 = (1, \alpha^{61}, \alpha^{76})$ and $\gamma_2 = (\alpha^2, \alpha^5, \alpha^8)$ from $\Gamma$ and $\Gamma^\perp$ respectively be transmitted over the noiseless 2-user GF $(3^4)$-Adder Channel. Then $r = \gamma_1 + \gamma_2 = (\alpha^{24}, \alpha^{78}, \alpha^{35})$ is the received vector.

Applying the orthogonal projector $\prod_\Gamma$ on r:

$$\gamma_1 = (\alpha^{24}, \alpha^{78}, \alpha^{35}) \prod_\Gamma = (1, \alpha^{61}, \alpha^{76}).$$

Then $\gamma_2 = r - (1, \alpha^{61}, \alpha^{76}) = (\alpha^2, \alpha^5, \alpha^8)$. Hence the codewords $\gamma_1 = (1, \alpha^{61}, \alpha^{76})$ and $\gamma_2 = (\alpha^2, \alpha^5, \alpha^8)$ are retrieved from $r = (\alpha^{24}, \alpha^{78}, \alpha^{35})$.

**Case (ii) Noisy 2-user F-Adder Channel**

Consider the pictorial representation of the noisy 2-user F-Adder Channel depicted in figure.

In this noisy channel, the received vector $r'$ may not always be the sum of the transmitted codewords $\gamma_1$ and $\gamma_2$, but may be a sum of the transmitted codewords along errors, i.e., $r' = \gamma_1 + \gamma_2 + e$, where e is an error-vector. Let $r[e; q] \leq \left\lfloor \frac{d-1}{2} \right\rfloor$ with d being a positive integer.

The problem here for the decoder is to employ an error-correcting decoding scheme to decode the received vector $r'$. If one employs the coding scheme described for the noiseless case s/he would not recover the transmitted codewords $\gamma_1$ and $\gamma_2$ from the erroneously received vector $r'$. Suppose $\Gamma$ is an LCD code such that $\gamma_1 \in \Gamma$, $\gamma_2 \in \Gamma^\perp$, then $\gamma_1 + \gamma_2 \in F^n = \Gamma \oplus \Gamma^\perp$ so that error correction is not possible in the erroneously received vector $r' = \gamma_1 + \gamma_2 + e$. Consider a situation wherein $\gamma_1 + \gamma_2$ is a codeword of an (n, k, d) F-ary linear code with the parity-check matrix H. Then employing the associated error-correcting decoding technique and using the syndrome $r'H^T$, one can decode the received vector r′ into $\gamma_1$ and $\gamma_2$.

If this is the situation, one can correct upto $\left\lfloor \frac{d-1}{2} \right\rfloor$ errors in the erroneously received vector r′. So, the prime motivation is to construct multiuser error-correcting codes that can be employed over the noisy 2-user F-Adder Channel effectively.

Chapter Five

# EFFECTIVE ERASURE CODES FOR RELIABLE COMPUTER COMMUNICATION PROTOCOLS USING CODES OVER ARBITRARY INTEGER RINGS

Linear algebraic codes can be defined using symbols chosen from a set of arbitrary size. However, most of the results of coding theory have been derived assuming that the code symbols are elements of a finite field especially the finite field $Z_2 = \{0, 1\}$. Recently linear codes over integer rings have raised a great interest for their role in algebraic coding theory and their successful applications in combined coding and modulation. [2, 3] has constructed cyclic codes over $Z_m$ (the ring of integers modulo m) where m is an interger of the form $p_1, p_2, \ldots, p_k$ where $p_i$'s are distinct primes, from cyclic codes over $Z_{p_i}$. In his later paper, Blake derived parity check matrices for codes over $Z_m$ analogous to Hamming codes and Reed-Solomon codes. Calderbank and Sloane [4], Priti Shanker [23], J.C. Interlando et al [13] extended the notion of cyclic codes, Reed-Solomon codes and BCH codes over GF(q) to class of codes over finite rings $Z_q$, with q a power of a prime. Most studies in algebraic coding theory deal with Hamming metric. Since the Hamming metric is not always well matched to the characteristics of the real channels EM. Gabidulin [9] introduced a new metric, called the rank metric, and he called the codes equipped with this metric as rank distance codes.

In this chapter we study the codes with rank metric over the ring of integers modulo 2m, where $2m = 2p_1, p_2, \ldots, p_t$, where $p_i$'s are distinct primes. These codes are proved to be better than the codes with Hamming metric as they can cater to the complex and unpredictable situations in the communicating channels. Further, they are found to have a better error correcting capability.

This chapter is divided into two sections. In the first section we define a new class of codes, called the integer rank distance codes using the ring of integers modulo 2m. In the second section we define the Maximum Integer Rank Distance codes (MIRD codes) using the main result - Singleton - style bound for integer rank distance codes. Further we give the method of coding and decoding algorithm for this new class of MIRD codes without which the construction will not be complete.

## 5.1. Integer Rank Distance Code

As said earlier $Z_{2m}$ is the ring of integers modulo 2m where $2m = 2p_1 p_2 \ldots p_t$, where $p_i$'s are distinct primes. Let $Z_{2m}[x]$ be the ring of polynomials in the indeterminate x. Let $p(x) \in Z_{2m}[x]$ be an irreducible polynomial of degree n over $Z_{2m}$. Let V be $Z_{2m}[x] / \langle p(x) \rangle$, where $\langle p(x) \rangle$ denotes the ideal generated by p (x). Clearly V is a module of dimension n over $Z_{2m}$. Any element $x \in V$ can be represented by $x = (x_1, x_2, \ldots, x_n)$ where $x_i \in Z_{2m}$ i.e., a polynomial is regarded as a n-tuple. The elements of $Z_{2m}$ can be treated as polynomials in 2 over $Z_2$. Hence an element $x_1 \in Z_{2m}$ has representation as a $N_0$-tuple $\left(x_{1i}, x_{2i}, \ldots, x_{N_{0i}}\right)$, $x_{ij} \in Z_2$. Hence, with each $x \in V$ we have an associated matrix,

$$m(x) = \begin{pmatrix} x_{11} & x_{12} & \cdots & x_{1n} \\ x_{21} & x_{22} & \cdots & x_{2n} \\ \vdots & \vdots & & \vdots \\ x_{N_0 1} & x_{N_0 2} & \cdots & x_{N_0 n} \end{pmatrix}$$

where the i-th column represents the i-th coordinate '$x_i$' of 'x' over $Z_2$.

**DEFINITION 5.1.1:** *Rank of an element $x \in V$ is defined as the rank of the matrix m(x) over $Z_2$.*

*Let n(x) denote the rank of the m(x). Then it is clear that the function $x \to r_I(x)$ is a norm on V. We call this as the integer rank norm and denote this by $r_I(x)$. The metric induced by the integer rank norm is defined as the integer rank metric on V. If $x, y \in V$, then the integer rank distance between x and y is $d_I(x, y) = r_I(x + y)$.*

We illustrate this by an example.

***Example 5.1.1:*** Let $V = Z_6[x] / \langle x^3 + 1 \rangle$. Then, V is a module over $Z_6$ and the elements of $Z_6$ can be treated as polynomials in 2 over $Z_2$. Let $x = (3, 5, 2) \in V$.

Then, $m(x) = \begin{pmatrix} 0 & 1 & 0 \\ 1 & 0 & 1 \\ 1 & 1 & 0 \end{pmatrix}$ over $Z_2$.

Clearly rank of m(x) is 3.

**DEFINITION 5.1.2:** *The module equipped with this integer rank metric, for convenience we call it as the integer rank distance space.*

**DEFINITION 5.1.3:** *A linear (n, k) integer rank distance code is a linear submodule of dimension k in the integer rank distance space V.*

*By C(n, k), we denote a linear (n, k) integer rank distance code. Let d = min {$d_I$ (x, y) : x, y ∈ C(n, k)}. This d is called the minimum integer rank distance of the integer rank distance code C(n, k). Now C(n, k) is a linear (n, k, d) integer rank distance code.*

**DEFINITION 5.1.4:** *A generator matrix G of a linear (n, k) integer rank distance code C(n, k) is a k × n matrix over $Z_{2m}$ whose rows form a basis for C(n, k). We can reduce the generator matrix G to the form G = [$I_k$, $A_{k, n-k}$] where $I_k$ is the k × k identity matrix and $A_{k, n-k}$ is some k × (n–k) matrix over $Z_{2m}$.*

**DEFINITION 5.1.5:** *If G is a generator matrix of a linear (n, k) integer rank distance code C[n, k], then a matrix H of order (n – k) × n over $Z_{2m}$ such that $GH^T$ = {0} is called a parity check matrix of C[n, k].*

H can be reduced to the form $H = \left[-A^{T}_{(n-k)\times k}, I_{n-k}\right]$.

We define a linear (n, k) integer rank distance code as :

1. The linear submodule generated by the rows of the generator matrix (or)
2. the solution space of the parity check matrix.

***Note:*** The usual Hamming distance between any two vectors, x, y ∈ V, which is the number of places in which x and y differ also induces a norm on each vector x ∈ V, called the Hamming norm, which is denoted throughout this paper as $r_H(x)$.

## 5.2. Maximum Integer Rank Distance Codes

In this section we mainly arrive at a singleton - style bound for integer rank distance codes, using which we define Maximum Integer Rank Distance codes, which we refer by MIRD codes here after. To achieve this we essentially make use of the following lemma. Further, we give a characterization theorem for finding the maximum integer rank distance of a integer rank distance code.

***Lemma 5.2.1***: Let $n_1$ and $n_2$ be any two norms defined on V with $n_1(x) \leq n_2(x)$, for all $x \in V$. Let $V_1(n, d)$ and $V_2(n, d)$ be the volumes of codes with maximum distance d with norms $n_1$ and $n_2$, respectively. Then, $V_1(n, d) \leq V_2(n, d)$.

***Proof :*** Since $n_1(x) \leq n_2(x)$, for all $x \in V$, any $y \in V_1(n, d)$ will also be in $V_2$ (n, d).

Hence $V_1(n, d) \leq V_2(n, d)$.

**THEOREM 5.2.1:** *(Singleton-style bound for integer rank distance codes) Any linear (n, k, d) integer rank distance code satisfies the inequality $d \leq n - k + 1$.*

***Proof :*** Choosing $n_1 = r_1(x)$ and $n_2 = r_H(x)$, we have $r_1(x) \leq r_H(x)$, we have by Lemma 5.2.1 and the singleton bound for codes over integer rings with Hamming metric, $d \leq n-k+1$.

**DEFINITION 5.2.1:** *Codes which attain equality in the singleton-style bound are called Maximum Integer Rank Distance codes (MIRD codes).*

The following theorem is used to find the minimum integer rank distance of the linear (n, k) integer rank distance code.

**THEOREM 5.2.2:** *Let C be a linear (n, k) integer rank distance code with generator matrix G and parity check matrix H. Then, C has rank distance d if and only if for any $(d-1) \times n$ matrix M of rank d–1 with elements from $Z_2$,*

$$r(MH^T; Z_{2m}) = d - 1 \qquad (5.2.1)$$

*and then there exists a $d \times n$ matrix $M_1$ of rank d with elements from $Z_2$ for which*

$$r(M_1 H^T; Z_{2m}) < d \qquad (5.2.2)$$

***Proof :*** Let $x = (x_1, x_2, \ldots, x_n)$ be any vector in V such that $r_1(x) \leq d$. Then, there exists a matrix $M_1$ such that $x = zM_1$, where $z = (z_1, z_2, \ldots, z_d)$, $z_i \in Z^{2m}$, $i = 1, 2, \ldots, d$ and $M_1$ is a $d \times n$ matrix of rank d with elements from $Z_2$. Assume that code C contains a codeword x with rank d. Then, $xH^T = (0)$, which implies $zM_1H^T = (0)$. Consequently, $r(M_1H^T; Z_{2m}) < d$. Since the codes has minimum distance d for any $(d - 1) \times n$ matrix N of rank $d - 1$ with elements from $Z_2$, the equation

$$(z_1, z_2, \ldots, z_{d-1})\, MH^T = (0) \qquad (5.2.3)$$

should have only a trivial solution i.e.,

$$r(MH^T; Z_{2m}) < d.$$

Sufficient part is direct from the definition.

**THEOREM 5.2.3:** *A code C is a linear MIRD (n, k) code if and only if for any $(n–k) \times n$ matrix M of rank $n – k$ with elements from $Z_2$;*

$$r(MH^T; Z_{2m}) = n–k \qquad (5.2.4)$$

***Proof:*** From theorem 5.2.1 we get $d \leq n - k + 1$. By theorem 5.2.1, $n – k + 1 \geq d$. Therefore, $d = n – k + 1$.

We define a class of MIRD codes of length $n \leq N_0$. These codes are analogous to the generalized Reed-Solomon codes over integer rings.

We introduce the notation $[i] = 2^i$; $i = 0, \pm 1, \pm 2, \ldots$

Assume that $h_i \in Z_{2m}$, $i = 1, 2, \ldots, n$ and assume that these elements are linearly independent over $Z_2$. Given the designed distance $d < n$, we generate the matrix

$$H = \begin{bmatrix} h_1 & h_2 & \cdots & h_n \\ h_1^{[1]} & h_2^{[1]} & \cdots & h_n^{[1]} \\ \vdots & \vdots & \ddots & \vdots \\ h_1^{[d-2]} & h_2^{[d-2]} & \cdots & h_n^{[d-2]} \end{bmatrix} \qquad (5.2.5)$$

**THEOREM 5.2.4:** *The linear (n, k) integer rank distance code C with parity check matrix H is a MIRD code of length n and minimum integer rank distance d.*

***Proof:*** By theorem 5.2.3 it is sufficient to show that for any $(d-1) \times n$ matrix M of rank d–1 with elements from $Z_2$, we have $r(HM^T ; Z_{2m}) = d - 1$. The square matrix $HM^T$ has the form,

$$HM^T = \begin{bmatrix} f_1 & f_2 & \cdots & f_n \\ f_1^{[1]} & f_2^{[1]} & \cdots & f_n^{[1]} \\ \vdots & \vdots & \ddots & \vdots \\ f_1^{[d-2]} & f_2^{[d-2]} & \cdots & f_n^{[d-2]} \end{bmatrix} \qquad (5.2.6)$$

where $(f_1, f_2, \ldots, f_{d-1}) = (h_1, h_2, \ldots, h_{d-1}) M^T$. The elements $f_1, \ldots, f_{d-1} \in Z_{2m}$ are linearly independent over $Z_2$, since otherwise, $h_1, h_2, \ldots, h_{d-1}$ would also be linearly dependent, in contradiction to our assumption. Clearly, $HM^T$ being a Vander Monde matrix is non singular, i.e., $r(HM^T ; Z_{2m}) = d - 1$.

**THEOREM 5.2.5:** *Let C be the code with parity check matrix H. Then, generating matrix G has the form,*

$$G = \begin{bmatrix} g_1 & g_2 & \cdots & g_n \\ g_1^{[1]} & g_2^{[1]} & \cdots & g_n^{[1]} \\ \vdots & \vdots & \ddots & \vdots \\ g_1^{[k-2]} & g_2^{[k-2]} & \cdots & g_n^{[k-2]} \end{bmatrix} \qquad (5.2.7)$$

*where $k = n - d + 1$ and the elements $g_1, g_2, \ldots, g_n$ are linearly independent over $Z_2$.*

***Proof:*** From theorem 5.2.4 with d = n–1, there exist elements $f_1$, $f_2$, …, $f_n \in Z_{2m}$ which are linearly independent over $Z_2$ and satisfy

$$\sum_{i=0}^{n} f_1 h_i^{[s]}, \quad s = 0, 1, \ldots, n-2. \tag{5.2.8}$$

Since $f_1$, $f_2$, …, $f_n$ are linearly independent over $Z_2$, $f_1^{[-k+1]}, f_2^{[-k+1]}, \ldots, f_n^{[-k+1]}$ are also linearly independent over $Z_2$. We take $g_1 = f_1^{[-k+1]}$, $g_2 = f_2^{[-k+1]}$, …, $g_n = f_n^{[-k+1]}$ to be the first row of the matrix (5.2.7).

As polynomials with coefficients from $Z_{2m}$ play an important role in the theory of maximum distance separable codes over $Z_{2m}$ the linearized polynomials play similar role in the theory of maximum rank distance codes over $GF(2^N)$ [29, 32]. In a analogous way we can use linearized polynomial with coefficients from $Z_{2m}$ to study MIRD codes over $Z_{2m}$.

A linearized polynomial is one of the form $F(z) = \sum_{i=0}^{n} f_1 z^{[1]}$, $[i] = 2^i$, where $f_i \in Z_{2m}$. We define the sums of the polynomials as

$$\begin{aligned} F(z) + G(z) &= \sum_{i=0}^{n} f_1 z^{[1]} + \sum_{i=0}^{n} g_1 z^{[1]} \\ &= \sum_{i=0} R_{2m}(f_i + g_i) z^{[i]} \end{aligned}$$

where $R_{2m}(f_i + g_i)$ is the least nonnegative integer when $f_i + g_i$ is divided by the integer 2m. The multiplication product is the symbolic product F * G = F (G(z)). This multiplication operation is non-commutative. The set of all linearized polynomials with coefficients from $Z_{2m}$ form a non-commutative ring with identity element $f_0(z) = z$.

Generalized Inversionless Euclidean algorithm for division (whether left or right) of one polynomial by another exists in this ring. Here we consider right division.

Let $F_0(z)$ and $F_1(z)$ be two linearlized polynomials with $\deg F_1(z) < \deg F_0(z)$. Then by generalized Euclidean algorithm we get a sequential chain of equalities

$$\begin{aligned}
F_0(z) &= G_1(z) * F_1(z) + F_2(z),\ \deg F_2(z) < \deg F_1(z) \\
F_1(z) &= G_2(z) * F_2(z) + F_3(z),\ \deg F_3(z) < \deg F_2(z) \\
\vdots \quad & \qquad\quad \vdots \qquad\qquad\qquad \vdots \\
F_{s-1}(z) &= G_s(z) * F_s(z) + F_{s+1}(z),\ \deg F_{s+1}(z) < \deg F_s(z) \\
F_s(z) &= G_{s+1}(z) * F_{s+1}(z)
\end{aligned} \qquad (5.2.9)$$

The last non zero remainder $F_{s+1}(z)$ in this chain is the right symbolic LCD of polynomials $F_0(z)$ and $F_1(z)$. If we introduce polynomials $U_i(z)$, $A_i(z)$, $V_i(z)$ and $B_i(z)$, defined recursively for $i \geq 1$, by

$$\begin{aligned}
U_i(z) &= U_{i-1}(z) * G_i(z) + U_{i-2}(z),\ U_0(z) = z,\ U_{-1}(z) = 0 \\
A_i(z) &= G_i(z) * A_{i-1}(z) + A_{i-1}(z),\ A_0(z) = z,\ A_{-1}(z) = 0 \\
V_i(z) &= V_{i-1}(z) * G_i(z) + V_{i-1}(z),\ V_0(z) = z,\ V_{-1}(z) = z \\
B_i(z) &= B_{i-1}(z) * G_i(z) + B_{i-2}(z),\ B_0(z) = z,\ B_{-1}(z) = z
\end{aligned} \qquad (5.2.10)$$

then,

$$\begin{aligned}
F_0(z) &= U_i(z) * F_i(z) + U_{i-1}(z) * F_{i+1}(z) \\
F_1(z) &= V_i(z) * F_i(z) + V_{i-1}(z) * F_{i+1}(z).
\end{aligned} \qquad (5.2.11)$$

In addition,

$$F_i(z) = (-1)^i (B_{i-1}(z) * F_0(z) - A_{i-1}(z) * F_i(z)). \qquad (5.2.12)$$

Consider the factor ring $R_{N_0}$ which is the ring of the linearlized polynomials over $Z_{2m}$ modulo $z^{[N_0]} - z$. The elements of this factor ring are also linearlized polynomials of degree $\leq [N_0] - 1$.

Let $F(z) = \sum_{i=0}^{N_0-1} f_1 z^{[i]} \in R_{N_0}$. Raising the polynomial to the power 2 we get $F^{[1]}(z) = f_{N_0-1}^{[1]} z^{[0]} + f_0^{[1]} z^{[1]} + \ldots + f_{N_0-2}^{[1]} z^{[N_0-1]}$ which is equivalent to raising all its coefficients to the power 2 and then performing a cyclical shift. This operation will be called as a 2-cyclical shift.

The ideals in $R_{N_0}$ are principal ideals and are generated by polynomials G(z) that are right divisors of $z^{[N_0]} - z$ i.e., the polynomial G(z) is such that $z^{[N_0]} - z = H(z) * G(z)$. The ideal {G(z)} is invariant under 2-cyclical shift.

The codes with generator matrix of the form (5.2.7) can be described in terms of the linearlized polynomials. Assume that $g_1, g_2, \ldots, g_n$ are specified elements that are linearly independent over $Z_2$. Then all vectors of the form

$g = (F(g_1), F(g_2), \ldots, F(g_n))$, where F(z) extends over all linearlized polynomials of degree less than or equal to $[k-1] = 2^{k-1}$ with coefficients from $Z_{2m}$ are codewords.

Code C is called 2-cyclical if a 2-cyclical shift of any codeword is also a codeword i.e., if $(g_0, g_1, \ldots, g_{n-1}) \in C$ then its 2-cyclic shift $g_{n-1}^{[1]}, g_0^{[1]}, \ldots, g_{n-2}^{[1]} \in C$. This is analogous to the usual cyclic codes over the integer ring.

Let us consider, for simplicity, codes with length $n = N_o$.

Assume that $G(z) = \sum_{i=0}^{r} G_i z^{[i]}$ is a right zero divisor of $z^{[N_0]} - z$. Then, its 2-cyclical code consists of all polynomials of the form $c(z) * G(z)$, where c(z) is an arbitrary linearlized polynomial of degree $\le (N_0 - r - 1)$. The dimension of the code is $k = N_0 - r$. Its generator matrix has the form

$$G = \begin{pmatrix} G_0 & G_1 & \cdots & G_r & 0 & 0 & \cdots & 0 \\ 0 & G_0^{[1]} & \cdots & G_{r-1}^{[1]} & G_r^{[1]} & 0 & \cdots & 0 \\ \vdots & \vdots & & \vdots & \vdots & \vdots & & \vdots \\ 0 & 0 & \cdots & G_0^{[k-1]} & G_1^{[k-1]} & G_2^{[k-1]} & \cdots & G_r^{[k-1]} \end{pmatrix}$$

(5.2.13)

If $z^{[N_0]} - z = G(z) * H(z)$, where G(z) is the generator polynomial then H(z) can be taken as the corresponding check polynomial.

We note that an element in V is a code vector if and only if the corresponding linearlized polynomial can be divided without remainder by the generating polynomial G(z).

In otherwords, an element g is a code vector if and only if the corresponding polynomial g(z) is such that

$$g(z) * H(z) \equiv 0 \bmod (z^{[N_0]} - z).$$

If $H(z) = \sum_{i=0}^{k} H_i z^{[i]}$ is a check polynomial, then the check matrix has the form

$$H = \begin{pmatrix} H_k & H_{k-1}^{[1]} & \cdots & H_0^{[k]} & 0 & 0 & \cdots & 0 \\ 0 & H_k^{[1]} & \cdots & H_1^{[k]} & H_0^{[k-1]} & 0 & \cdots & 0 \\ \vdots & \vdots & & \vdots & \vdots & \vdots & & \vdots \\ 0 & 0 & \cdots & 0 & H_k^{[r-1]} & H_{k-1}^{[r]} & \cdots & H_0^{[N_0-1]} \end{pmatrix}$$

(5.2.14)

$$= \begin{pmatrix} h_0 & h_1 & \cdots & h_k & 0 & 0 & \cdots & 0 \\ 0 & h_0^{[1]} & \cdots & h_{k-1}^{[1]} & h_k^{[1]} & 0 & \cdots & 0 \\ \vdots & \vdots & & \vdots & \vdots & \vdots & & \vdots \\ 0 & 0 & \cdots & h_0^{[r-1]} & h_1^{[r-1]} & h_2^{[r-1]} & \cdots & h_k^{[r-1]} \end{pmatrix}$$

where $H_i^{[k-1]} = h_{k-i}$.

Now we proceed onto give the coding and decoding techniques of MIRD codes.

As in the case of cyclic codes over the integer ring, systematic coding can be effected either by means of a check polynomial or by means of a generating polynomial.

If $H(z) = \sum_{i=0}^{k} H_i z^{[i]}$ is a check polynomial, then each

$g(z) = \sum_{i=0}^{N_0-1} g_i z^{[i]}$ satisfies g(z) * H(z) = 0, which gives

$$\sum_{j=0}^{k} g_{N_0-i-j+1} H_j^{[N_0-i-j+1]} = 0, \; i = 0, 1, \ldots, N_0 - 1. \qquad (5.2.15)$$

If we assume that $g_{N_0-1}, \ldots, g_{N_0-k}$ are information symbols then we can determine the check symbols $g_{N_0-k-1}, \ldots, g_0$.

Assume that we are given the generating polynomial

$G(z) = \sum_{i=0}^{N_0-1} G_i z^{[i]}$. We divide the polynomial

$G_0(z) = g_{N_0-1} z^{N_0-1} + \ldots + g_{N_0-k} z^{N_0-k}$ on the right by G(z) to get

$$G_0(z) = Q(z) * G(z) + R(z), \text{ deg. } R(z) < [N_0 - k] \qquad (5.2.16)$$

The coefficients $g_{N_0-1}$ for degrees $[N_0 - i]$, $i = 1, \ldots, N_0$ of the remainder are the check symbols.

MIRD codes with check matrix (5.2.5) can be decoded using an algorithm that is similar to the algorithm for maximum rank codes over fields [32] with appropriate modifications due to the presence of zero divisors. Let $g = (g_1, g_2, \ldots, g_n)$ be the code vector. $e = (e_1, e_2, \ldots, e_n)$ be the error vector and $y = g + e$ be the received vector. We first calculate the syndrome

$$s = (s_0, s_1, \ldots, s_{d-2}) = yH^T = eH^T. \quad (5.2.17)$$

The decoders problem is to determine the error vector e on the basis of the known syndrome vector s. Assume that the rank norm of the error vector is m. Thus we have

$$e = EN = (E_0, E_1, \ldots, E_m)N \quad (5.2.18)$$

where $E_0, E_1, \ldots, E_m$ are linearly independent over $Z_2$ and $N = (N_{ij})$ is a $m \times n$ matrix of rank m with elements from $Z_2$. Thus (5.2.17) can be written as

$$s = ENH^T = EX \quad (5.2.19)$$

where the matrix $X = NH^T$ has the form,

$$X = \begin{bmatrix} x_1 & x_1^{[1]} & \ldots & x_1^{[d-2]} \\ x_2 & x_2^{[1]} & \ldots & x_2^{[d-2]} \\ \vdots & \vdots & \ddots & \vdots \\ x_m & x_m^{[1]} & \ldots & x_m^{[d-2]} \end{bmatrix}$$

where $x_p = \sum_{j=1}^{n} N_{pj} h_j$, $p = 1, 2, \ldots, m$ (5.2.20)

are linearly independent over $Z_2$. (5.2.19) is equivalent to the system of equations in the unknowns $E_0, E_1, \ldots, E_m, x_1, x_2, \ldots, x_m$.

$$\sum_{i=1}^{m} E_i x_i^{[p]} = s_p,\ p = 1, 2, \ldots, d-2. \quad (5.2.21)$$

Assuming that the solution of this system has been found, we can determine the matrix N and the error vector e from (5.2.20) and (5.2.18). For $m \le (d-1)/2$ all the solutions lead to the error vector e. Thus, the decoding problem reduces to the solution of system (5.2.20) for the smallest possible value of m. The solution to the system (5.2.21) is unique if and only if all the error magnitudes $E_0, E_1, \ldots, E_m$ are units analogous to the case of the maximum distance codes over the integer rings [32].

We introduce the polynomial $s(z) = \sum_{j=0}^{d-2} s_j z^{[j]}$, corresponding to the syndrome s. Assume that $\Lambda(z) = \sum_{p=0}^{m} \Lambda_p z^{[p]}$, $\Lambda_m = 1$, denotes a polynomial whose roots are all possible linear combinations of $E_0, E_1, \ldots, E_m$ with coefficients from $Z_2$. Let $F(z) = \sum_{i=0}^{m-1} F_i z^{[i]}$, where $F_i = \sum_{p=0}^{i} \Lambda_p s_{i-p}^{[p]}$, $i = 0, 1, \ldots, m-1$.

We have the equality

$$F(z) = \Lambda (z) * s (z) \bmod z^{[d-1]} \qquad (5.2.22)$$

Indeed,

$$\Lambda (z) * s (z) = \sum_{p=0}^{m} \Lambda_p (s(z))^{[p]}$$

$$= \sum_{i=0}^{m+d-1} z^{[i]} \sum_{p+j-i} \Lambda_p s_j^{[p]} .$$

But for $m \leq i \leq d-2$, we have

$$\sum_{p+j-i} \Lambda_p s_j^{[p]} = \sum_{p=0}^{m} \Lambda_p s_{i-p}^{[p]}$$

$$= \sum_{p=0}^{m} \Lambda_p \left( \sum_{j=1}^{m} E_j x_j^{[j-p]} \right)^{[p]}$$

$$= \sum_{j=1}^{m} x_j^{[i]} \Lambda\left(E_j\right) = 0,$$

since $\Lambda (E_j) = 0$, $j = 1, 2, \ldots, m$.

If the coefficients of polynomial F(z) are known, then the coefficients of polynomial Λ (z).

$$\Lambda_0 = \frac{F_j}{s_j}$$

$$\Lambda_p = \frac{F_{j+p} - \sum_{i=0}^{p-1} \Lambda_i s_{p-j-i}^{[i]}}{s_j^{[p]}}, \ p = 1, 2, \ldots \qquad (5.2.23)$$

where for $j + p \geq m$ we set $F_{j+p} = 0$.

Now we assume that $E_0$, $E_1$, …, $E_m$ as well as the coefficients of $\Lambda(z)$, are known. We consider the following truncated system in the unknows:

$$\sum_{j=1}^{m} E_j x_j^{[p]} = s_p, \ p = 1, 2, \ldots, m-1. \qquad (5.2.24)$$

We will solve (5.2.23) using the method of successive elimination of variables. We set $A_{ij} = E_j$, $Q_{1p} = s_p$; we multiply the (p+1)-th equation of the system by $A_{11}$, we extract the root of degree 2, and we subtract the p-th equation. As a result we obtain a system that does not contain $x_1$ :

$$\sum_{j=1}^{p} A_{1j} x_j^{[p]} = Q_{1p}, \ p = 0, 1, 2, \ldots, m-2 \qquad (5.2.25)$$

where

$$A_{1j} = A_{1j} = \left(\frac{A_{1j}}{A_{11}}\right)^{[-1]} A_{11}, \ j = 2, \ldots, m.$$

$$Q_{1p} = A_{1p} - \left(\frac{Q_{1p+1}}{A_{11}}\right)^{[-1]} A_{11}, \ p = 0, 1, \ldots, m-2. \qquad (5.2.26)$$

Repeating this process m – 1 times, and retaining the first equations obtained from the systems at each step, we arrive at a system of linear equations with a upper triangular coefficient matrix:

$$\sum_{j=1}^{m} A_{1j}x_j = Q_{i0}, \; p = 0, 1, 2, \ldots, m. \qquad (5.2.27)$$

where

$$A_{1j} = E_j, \; j = 1, 2, \ldots, m$$

$$A_{1j} = \begin{cases} 0 & \text{if } j < i \\ A_{(i-1)j} - \left(\dfrac{A_{(i-1)j}}{A_{(i-1)(i-1)}}\right)^{[-1]} A_{(i-1)(i-1)} & \text{if } j \geq i; \; i = 2, 3, \ldots, m \end{cases}$$

$$(5.2.28)$$

$$Q_{1p} = s_p, \; p = 0, 1, \ldots, m-1$$

$$Q_{1j} = Q_{(i-1)p} \left(\frac{Q_{(i-1)(p+1)}}{Q_{(i-1)(i-1)}}\right)^{[-1]} A_{(i-1)(i-1)} \; p = 0, 1, \ldots, m-i;$$

$$i = 2, \ldots, m \qquad (5.2.29)$$

The solution of (5.2.25) can be found by back substitution i.e.,

$$x_m = \frac{Q_{m0}}{A_{mm}}$$

$$x_{m-i} = \frac{Q_{(m-i)i} - \sum_{j=m-i+1}^{m} A_{(m-i)j}x_j}{A_{(m-i)(m-i)}}, \; i = 1, 2, \ldots, m-1. \; (5.2.30)$$

We now state the above in the form of a decoding algorithm.

**Step 1:** We calculate the syndrome s = $(s_0, s_1, \ldots, s_{d-2})$ and the corresponding polynomial $s(z) = \sum_{j=0}^{d-2} s_j z^{[j]}$.

**Step 2:** We set $F_0(z) = z^{[d-1]}$, $F_1(z) = s(z)$ and employ generalized Euclidean algorithm until we reach a $F_{1+1}(z)$ such that deg $F_1(z) \geq 2^{(d-1)/2}$, deg $F_{1+1}(z) < 2^{(d-1)/2}$. Then

$$\Delta(z) = \gamma A_m(z)$$

$$F(z) = \gamma (-1)^m F_{m+1}(z) \qquad (5.2.31)$$

where, $\gamma$ is chosen such that the coefficient of $\Lambda_m$ is equal to 1.

Polynomial $\Lambda(z)$ can be determined either on the basis of the first formula in (5.2.22), if polynomials $A_i(z)$, i = 1, 2, …, are calculated in parallel in the course of Euclidean algorithm, or using (5.2.23), which employ the coefficients of the remainder $F_{m+1}(z)$ calculated in the course of the algorithm. Then roots $E_0, E_1, \ldots, E_m$ of $\Lambda(z)$ that are linearly independent over $Z_2$ are determined.

**Step 3:** Using 5.2.27 to 5.2.30 the known $E_0, E_1, \ldots, E_m$ are used as a basis for determining $x_1, x_2, \ldots, x_m$. Representing these quantities in the form (5.2.20), we can obtain matrix N. Finally, we calculate the error vector using (5.2.29).

***Example 5.2.1:*** The code with the following parity check matrix $H = \begin{pmatrix} 1 & 4 & 2 \\ 1 & 2 & 4 \end{pmatrix}$ has length n = 3 and the designed distance d = 3. Let y = (3, 2, 1) be the received codeword. Then, the syndrome is $s = (s_0, s_1) = yH^T = (5, 5)$. Hence, by modified inversionless Euclid's algorithm, $\Lambda(z) = -\ (s_0^{[1]}\ /\ s_1) + z^{[1]} = -5z + z^2$. The nonzero root of $\Lambda(z)$ is 5. Hence, E = 5. From the single equation of the system (5.2.26) we determine

$x = (s_1 / s_0) = 1 = (1 \times 1) + (0 \times 4) + (0 \times 2) = y_1\, h_1 + y_2\, h_2 + y_2\, h_2)$ which gives $(y_1, y_2, y_3) = (1, 0, 0)$. The error vector is $e = (y_1\, E, y_2E, y_3E) = (5, 0, 0)$

We have constructed a new class of effective erasure codes over the ring of integers using the rank distance metric. This class of codes can be suited in different applications like reliable computer communication protocols, ARQ protocols in satellites communications and in amateur radios.

**Chapter Six**

# CONCATENATION OF ALGEBRAIC CODES

In this chapter we study concatenation of linear block code with another linear block code or concatenation of three or more linear block codes. We describe this in section one. This chapter has two sections. In section two we describe concatenation of RD codes with CR-matric. We give the probable ways of building concatenated bicodes, biconcatination of codes and quasi concatenated bicodes.

## 6.1 Concatenation of Linear Block Codes

The construction and the encoding procedure of the concatenated codes is described in this section which is as follows:

Let $S = ((n_1, n_2, \ldots, n_t), (k_1, \ldots, k_t))$ be the special supercode defined over the Galois field $Z_2 = GF(2)$ and the inner codes $C_1 = C_1(n_1, k_1)$, $C_2 = C_2(n_2, k_2)$, … and $C_t = C_t(n_t, k_t)$ be codes with Hamming metric defined over the Galois field GF(2).

Suppose $m = (x_1 \mid x_2 \mid \ldots \mid x_t)$ be the message to be encoded where each $x_i$ is a $k_i$ tuple of the form $(a^i_1, a^i_2, \ldots, a^i_{k_i})$; $1 \leq i \leq t$, thus we get $x = (b_1 \mid b_2 \mid \ldots \mid b_t)$ where $b_i = (a^i_1, a^i_2, \ldots, a^i_{n_i})$; $1 \leq i \leq t$ where $a^i_j \in GF(2)$; $1 \leq j \leq n_i$; $1 \leq i \leq t$.

Thus we have a $k_1 + k_2 + \ldots + k_t$ dimensional super vector space associated with this concatenated code or the super code is a $k_1 + k_2 + \ldots + k_t$ dimensional vector subspace of a $n_1 + n_2 + \ldots + n_t$ dimensional vector space over $Z_2$. Infact S has $2^{k_1+k_2+\ldots+k_t}$ elements.

We can consider a map;

$g : S \to (C_1 \times C_2 \times \ldots \times C_t) = \{(x_1, x_2, \ldots, x_t) \mid x_i \in C_i,$ $1 \leq i \leq t\}$ given by $g\,((x_1 \mid \ldots \mid x_t)) = (x_1, x_2, \ldots, x_t)$ where $x_i \in C_i$; $1 \leq i \leq t$. This map g is one to one linear transformation of S to $C_1 \times C_2 \times \ldots \times C_t$.

So $S \cong (C_1 \times C_2 \times \ldots \times C_t)$ just '|' replaced by $\times$.

The concatenated code is comprised of the outer code and t number of inner codes. The outer code is a super code and the t-inner codes are the row submatrices of the outer code. That is the outer code which is the super code has its sub row vectors from these t inner codes. The t inner codes can be distinct or otherwise. The following is the concatenated coding system.

Thus we see the concatenated code word which is a super code takes the form of a super row matrix $x = (x_1^1, x_2^1, \ldots, x_{n_1}^1 \mid x_1^2, x_2^2, \ldots, x_{n_2}^2 \mid \ldots \mid x_1^t, x_2^t, \ldots, x_{n_t}^t)$ with $x_p^j \in Z_2$; $1 \leq j, p \leq t$ and the concatenated code

$S = \{x = (x_1^1, x_2^1, \ldots, x_{n_1}^1 \mid x_1^2, x_2^2, \ldots, x_{n_2}^2 \mid \ldots \mid x_1^t, x_2^t, \ldots, x_{n_t}^t \mid x_1^j, x_2^j, \ldots, x_{n_j}^j) \in C_j = C_j(n_j, k_j);\ 1 \leq j \leq t;\ x_p^j \in Z_2;\ 1 \leq j, p \leq t\}$.

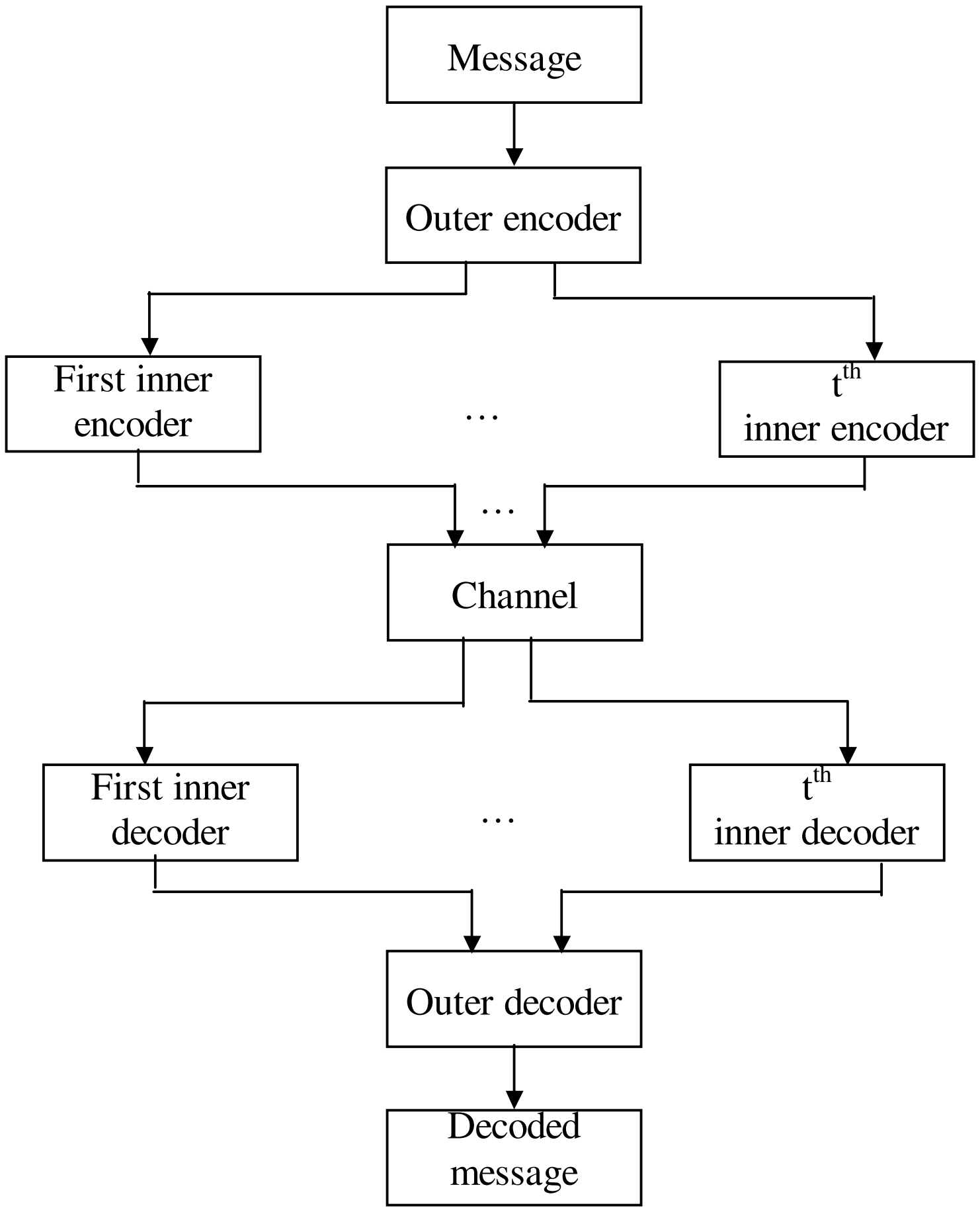

For any two x, y ∈ S we can define $d_S(x, y) = \sum_{j=1}^{t} d_j(x_j, y_j)$ where $x_j, y_j \in C_j(n_j, k_j)$; $1 \le j \le t$. It is easily verified $d_S(x, y) = d_S(x+y)$ denotes the super Hamming metric on S over GF(2) or t-concatenated Hamming metric on S. All properties of Hamming metric can be derived in case of t-concatenated Hamming metric for the concatenated code S.

We illustrate this by an example.

***Example 6.1.1:*** Let S = {($x_1$ | $x_2$ | $x_3$ | $x_4$) | $x_1 \in$ C (2, 3), $x_2 \in$ C(5, 2), $x_3 \in$ C (7, 4) and $x_4 \in$ C (4, 2)} be a concatenated code with the codes $C_1$ = C (8, 3), $C_2$ = (5, 2), $C_3$ = C (7, 4) and $C_4$ = C (4, 2).

We see error detection can be done using the parity check matrices of $C_1$, $C_2$, $C_3$ and $C_4$.

Error correction can be carried out using coset leader method [16]. This type of concatenated super codes has the following advantages.

(i) The rate of transmission is increased.

(ii) This concatenated codes saves time during transmission as insteaded of sending t-vectors of each $n_t$ tuples we send a lengthened code.

(iii) We can retrieve the messages.

Further $\sum_{i=0}^{t}\left\lfloor \frac{d_{min}-1}{2} \right\rfloor_i \leq \left\lfloor \sum_{t} \frac{d_{min}-1}{2} \right\rfloor$.

So in a Super code which is the concatenated code using linear codes with Hamming metric we can detect and correct errors, with in the limits of compatibility.

Now we describe another type of concatenated code of a linear codes with Hamming metric.

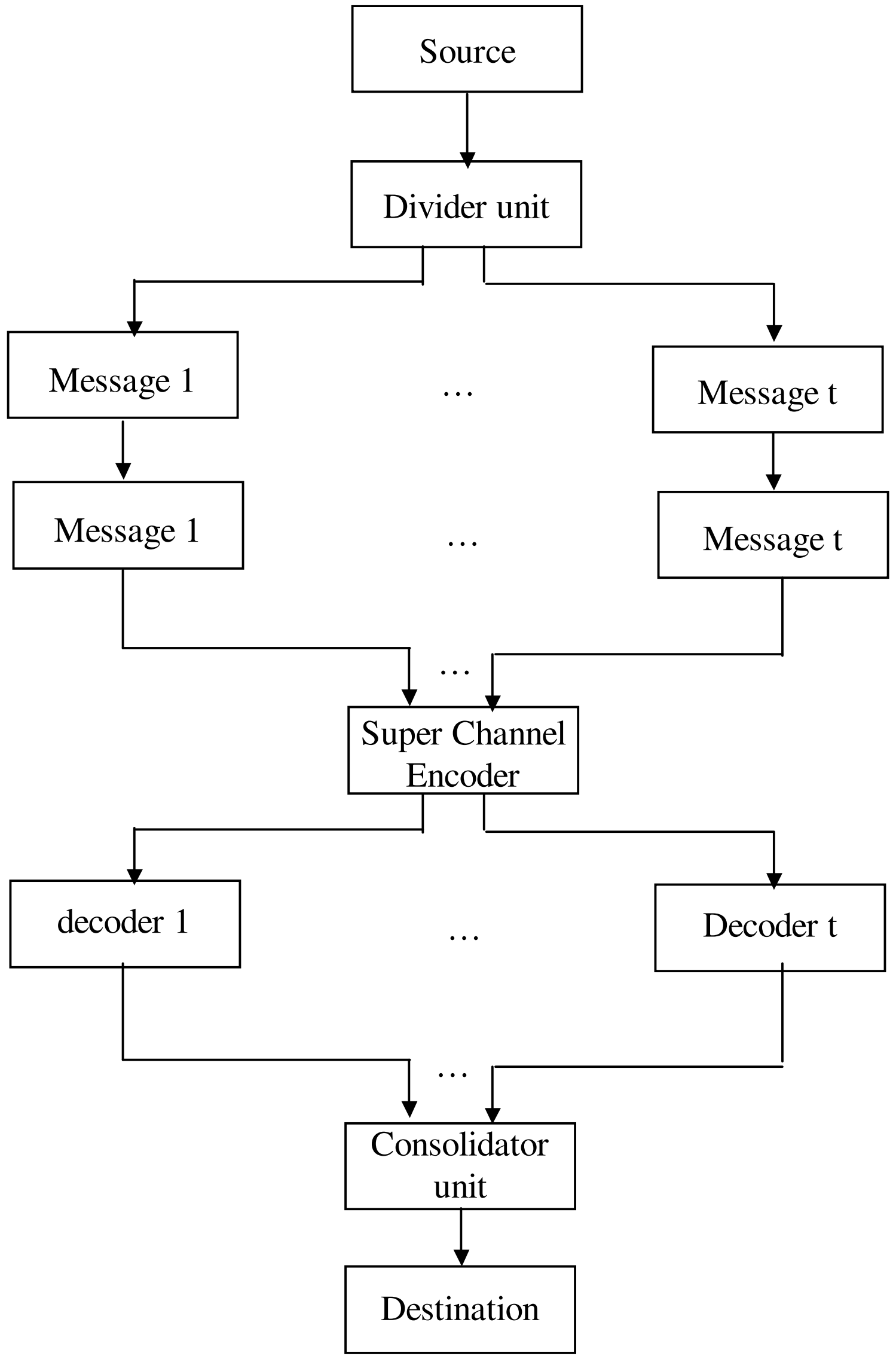

The message x = $(x_1^1, x_2^1, \ldots, x_{n_1}^1,\ x_1^2, x_2^2, \ldots, x_{n_2}^2, \ldots x_1^t, x_2^t, \ldots, x_{n_t}^t)$ is sent to the divider unit which devides them into t messages and the messages are sent to t encoder and they pass the super channel and reach as t decoders, from the decoder units reaches the consolidation unit and then the destination. These codes will be useful when the bulk work is distributed to different units and the result is consolidated to get the final message.

Let $C_1, C_2, \ldots, C_t$ be t number of linear codes distinct or otherwise of length $n_1, n_2, \ldots, n_t$ and messages $k_1, k_2, \ldots, k_t$ respectively. We denote by P = $(x_1^1, x_2^1, \ldots, x_{n_1}^1,\ x_1^2, x_2^2, \ldots, x_{n_2}^2, \ldots x_1^t, x_2^t, \ldots, x_{n_t}^t) \mid x_i^j \in Z_2;\ 1 \leq j \leq t$ and $1 \leq i \leq n_1, n_2, \ldots, n_t\}$. We say P is concatenated using t inner codes, however P is not a super code, just a code. By using usual transmission the rate of transmission is certainly less or equal to the sum of the rates of transmissions.

The diagram is self explanatory.

We can use error connection technique or erasure correction technique.

Finally we give yet another type of concatenated linear codes which we will define in the following. These will also be known as special blank concatenated codes.

Suppose $C_1 = C(n, k)$ code and $C_2 = C(n, k')$ or $C(n–1, k_1)$ code. Then we form the concatenated code $C = C_1\ C_2 = \{(x_1y_1\ x_2y_2 \ldots x_n\ y_n) \mid (x_1, x_2, \ldots, x_n)$ is in C (n, k) and $(y_1, y_2, \ldots, y_n) \in C(n, k')\}$ $((x_1\ y_1\ x_2\ y_2 \ldots x_{n-1}\ y_{n-1}\ x_n)$ where $(x_1, \ldots, x_n) \in C(n, k)$ and $(y_1, \ldots, y_{n-1}) \in C\ (n-1, k_1))$ will be defined as the alternate concatenated two code or special blank concatenated code and the transmission is done using the 'special blank'; however 'special blanks' are not erasures.

We will describes how the transmission takes place.

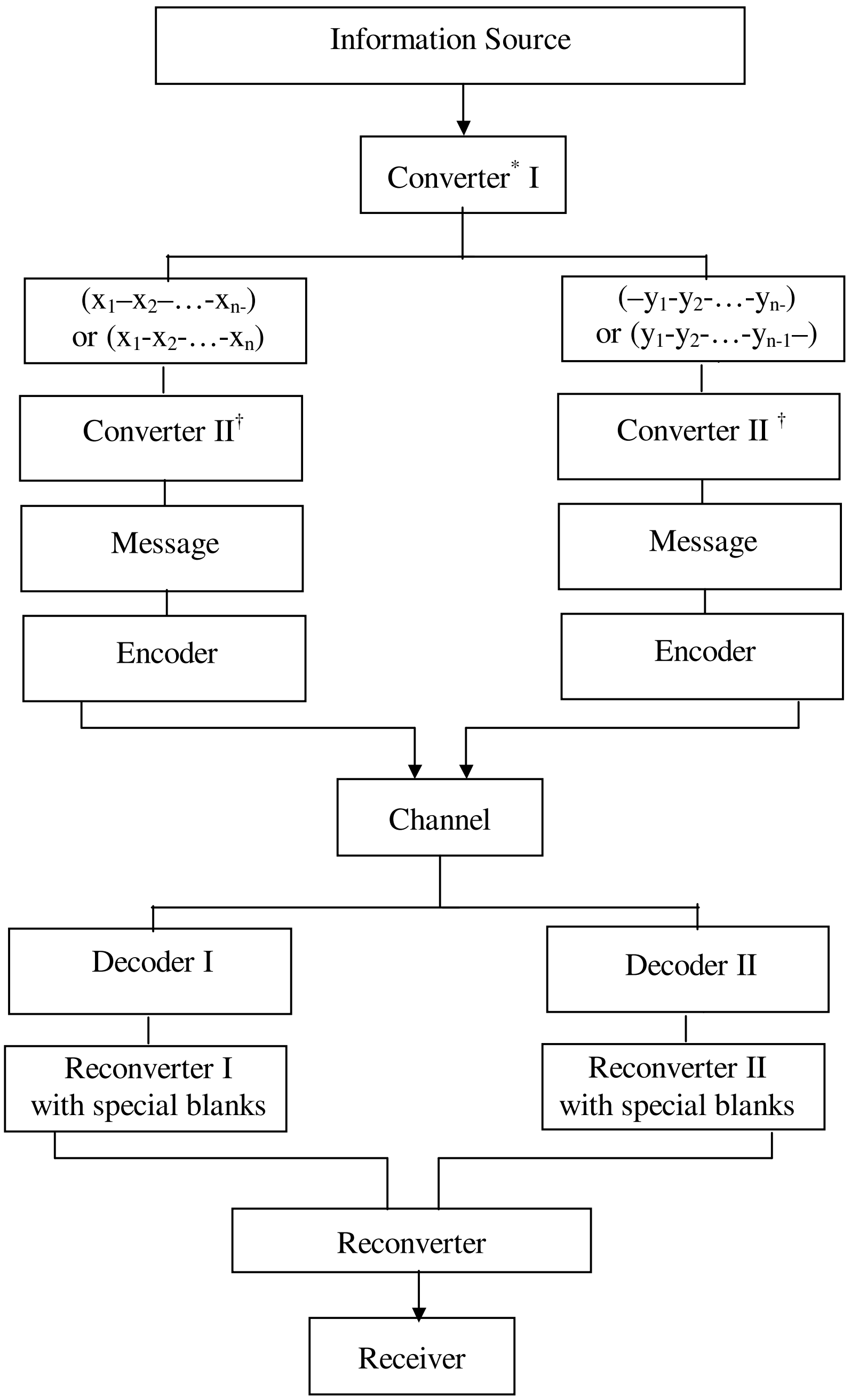

* converter converts this concatenated code into two codes with special blanks.
† the special blanks are omitted and the message are sent by converter I and converter II.

This sort of coding will be used in networking where some secrecy is to be maintained so that processed as two units by two separate systems and the consolidated result is availed by the receiver. Likewise we can define m-concatenated codes. Deconverter converts the code words with special blanks into the concatenated code.

Suppose we have 3 codes say $C_1 = C(5, 2)$, $C_2 = C(6, 3)$ and $C_3 = C(5, 3)$; the 3-concatinated code $C = \{x_1y_1z_1\ x_2y_2z_2\ x_3y_3z_3\ x_4y_4z_4\ x_5y_5z_5\ y_6) \mid (x_1\ x_2\ \ldots\ x_5) \in C_1, (y_1\ y_2\ y_3\ y_4\ y_5\ y_6) \in C_2$ and $(z_1, z_2, \ldots, z_5) \in C_3\}$.

Now this code C of length 16 is sent and converter divides this into 3 codes of length 16 with special blanks codes writer I, writer II and writer III, they convert them into messages of length 5, 6 and 5 respectively and coded and sent via three channels. We can have m number of codes of lengths $n_1, n_2, \ldots, n_m$ with message symbols $k_1, k_2, \ldots, k_m$ respectively. That is $C_i = (n_i, k_i)$, $1 \leq i \leq m$.

Hence C the concatenated m code $C_1, \ldots, C_m$ is given by $C = \{ x_1^1x_1^2 \ldots x_1^m\ x_2^1x_2^2 \ldots x_2^m\ x_3^1x_3^2 \ldots x_3^m \ldots x_1^{p_1}x_2^{p_2} \ldots x_m^{p_m} \}$ where some $p_i$'s may be zero and some $p_j$'s are $n_j$'s; $1 \leq i, j \leq m$.

Now when any x from C is transmitted x is decomposed into m units say 1, 2, …, m where the first unit the code words will be of the form $y = (x_1^1 \ldots x_1^2 \ldots\ x_1^{n_1} \ldots)$, the second unit $y_2 = (-x_2^1\ \ldots -x_2^2\ \ldots\ x_1^{n_2} \ldots)$ and so on. The blanks are special blanks. From the m-units the m-converters will convert these $y_1, y_2, \ldots, y_m$ code words into $z_1 = (x_1^1x_1^2 \ldots x_1^{n_1})$, $z_2 = (x_2^1x_2^2 \ldots x_2^{n_2})$, …, $z_m = (x_m^1x_m^2 \ldots x_m^{n_m})$ by removing the special blanks. Now $z_1, z_2, \ldots, z_m$ are sent and as in case of usual transmission we get the received message say as $P_1, P_2, \ldots, P_m$ then they are passed through the reconverter where they get the special blanks and then they pass through the consolidation channel. Then the receiver receives the message which is a concatenated m-linear code.

We just describe the transmission by the following diagram:

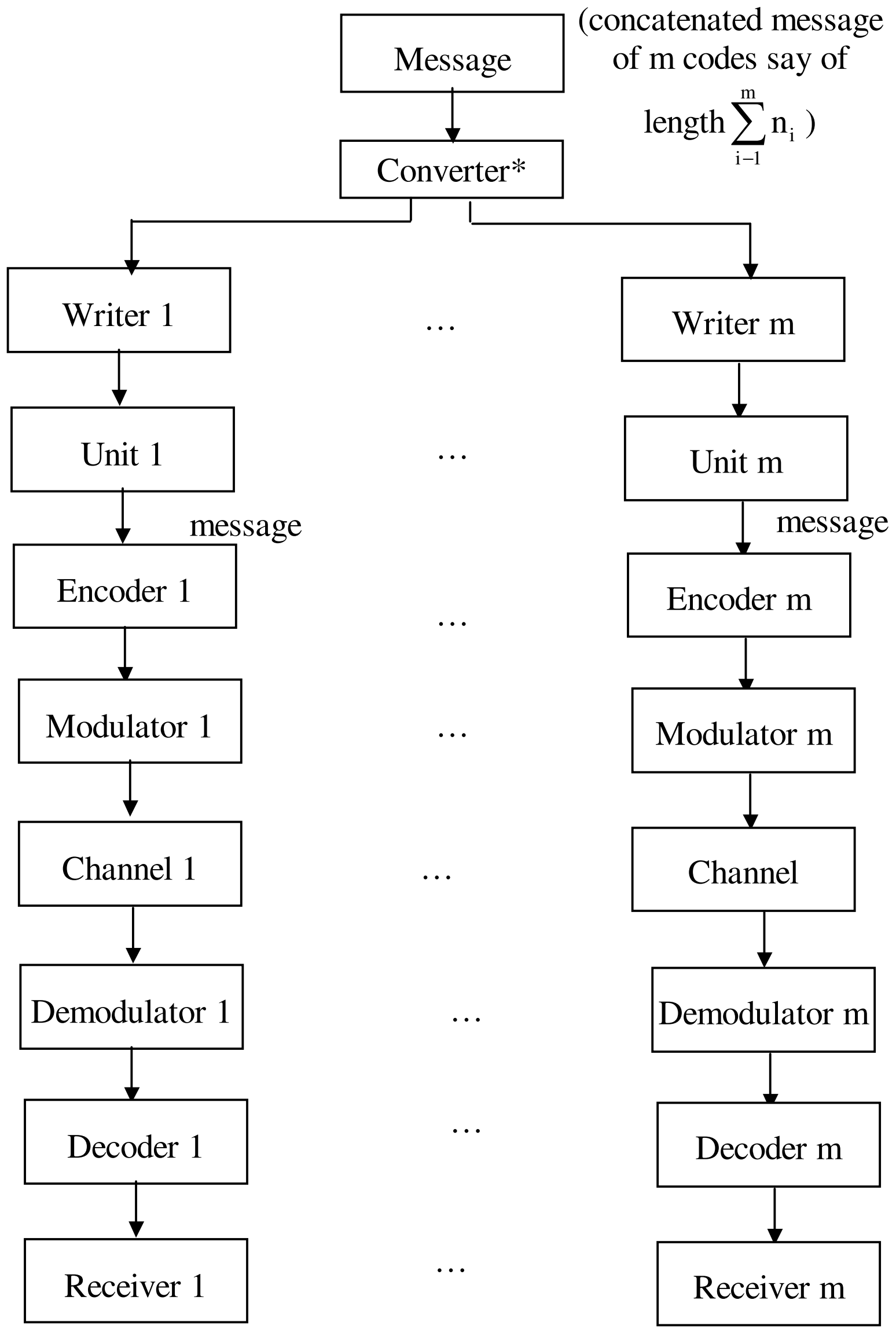

* Converter converts them into m codes with special blanks.

When the concatenated code of the length $\sum_{i=1}^{m} n_i$ with concatenated message symbols $\sum_{i=1}^{m} k_i$ is sent, then the converter divides them into m codes of each length $\sum_{i=1}^{m} n_i = n_1+n_2 + \ldots + n_m$ with special blanks.

Now these are sent to m writers writer 1, writer 2, … writer m. These writers write the codes by deleating the special blanks spaces.

These units acts as the sender of the message, these messages are encoded by encoder 1, encoder 2, …, encoder m. From the encoder it passes thro' the modulator, and from the modulators to the channels from channels to decoders and from the decoders to the m receiver if the purpose of the concatenated coding is to send the message to m receivers from one sender.

But on the other hand if the purpose is to send to receiver but each channel or system does a different type of m-jobs and should reach the single receiver then we have the decoded m-messages would be sent to m-reversing units and then m of these reversing units to m-rewriters who include the blanks and from it reaches a single reconverter and from the reconverter the receiver receives the message then in that case of extended diagram would be as follows.

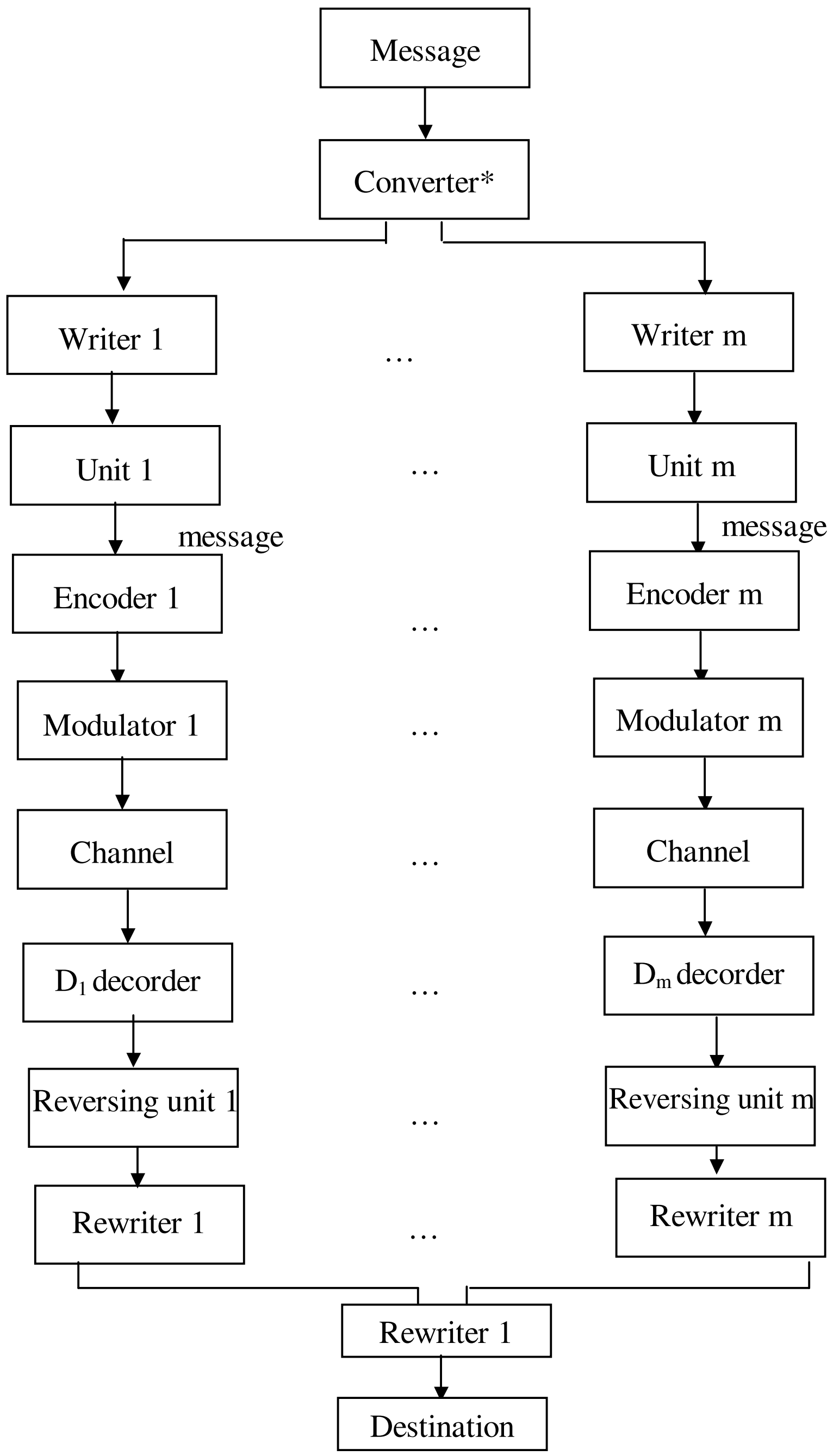

Thus it is a 12 step process.

### 6.2 Concatenation of RD Codes with CR-metric

The concatenated code consists of an outer and an inner code. The outer code is a RD code and the inner code is a binary code and the outer code is a code over the inner code. The concatenated code consists of the codewords of the outer code expressed in terms of the alphabets of the inner code. This concatenated coding system is depicted in the following figure:

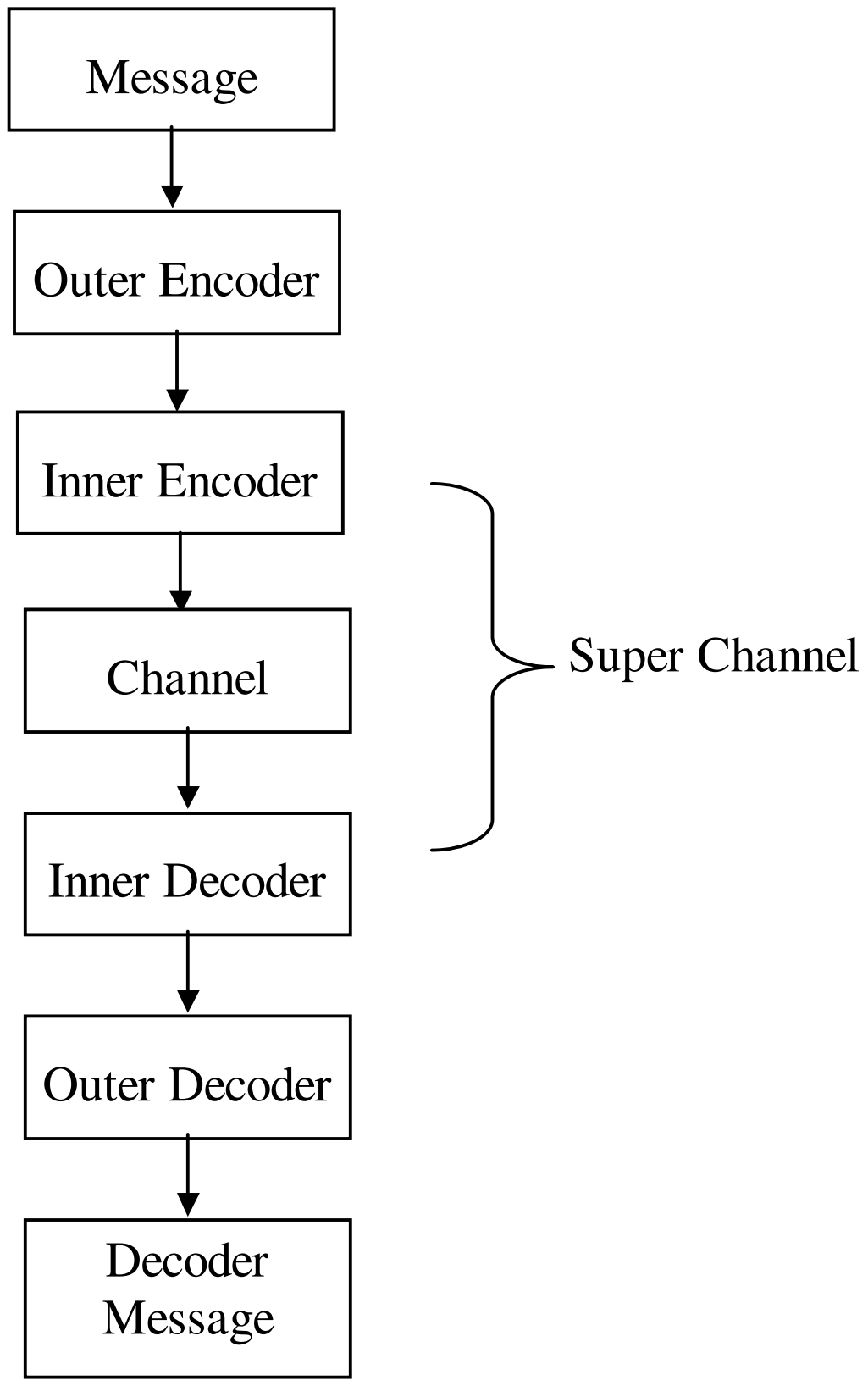

From the above figure we see that the encoder of a concatenated code consists of an outer encoder and an inner encoder corresponding to the outer and inner codes respectively.

Similarly the decoder of a concatenated code consists of an outer decoder and an inner decoder corresponding to the outer and inner codes respectively. Here the combination of the inner encoder, channel and the inner decoder can be thought of as forming a new channel called a super-channel. The super-channel transmits the codewords of the outer code.

The encoding procedure of concatenated codes is given in the following section.

Let the outer code A be a linear $(n_a, k_a, d_a)_{2^{k_b}}$ RD code defined over, the Galois field GF( $2^{k_b}$ ) and the inner code B be a linear $(n_b, k_b, d_b)_2$ code defined over the Galois field GF (2) of order 2. Throughout this section, we consider the inner code with Hamming metric. Let m = $(a_1, a_2, \ldots, a_{k_a}) \in [GF(2^{k_b})]^{k_a}$ be the message to be encoded, where each $a_i \in GF(2^{k_b})$. The procedure of concatenation of the outer code A and the inner code B is given in the following three steps:

**Step 1:** The message m = $(a_1, a_2, \ldots, a_{k_a}) \in$ GF $(2^{k_b})]^{k_a}$ where each $a_i \in GF(2^{k_b})$ is encoded with the outer code A into a codeword. Thus we get a $n_a$-tuple a = $(a_1, a_2, \ldots, a_{n_a})^T$ where each $a_i \in GF(2^{k_b})$.

Now a is the codeword of MRD code which is to be transmitted after applying the following steps:

**Step 2:** $GF(2^{k_b})$ is a $k_b$ - dimensional vector space over the field GF(2). Let g : $GF(2^{k_b}) \to [GF(2)]^{k_b}$ be the mapping of $GF(2^{k_b})$ onto $[GF(2)]^{k_b}$ such that for each symbol $a_i \in GF(2^{k_b})$, i = 1, 2, …, $n_a$ is mapped into a $k_b$-dimensional vector with symbols in GF(2); that is $g(a_i) = (b_{1i}, b_{2i}, \ldots, b_{k_b i})$ where $b_{ji} \in GF(2)$, j = 1, 2, …, $k_b$ and i = 1, 2, …, $n_a$.

Clearly the map g is a one-one linear transformation of the vector space GF( $2^{k_b}$ ) onto $[GF(2)]^{k_b}$.

Hence, GF( $2^{k_b}$ ) = $[GF(2)]^{k_b}$.

After applying steps 1 and 2 to the message m we denote the result in the form of a matrix;

$$\begin{bmatrix} g(a_1) \\ g(a_2) \\ \cdots \\ g(a_{n_a}) \end{bmatrix} = \begin{pmatrix} b_{11} & b_{21} & \cdots & b_{k_b 1} \\ b_{12} & b_{22} & \cdots & b_{k_b 2} \\ \cdots & \cdots & \cdots & \cdots \\ b_{1n_a} & b_{2n_a} & \cdots & b_{k_b n_a} \end{pmatrix}, \qquad (6.2.1)$$

where $b_{ij} \in$ GF (2), i = 1, 2, …, $k_b$ and j = 1, 2, …, $n_a$.

**Step 3:** Each g ($a_i$) = ($b_{1i}$, $b_{2i}$, …, $b_{n_b i}$ ) is encoded with the inner code resulting in a codeword from the inner code given by ($b_{1i}$, $b_{2i}$, …, $b_{n_b i}$ ), where $b_{ji} \in$ GF(2), j = 1, 2, …, $n_b$ and i = 1, 2, …, $n_a$.

After applying steps 1, 2, and 3 of the encoding procedure described above to the message m and encoding each row of the matrix (6.2.1) by using the inner code we get the codeword of the concatenated code represented by the following matrix.

$$\begin{pmatrix} b_{11} & b_{21} & \cdots & b_{n_b 1} \\ b_{12} & b_{22} & \cdots & b_{n_b 2} \\ \cdots & \cdots & \cdots & \vdots \\ b_{1n_a} & b_{2n_a} & \cdots & b_{n_b n_a} \end{pmatrix} \qquad (6.2.2)$$

where $b_{ij} \in$ GF(2), i = 1, 2, …, $n_b$ and j = 1, 2, …, $n_a$.

From here onwards the matrix (6.2.2) will be known as the Concatenated Code Matrix (CC matrix) and will be denoted by CCM (m), that is, it is the related concatenated code matrix of the message m.

**DEFINITION 6.2.1:** *The concatenated code obtained by concatenating the outer RD code, a defined over GF(* $2^{k_b}$ *) and the inner binary code B is the set of all concatenated code matrices. Let CCM =* $\{CCM(m) / m \in [GF(2^{k_b})]^{k_a}\}$ *denote the collection of all* $(n_a \times n_b)$ *CC matrices or equivalently.*

$CCM = \{(b_{ij}): b_{ij} \in GF(2), i = 1, 2, \ldots, n_a, j = 1, 2, \ldots, n_b$ *and* $(b_{ij})\}$ *is a CC matrix.*

**DEFINITION 6.2.2:** *Let CCM be the concatenated code obtained by concatenating the outer RD code A and the inner binary code B. For two concatenated code matrices X, Y* $\in$ *CCM define* $d_c(X, Y) = r(X + Y)$*, where* $r(X+Y)$ *denotes the rank of the matrix X + Y over GF (2) obtained by adding the CC matrices X and Y, using the usual matrix addition modulo 2. That is, if* $X = (a_{ij})$ *and* $Y = (b_{ij})$ *then* $X + Y = (a_{ij} + b_{ij}) \bmod 2$*.*

By the usual properties of the rank of the matrix, for every X, Y, Z $\in$ CCM.

(i) $d_c(X,Y) \geq 0$
(ii) $d_c(X, Y) = 0$ if and only if X = Y,
(iii) $d_c(X, Y) = d_c(Y, X)$
(iv) $d_c(X, Y) \leq d_c(X, Z) + d_c(Z, Y)$.

Thus $d_c$ is a metric on the set of all concatenated code matrices on CCM and we define $d_c$ as the concatenated rank metric (CR-metric).

**DEFINITION 6.2.3:** *The set of all concatenated code matrices (CCM) equipped with the concatenated rank metric (CR-metric)* $d_c$ *is called the concatenated Rank Metric code (CRM code).*

Throughout we denote a CRM code by (CCM, $d_c$). It is clear from the very construction that the CRM code is not a RD code or a binary code.

**DEFINITION 6.2.4:** *Let (CCM, $d_c$) be a CRM code. The minimum distance of the concatenated rank metric code (CCM, $d_c$) is defined as*

$$d = \min \{r(X+Y) / X, Y \in CCM, X \neq Y\}.$$

For X, Y ∈ CCM we have X + Y ∈ CCM we can restate the above definition as follows:

**DEFINITION 6.2.5:** *Let (CCM, $d_c$) be a CRM code with $d_c$ the concatenated rank metric. The minimum distance of the concatenated rank metric code (CCM, $d_c$) is defined as* $d = \min \{r(X) : X \in CCM - \{0\}\}$.

**DEFINITION 6.2.6:** *Let (CCM, $d_c$) be a CRM code with $d_c$ the concatenated rank metric. Let X ∈ CCM. We say r (X) = 0 if and ony if X = 0.*

**DEFINITION 6.2.7:** *Let (CCM, $d_c$) be a CRM code with $d_c$ the concatenated rank metric. Let m be the transmitted message with the corresponding CC matrix viz., X = CCM (m). Let Y be the received matrix. If r (X+Y) = 0, that is, X+Y = 0 which implies X = Y, then Y is the correct message since addition of X and Y is under addition modulo 2.*

*If r (X+Y) > s, s > 0, then Y has an error of rank s.*

A relation between the minimum distance of the concatenated rank metric code and the minimum distance of the outer RD code is given by the following theorem.

**THEOREM 6.2.1:** *Let the outer code A be a linear* $(n_a, k_a, d_a)_{2^{k_b}}$ *RD code defined over GF($2^{k_b}$) and the inner code B be a binary linear ($n_b$, $k_b$, $d_b$)$_2$ code defined over GF(2). Let (CCM, $d_c$) be the concatenated rank metric code with the CR-metric $d_c$* *obtained by using the outer code A and the inner code B. Then the minimum distance of the concatenated rank metric code (CCM, $d_c$) is $d_a$ : $d_a$ is just the minimum distance of the outer RD code.*

***Proof :*** Let the outer RD code A be a $(n_a, k_a, d_a)_{2^{k_b}}$ code with minimum distance $d_a$ and let the inner code B be a binary linear $(n_b, k_b, d_b)_2$ code defined over GF(2). Let (CCM, $d_c$) be the concatenated rank metric code obtained by using the outer code A and the inner code B. Let the message m = $(a_1, a_2, \ldots, a_{k_a}) \in [GF(2^{k_b})]^{k_b}$ where each $a_i \in GF(2^{k_b})$ be encoded with the outer code A into a codeword a = $(a_1, a_2, \ldots, a_{n_a})^T \in$ A. Since the minimum distance of the outer RD code is $d_a$, so $r(a) \geq d_a$. This implies atleast $d_a$ columns of the matrix $a^T$ are linearly independent over GF(2). Without loss of generality let us assume that the first $d_a$ columns of the matrix $a^T$ are linearly independent over GF(2). That is, $a_1, a_2, \ldots, a_{d_a}$ are linearly independent over GF(2). By the method of concatenating the codes A and B given earlier and after applying the mapping g we get the matrix

$$B = \begin{pmatrix} b_{11} & b_{21} & \cdots & b_{k_b 1} \\ b_{12} & b_{22} & \cdots & b_{k_b 2} \\ \cdots & \cdots & \cdots & \cdots \\ b_{1d_a} & b_{2d_a} & \cdots & b_{k_b d_a} \\ \cdots & \cdots & \cdots & \cdots \\ b_{1n_a} & b_{2n_a} & \cdots & b_{k_b n_a} \end{pmatrix}$$

where $b_{ij} \in GF(2)$, i = 1, 2, …, $k_b$ and j = 1, 2, …, $n_a$.

Since $r(a) \geq d_a$, we get $r(b) \geq d_a$. We have, after applying the inner encoder to the matrix B, let the concatenated code matrix corresponding to the message m be CCM(m) = $b_1$, that is,

$$b_1 = \begin{pmatrix} b_{11} & b_{21} & \cdots & b_{n_b 1} \\ b_{12} & b_{22} & \cdots & b_{n_b 2} \\ \cdots & \cdots & \cdots & \cdots \\ b_{1d_a} & b_{2d_a} & \cdots & b_{n_b d_a} \\ \cdots & \cdots & \cdots & \cdots \\ b_{1n_a} & b_{2n_a} & \cdots & b_{nfc_b n_a} \end{pmatrix}$$

where $b_{ij} \in GF(2)$, $i = 1, 2, \ldots, n_b$ and $j = 1, 2, \ldots, n_a$.

Then, $r(b_1) \geq d_a$ since the addition of parity bits by the inner code B to the matrix b does not change the rank of the matrix b. Hence the minimum distance of the concatenated rank metric code of $d_a$.

The above theorem is illustrated by the following example.

***Example 6.2.1:*** Let the outer code A be a (2, 1, 1) RD code with minimum distance 1 defined over $GF(2^2) = \{0, 1, \alpha, \alpha^2\}$ where $\alpha$ is the root of the primitive polynomial $x^2 + x + 1$ of $GF(2^2)$. Let the generator matrix of the outer code A be G = (1 0) where $0, 1 \in GF(2^2)$. Let the inner code B be a binary (4, 2, 3) code with minimum distance 3 be defined over GF(2) having the parity check matrix

$$H = \begin{pmatrix} 1 & 0 & 1 & 0 \\ 1 & 1 & 0 & 1 \end{pmatrix}.$$

The mapping $g : GF(2^2) \to [GF(2)]^2$ is the mapping given by the primitive polynomial of $GF(2^2)$, defined as $0 \to 00$, $1 \to 01$, $\alpha \to 10$, $\alpha^2 \to 11$.

The concatenated rank metric code

$$\text{CCM} = \left\{ \begin{pmatrix} 0 & 0 & 0 & 0 \\ 0 & 0 & 0 & 0 \end{pmatrix}, \begin{pmatrix} 0 & 1 & 0 & 1 \\ 0 & 0 & 0 & 0 \end{pmatrix}, \right.$$

$$\left. \begin{pmatrix} 1 & 0 & 1 & 1 \\ 0 & 0 & 0 & 0 \end{pmatrix}, \begin{pmatrix} 1 & 1 & 1 & 0 \\ 0 & 0 & 0 & 0 \end{pmatrix} \right\}.$$

The CCM has minimum distance d = 1, (since rank of reach of the CC matrices in CCM is either zero or one).

In general, in the case of concatenated codes with Hamming matrix the true minimum distance of the code cannot be obtained but only a lower bound can be obtained. Whereas for the CRM code constructed by us the minimum distance calculated in the above theorem is the true minimum distance and not a lower bound. We prove this in the following corollary.

***Corollary 6.2.1:*** Let the outer code A be a linear $(n_a, k_a, d_a)_{2^{k_b}}$ RD code with minimum distance $d_a$ defined over GF ($2^{k_b}$) and the inner code B be a binary linear $(n_b, k_b, d_b)_2$ code with minimum distance $d_b$ defined GF (2). Let (CCM, $d_c$) be the CRM code obtained by using the codes A and B. The minimum distance of the CRM code is the true minimum distance and it is just the minimum distance of the outer code.

***Proof:*** Let the outer code A be a linear $(n_a, k_a, d_a)_{2^{k_b}}$ RD code with minimum distance $d_a$ defined over GF($2^{k_b}$) and the inner code B be a binary linear $(n_b, k_b, d_b)_2$ code defined over GF(2). Let(CCM, $d_c$) be the concatenated rank metric code obtained by using the outer code A and the inner code B. From theorem 6.2.1 we know that the minimum distance of the concatenated rank metric code (CCM, $d_c$) is $d_a$. Hence the minimum distance of the concatenated rank metric code is the true minimum distance and not a lower bound for the minimum distance.

It is interesting to note that, in the case of (CCM, $d_c$) code or CRM codes if we want to have a desired minimum distance we can without hesitation take the outer RD code with the desired minimum distance and use any convenient binary inner code, for the minimum distance is independent of he choice of the inner code. Thus this concatenation technique helps one to construct any CRM code of desired minimum distance which is not enjoyed by any other class of codes.

**THEOREM 6.2.2:** *(Error detection theorem for CRM codes) Let (CCM, $d_c$) be the concatenated rank metric code with CR-metric $d_c$ obtained by using the codes A and B where the outer code A is a linear $(n_a, k_a, d_a)_{2^{k_b}}$ RD code defined over GF($2^{k_b}$) and the inner code B is a binary linear $(n_b, k_b, d_b)_2$ code defined over GF(2). Let X be the transmitted CC matrix and Y be the received matrix. If $r(Y) < d_a$ then an error has occurred in the matrix Y during transmission.*

***Proof :*** Let m be the message transmitted. Let X be the associated concatenated code matrix of m, that is X = CCM(m). Let Y be the received concatenated code matrix which contains errors occurred during transmission. Suppose $r(Y) < d_a$. From corollary 6.2.1 the minimum distance of the concatenated rank metric code is the true minimum distance and hence $r(X) \geq d_a$ for every $X \in$ CCM. So, $r(Y) < d_a$ indicates that an error has occurred during transmission.

It is very important to note that from the above theorem we can just by computing the rank of the received matrix, we can immediately conclude that an error has occurred during transmission.

We give the decoding procedure for CRM codes.

Our decoding procedure for CRM codes is nothing but to obtain a decoding method for the received concatenated code matrix which is done by first decoding with the inner code and then decoding with the outer code. The reason for doing so is due to the systematic construction of CRM codes explained earlier. We know by the very construction of the CRM codes that the received concatenated code matrix has its entries from the alphabets of the inner code. Hence the received concatenated code matrix is first decoded by using the inner code and then it is decoded with the decoding algorithm of the outer code. Hence the decoder of a concatenated code consists of an inner decoder and an outer decoder. Now we explain the decoding procedure for concatenated rank metric codes (CCM, $d_c$) with CR-matrix $d_c$.

Let Y be the received message. Clearly Y is the matrix given by

$$Y = \begin{pmatrix} y_{11} & y_{21} & \cdots & y_{n_b 1} \\ y_{12} & y_{22} & \cdots & y_{n_b 2} \\ \cdots & \cdots & \cdots & \vdots \\ y_{1n_a} & y_{2n_a} & \cdots & y_{n_b n_a} \end{pmatrix} \qquad (6.2.3)$$

where $y_{ij} \in GF(2)$, $i = 1, 2, \ldots, n_b$ and $j = 1, 2, \ldots, n_a$.

We describe the decoding procedure of a received concatenated code matrix Y in the following three steps.

**Step 1:** For each $i = 1, 2, \ldots, n_b$ decode $(y_{1i}, y_{2i}, \ldots, y_{n_b i})$ using the inner code decoder and obtain a $k_b$ - tuple $(y_{1i}, y_{2i}, \ldots, y_{k_b i})$ where $b_{ji} \in GF(2)$, $j = 1, 2, \ldots, n_b$ and obtain the matrix.

$$y = \begin{pmatrix} y_{11} & y_{21} & \cdots & y_{k_b 1} \\ y_{12} & y_{22} & \cdots & y_{k_b 2} \\ \cdots & \cdots & \cdots & \vdots \\ y_{1n_a} & y_{2n_a} & \cdots & y_{k_b n_a} \end{pmatrix}$$

where $y_{ij} \in GF(2)$, $I = 1, 2, \ldots, k_b$ and $j = 1, 2, \ldots, n_a$. Here when decoding with the inner code the number of columns in the received concatenated matrix Y is reduced and the number of row remains unaltered.

**Step 2:** Applying the mapping $g^{-1}$: $[GF(2)]^{k_b} \to GF(2^{k_b})$ on each row of the matrix y we get $z_i \in GF(2^{k_b})$. Hence corresponding to the matrix y we obtain a $n_a$ tuple $(z_1, z_2, \ldots, z_{n_a}) \in [GF(2^{k_b})]^{n_a}$.

**Step 3:** Now decode the $n_a$-tuple
$(z_1, z_2, \ldots, z_{n_a}) \in [GF(2^{k_b})]^{n_a}$ using the outer decoder. Thus the transmitted message corresponding to the received concatenated code matrix is given by
$(z_1, z_2, \ldots, z_{k_a}) \in [GF(2^{k_b})]^{k_a}$ which gives only the message symbols.

The following examples illustrates the decoding procedure described above.

***Example 6.2.2:*** Let the outer code A be a (2, 1,1) RD code defiend over $GF(2^2) = \{0, 1, \alpha, \alpha^2\}$ where $\alpha$ is the root of the primitive polynomial $x^2 + 2 + 1$ of $GF(2^2)$. Let the generator matrix of the outer code A be G = (1 0) where $0, 1 \in GF(2^2)$. Let the inner code B be a binary (4,2,3) code defined over GF(2) having the parity check matrix

$$H = \begin{pmatrix} 1 & 0 & 1 & 0 \\ 1 & 1 & 0 & 1 \end{pmatrix}.$$

The mapping $g : GF(2^2) \to [GF(2)]^2$ is the mapping given by the primitive polynomial of $GF(2^2)$, defined as $0 \to 00$, $1 \to 01$, $\alpha \to 10$, $\alpha^2 \to 11$.

The concatenated rank metric code

$$CCM = \left\{ \begin{pmatrix} 0 & 0 & 0 & 0 \\ 0 & 0 & 0 & 0 \end{pmatrix}, \begin{pmatrix} 0 & 1 & 0 & 1 \\ 0 & 0 & 0 & 0 \end{pmatrix}, \right.$$

$$\left. \begin{pmatrix} 1 & 0 & 1 & 1 \\ 0 & 0 & 0 & 0 \end{pmatrix}, \begin{pmatrix} 1 & 1 & 1 & 0 \\ 0 & 0 & 0 & 0 \end{pmatrix} \right\}.$$

Suppose the matrix $Y = \begin{pmatrix} 1 & 1 & 1 & 0 \\ 1 & 0 & 0 & 0 \end{pmatrix}$ is received. The received matrix Y is decoded using the decoding procedure of CRM codes described above.

By step 1, we first use the inner decoder to decode the received matrix Y. The inner decoder decodes (1 1 1 0) to (1 1) and (1 0 0 0) to ( 0 0) to obtain the matrix $y = \begin{pmatrix} 1 & 1 \\ 0 & 0 \end{pmatrix}$.

Then by step 2, on applying the map $g^{-1} : [GF(2)]^2 \to GF(2^2)$ to each row of matrix y, we get $y_2 = \begin{pmatrix} \alpha^2 \\ 0 \end{pmatrix}$. By step 3, the outer decodes this column vector $y_2$ to $\alpha^2$ the transmitted message symbol (since in this example the RD code has only one message symbol).

Now from these concatenated MRD codes we can define biconcatenated MRD codes or concatenated MRD bicodes, $C = C_1 \cup C_2$ where each $C_i$ is a concatenated CRM codes. We can also extend this to n-concatenated CRM codes or concatenated n-CRM codes as $C = C_1 \cup C_2 \cup \ldots \cup C_n$ where each $C_i$ is a CRM code; $1 \leq i \leq n$. This form of n-codes will help in the simultaneous processing using n-units which will be time saving.

We can also define quasi n-concatenated CRM codes, when some t of them are concatenated CRM codes and the rest are MRD codes. Further the n-concatenated special super codes; $C = C_1 \cup C_2 \cup \ldots \cup C_n$ if each $C_i$ is a special super code; i = 1, 2, …, n. If n = 2 we get the biconcatenated special super code.

We can also have biconcatenated CRM code with a special super codes defined as mixed biconcatenated code. That is $C = C_1 \cup C_2$ where $C_1$ is a super special code and $C_2$ is a CRM code. We can extend this to mixed n-concatenated codes where t of them are super special codes ($t < n$) and rest are just CRM codes. We can also have the notion of quasi mixed n-concatenated codes where $t_1$ of them are CRM codes, $t_2$ of them are just MRD codes, $t_3$ of them special super codes and the rest of them are just linear codes. Thus we can have quasi mixed n-concatenated codes. These codes will help when bulk messages are sent and this will save time and economy.

# FURTHER READINGS

# INDEX

**T**

**V**

# ABOUT THE AUTHORS

**Dr.W.B.Vasantha Kandasamy** is an Associate Professor in the Department of Mathematics, Indian Institute of Technology Madras, Chennai. In the past decade she has guided 13 Ph.D. scholars in the different fields of non-associative algebras, algebraic coding theory, transportation theory, fuzzy groups, and applications of fuzzy theory of the problems faced in chemical industries and cement industries. She has to her credit 646 research papers. She has guided over 68 M.Sc. and M.Tech. projects. She has worked in collaboration projects with the Indian Space Research Organization and with the Tamil Nadu State AIDS Control Society. She is presently working on a research project funded by the Board of Research in Nuclear Sciences, Government of India. This is her 66th book.

On India's 60th Independence Day, Dr.Vasantha was conferred the Kalpana Chawla Award for Courage and Daring Enterprise by the State Government of Tamil Nadu in recognition of her sustained fight for social justice in the Indian Institute of Technology (IIT) Madras and for her contribution to mathematics. The award, instituted in the memory of Indian-American astronaut Kalpana Chawla who died aboard Space Shuttle Columbia, carried a cash prize of five lakh rupees (the highest prize-money for any Indian award) and a gold medal.
She can be contacted at vasanthakandasamy@gmail.com
Web Site: http://mat.iitm.ac.in/home/wbv/public_html/
or http://www.vasantha.in

---

**Dr. Florentin Smarandache** is a Professor of Mathematics at the University of New Mexico in USA. He published over 75 books and 200 articles and notes in mathematics, physics, philosophy, psychology, rebus, literature. In mathematics his research is in number theory, non-Euclidean geometry, synthetic geometry, algebraic structures, statistics, neutrosophic logic and set (generalizations of fuzzy logic and set respectively), neutrosophic probability (generalization of classical and imprecise probability). Also, small contributions to nuclear and particle physics, information fusion, neutrosophy (a generalization of dialectics), law of sensations and stimuli, etc. He got the 2010 Telesio-Galilei Academy of Science Gold Medal, Adjunct Professor (equivalent to Doctor Honoris Causa) of Beijing Jiaotong University in 2011, and 2011 Romanian Academy Award for Technical Science (the highest in the country). Dr. W. B. Vasantha Kandasamy and Dr. Florentin Smarandache got the 2011 New Mexico Book Award for Algebraic Structures. He can be contacted at smarand@unm.edu

---

**Dr. R.Sujatha** has obtained her doctorate degree on Algebraic Coding Theory from Indian Institute of Technology under the guidance of Dr. W.B.Vasantha, Department of mathematics, IIT (madras) Chennai – 36. She is presently working in SSN college of engineering. She has 14 years of teaching and research experience. She has published over 35 research articles and has authored two books.She can be contacted at sujathar@ssn.edu.in

---

**Dr. R.S.Raja Durai** is an Assistant Professor of Jaypee University of Information Technology, has obtained his doctorate from the Indian Institute of Technology, Chennai – 36, under the guidance of Dr.W.B.Vasantha Kandasamy. He has published over 10 research papers and is currently guiding Ph.D students.